\documentclass{article}%
\usepackage{amsfonts}
\usepackage{amsmath}
\usepackage{amssymb}
\usepackage{graphicx}%
\providecommand{\U}[1]{\protect\rule{.1in}{.1in}}
\newtheorem{theorem}{Theorem}[section]

\newtheorem{corollary}[theorem]{Corollary}

\newtheorem{definition}[theorem]{Definition}
\newtheorem{example}[theorem]{Example}

\newtheorem{lemma}[theorem]{Lemma}

\newtheorem{proposition}[theorem]{Proposition}
\newtheorem{remark}[theorem]{Remark}

\newenvironment{proof}[1][Proof]{\noindent\textbf{#1.} }{\ \rule{0.5em}{0.5em}}
\begin{document}

\title{Rough Paths on Manifolds }
\author{Thomas Cass\thanks{Research supported by EPSRC grant EP/F029578/1}
\footnotemark[3] , Christian Litterer\thanks{Research supported by a grant of
the Leverhulme trust (F/08772/E).} \thanks{Part of this research was carried
out during a stay at the Hausdorff Institute, Bonn} and Terry Lyons\\Mathematical Institute and \\Oxford-Man Institute of Quantitative Finance \\University of Oxford \\24-29 St Giles', Oxford, OX1 3LB, UK }
\date{\today }
\maketitle

\begin{abstract}
We develop a fundamental framework for and extend the theory of rough paths to
Lipschitz-$\gamma$ manifolds.

\end{abstract}

\section{Introduction}

The theory of rough paths is an extension of classical Newtonian calculus
aimed at allowing models for the interaction of highly oscillatory and
potentially non-differentiable systems. These models take the form of
differential equations driven by rough paths.

An example might be that one seeks the trajectory
\[
dy_{t}=\sum_{i=1}^{N}f^{i}\left(  y_{t}\right)  dx_{t}^{i},
\]
where $x_{t}=\left(  x_{t}^{i}\right)  _{i=1}^{N}$is a highly oscillatory
path. The thrust of the theory is that one can indeed identify a family of
paths (known as rough paths) for which such equations make good sense and
which generalise the classical notion of a differentiable path to a sufficient
extent that it allows $x_{t}$ to be chosen randomly. Rough paths are rich
enough to include many stochastic processes that are nowhere differentiable
(being a bit more technical semimartingale paths almost surely extend to rough
paths, but there are many other examples).

The theory introduced in Lyons \cite{lyons2}, \cite{lyons3},
\cite{lyons-revista} has, building strongly on ideas of Chen \cite{chen},
\cite{chen2}, Magnus \cite{magnus}, Fliess \cite{fliess} as well as analytic
estimates of Young \cite{young}, produced a robust class of objects known as
$p$ -rough paths. An introduction to the basic theorems can be found in Lyons,
Caruana, Levy \cite{lyons-stflour} and a recent comprehensive account in the
finite dimensional case can be found in Friz, Victoir \cite{friz}.

However, the theory to date does not properly identify the fully geometric
framework of rough paths evolving on manifolds - even though many of the
important controlled differential equations are of the type%
\begin{equation}
dy_{t}=f\left(  y_{t}\right)  dx_{t} \label{first}%
\end{equation}
and occur in a geometric context (e.g. parallel translation and development).
One of the breakthroughs in probability on differential manifolds was the
description of Brownian motion as a projection of the solution to a globally
defined differential equation on the orthogonal frame bundle of the kind
$\left(  \ref{first}\right)  ,$ see Elworthy \cite{elworthy}$.$

The roughness of a path $X$ is controlled by a real number $p\in
\lbrack1,\infty).$ $1$ -rough paths correspond to paths of bounded variation.
Semimartingale paths are almost surely $p$ -rough paths for all $p>2.$ It is
well known (see e.g. Lyons \cite{lyons-revista}) that if $X$ is a $p$ -rough
path on a Banach space $V$ and~$f$ is a Lip-$\left(  p+\varepsilon\right)  $
function defined around the trajectory of $X$ and taking values in a Banach
space $W,$ then $f\left(  X\right)  $ is a $p-$ rough path on $W.$ Some
readers may wonder what it means for a function to be Lip-\thinspace$p$ for
$p>1,$ as the naive definition would lead one to the view that such functions
are always constant. The theory of rough paths relies on a definition
introduced in Stein \cite{stein}, where the author also proves the appropriate
Whitney extension theorem.

A crucial feature of rough path theory is that it relies on uniform (as
opposed to asymptotic) estimates. It is these estimates that allow much of the
theory to be moved to the infinite dimensional Banach context. To make a sound
geometric framework for discussing rough paths on manifolds we introduce the
class of Lipschitz-$\gamma$ manifolds with constant $L.$ Our definition is
distinctive - because it contains some global rigidity. A single smooth
non-compact manifold can admit different Lip-$\gamma$ atlases and different
classes of Lip-$\gamma$ functions and one does not expect maximal atlases.

Given a Lipschitz manifold $M,$ it is immediate that one can place a norm on
the space of currents on $M.$ It is straightforward to develop a theory of
rough paths in the space of one-currents using the results of Lyons
\cite{lyons-revista} in a routine way. However, in general such objects have
no natural locality. They are of interest in their own right, for example
there are rough current based systems that cannot be detected locally although
their effect can be felt globally.

If $\alpha$ is a Lip$\left(  p-1+\varepsilon\right)  $ smooth one form on a
Banach space $V$ and $X$ is a based $p-$rough path on $V,$ then classical
theory shows us that there is no difficulty in integrating $\alpha$ against
the path for any $t.$ This produces, in a natural way, a map
\[
\theta_{t}:X\rightarrow\int_{0}^{t}\alpha\left(  X_{s}\right)  dX_{s}.
\]
The map $\theta$ corresponds to a current valued path, which with further
study can be interpreted as a $p-$ rough path in the currents.

This paper considers objects on a Lip-$\gamma$ manifold $M$ which have a
formal integral against vector valued one forms yielding rough paths as their
integrals. The challenge in the paper is to identify axiomatic conditions that
pick from the general class of current valued rough paths (or even more
generally objects that integrate one forms) the ones that localise to $p$
-rough paths on $M$ and it turns out that doing this precisely is a little delicate.

We also establish existence for the solution of differential equations on
Lipschitz manifolds driven by rough paths. We characterise the differential
equations in terms of generalised connections. It is clear from general
results that (although we do not explore it in the current paper) the solution
process, which we call the It\^{o} map, is uniformly continuous in appropriate contexts.

\section{Background: Rough Paths on Banach spaces}

In this section we give a brief overview of some of the core results of the
theory of rough paths on Banach spaces as it was developed in Lyons
\cite{lyons-revista}. To distinguish them from our new definitions on the
manifold we will where necessary to avoid confusion refer to rough paths on
Banach spaces as Banach space valued or "classical" rough paths. We focus on
those results that are particularly relevant to the development of the theory
on manifolds in the following sections. For a detailed introduction to rough
paths on Banach spaces and proofs of the results in this section see Lyons,
Caruana, Levy \cite{lyons-stflour}. Let $V$ be a Banach space and $T\left(
\left(  V\right)  \right)  $ be the space of formal series of tensors of $V.$
Then $T\left(  \left(  V\right)  \right)  $ is a real non-commutative unital
algebra. We may consider the quotient $T^{(n)}\left(  V\right)  $ of $T\left(
\left(  V\right)  \right)  $ by the ideal $\bigoplus_{i=n+1}^{\infty
}V^{\otimes i}$ and identify it with the space%

\[
T^{\left(  n\right)  }\left(  V\right)  =\bigoplus_{i=0}^{n}V^{\otimes i}.
\]
Let $\left\Vert \cdot\right\Vert $ be the norm on $V$ and let $V^{\otimes i} $
be equipped with the projective tensor norm also denoted by $\left\Vert
\cdot\right\Vert $. We refer to \cite{lyons-stflour} for the requirements a
general tensor norm has to satisfy to be compatible. In the following let
$\Delta_{T}=\left\{  \left(  s,t\right)  \in\left[  0,T\right]
:s\,<t\right\}  .$

\begin{definition}
Let $n\geq1$ be an integer and $X_{\cdot,\cdot}:\Delta_{T}\rightarrow
T^{(n)}(V)$ a continuous map. Then $X$ is a multiplicative functional of
degree $p$ if for all $(s,t)\in\Delta_{T}$ we have $X_{s,t}^{0}=1$ and%
\[
X_{s,u}\otimes X_{u,t}=X_{s,t}\text{ for all }0\leq s\leq u\leq t\leq T.
\]

\end{definition}

Multiplicative functionals of degree $n$ are fully characterised by
specifying
\[
X_{s,t}=\left(  1,X_{s,t}^{1},X_{s,t}^{2},\ldots,X_{s,t}^{n}\right)  \in
T^{\left(  n\right)  }\left(  V\right)  ,
\]
for all $(s,t)\in\Delta_{T}.$ We will refer to $X_{s,t}^{1}$ in the following
also as level one or the trace of a multiplicative functional or rough path.
Particularly important examples of multiplicative functionals are the lifts of
continuous bounded variations paths by the signature map. Let $\phi\in
C_{bv}([0,T],V)$; $s,t\in\Delta_{T}.$ We define the signature $S_{s,t}%
:C_{bv}([0,T],R^{d})\rightarrow T(\left(  V\right)  )$ by the relation%
\[
S_{s,t}(\phi)=\sum_{k=0}^{\infty}\int_{s<t_{1}<\cdots<t_{k}<t}d\phi
(t_{1})\otimes\cdots\otimes d\phi(t_{k}).
\]
Chen's theorem tells us that the signature and its truncations at level $n$
are multiplicative, i.e. satisfy for $0\leq s\leq r\leq t\leq T$%
\[
S_{s,r}(\phi)\otimes S_{r,t}(\phi)=S_{s,t}(\phi).
\]

\begin{definition}
A control $\omega$ is a non-negative continuous function on $\Delta_{T}$ that
is super-additive, i.e. satisfies%
\[
\omega(s,u)+\omega(u,t)\leq\omega\left(  s,t\right)  ,\text{ }0\leq s\leq
u\leq t\leq T
\]
and vanishes on the diagonal.
\end{definition}

Intuitively a control is designed to abstract the important properties of the
notion of $p$-variation. In particular, if $y_{t}$ is a continuous path of
bounded $p-$variation the $p$th power of its \thinspace$p-$ variation%
\[
\sup_{\mathcal{D\in}\left[  0,T\right]  }\sum_{\mathcal{D}}\left\Vert
y_{t_{i+1}}-y_{t_{i}}\right\Vert ^{p}%
\]
is a control.

\begin{definition}
We say a multiplicative functional of degree $n$ has finite $p\geq1$ variation
if%
\begin{equation}
\left\Vert X_{s,t}^{i}\right\Vert \leq\frac{\omega(s,t)^{i/p}}{\beta\left(
\frac{i}{p}\right)  !},\text{ }i=1,\ldots,n,~\left(  s,t\right)  \in
\Delta_{T,} \label{finite-variation}%
\end{equation}
where
\[
\beta=p\left(  1+\sum_{r=3}^{\infty}\left(  \frac{2}{r-2}\right)
^{\frac{\lfloor p\rfloor+1}{p}}\right)  .
\]

\end{definition}

\begin{remark}
The $p$ rather than the more typical $p^{2}$ in the definition of $\beta$ is a
recent innovation due to Hara and Hino \cite{hino}.
\end{remark}

{\ Multiplicative functionals of any degree with finite }$p-$variation are
determined by their values on $T^{\left(  \lfloor p\rfloor\right)  }\left(
V\right)  .$

\begin{theorem}
[Extension Theorem]Let $p\geq1$ and $X$ a multiplicative functional of degree
$\lfloor p\rfloor$ with finite $p-$ variation controlled by some control
$\omega$. For any $m\geq\lfloor p\rfloor+1$ the multiplicative functional $X$
has a unique extension to a multiplicative functional of degree $m$ that is
controlled by $\omega.$
\end{theorem}

In particular the extension theorem implies that the lift of bounded variation
paths to a multiplicative functional on $T\left(  \left(  V\right)  \right)  $
is unique. With the extension theorem in mind we are ready to introduce the
notion of rough paths.

\begin{definition}
A $p$- rough path is a multiplicative functional of degree $\lfloor p\rfloor$
with finite $p$ -variation.
\end{definition}

Let $\Omega_{p}\left(  V\right)  $ denote the space of all rough paths on $V.$
We can define a metric on $p-$rough paths as follows.

\begin{definition}
We define the $p$-variation distance on rough paths by%
\[
d_{p}(X,Y)=\max_{1\leq i\leq\lfloor p\rfloor}\sup_{D\subset\left[  0,T\right]
}\left(  \sum_{D}\left\Vert X_{t_{l-1,}t_{l}}^{i}-Y_{t_{l-1,}t_{l}}%
^{i}\right\Vert ^{p/i}\right)  ^{i/p}.
\]

\end{definition}

It can be shown that for any $p\geq1$ the rough path spaces $\Omega_{p}\left(
V\right)  $ together with the $p-$variation distance $d_{p}$ form complete
metric spaces. and we refer to the topology induced by $d_{p}$ on $\Omega
_{p}\left(  V\right)  $ as the $p$ - variation topology. A particularly
important class of rough paths are the so called geometric rough paths.

\begin{definition}
The geometric $p$-rough paths are the closure of the $1$-rough paths in the
$p$-variation metric. Let $G\Omega_{p}\left(  V\right)  $ denote the geometric
$p-$rough paths on $V.$
\end{definition}

Just as $\left(  \Omega_{p}\left(  V\right)  ,d_{p}\right)  $ the space
$\left(  G\Omega_{p}\left(  V\right)  ,d_{p}\right)  $ is a complete metric
space. The geometric $p-$rough paths are the closure of the (lifts of) bounded
variation paths. Note that by the extension theorem the lift of a bounded
variation path to a $1$ -rough path provided by the signature is unique. It
turns out that nearly all example of practical relevance are geometric rough
paths and for the remainder of this paper we will exclusively work with this
class of paths. A cornerstone of the theory of rough is the development of a
rough integral of $p$ -rough paths $X$ against sufficiently regular Banach
space valued one forms. A crucial role is played by Lipschitz functions
defined by Stein \cite{stein} on subsets of a Banach space.

\begin{definition}
\label{lipschitz functions}Let $V$,$W$ be Banach spaces and $f:F\rightarrow W
$ a function from a Borel set $F\subseteq V$ to $W.$ Let $\gamma>0$ and $k $
be the unique integer such that $\gamma\in(k,k+1].$ For $j=0,\ldots,k$ let
$f^{j}:F\rightarrow L(V^{\otimes j},W).$ The collection $(f^{j})$ is in
Lip-$\gamma$($F,W$ ) if there is a constant $L$ such that
\[
\sup_{x\in F}%
%TCIMACRO{\U{a6}}%
%BeginExpansion
\vert
%EndExpansion
f^{j}(x)%
%TCIMACRO{\U{a6}}%
%BeginExpansion
\vert
%EndExpansion
\leq L
\]
and there exist functions $R_{j}:V\times V\rightarrow L(V^{\otimes j},W)$ such
that for all $x,y$ $\in F,$ $v\in V^{\otimes j}$%
\[
f^{j}(y)(v)=\sum_{l=0}^{k-j}\frac{1}{l!}f^{j+l}(x)(v\otimes(y-x)^{\otimes
l})+R_{j}(x,y)(v)
\]
and
\[
\left\vert R_{j}(x,y)\right\vert \leq L%
%TCIMACRO{\U{a6}}%
%BeginExpansion
\vert
%EndExpansion
x-y%
%TCIMACRO{\U{a6}}%
%BeginExpansion
\vert
%EndExpansion
^{\gamma-j}.
\]

\end{definition}

The notion of Lip-$\gamma$ intuitively captures how well a function is
approximated by polynomials. The multi-linear forms $f^{j}$ in Definition
\ref{lipschitz functions} play a role analogous of derivatives but are in
general (for example if $F$ contains isolated points) not unique. If $F$ is
open the $f^{j}$ are uniquely determined by the derivatives of $f.$ Stein
proves the following important extension theorem if $W=%
%TCIMACRO{\U{211d} }%
%BeginExpansion
\mathbb{R}
%EndExpansion
^{d}.$

\begin{theorem}
\label{whitney-extension}Let $\gamma\geq1,d$ and $k$ positive integers and
$F\subseteq%
%TCIMACRO{\U{211d} }%
%BeginExpansion
\mathbb{R}
%EndExpansion
^{d}$ a closed set. Then there exists a bounded linear extension operator%
\[
E_{\gamma}:Lip-\gamma\left(  F,%
%TCIMACRO{\U{211d} }%
%BeginExpansion
\mathbb{R}
%EndExpansion
^{k}\right)  \rightarrow Lip-\gamma\left(
%TCIMACRO{\U{211d} }%
%BeginExpansion
\mathbb{R}
%EndExpansion
^{d},%
%TCIMACRO{\U{211d} }%
%BeginExpansion
\mathbb{R}
%EndExpansion
^{k}\right)  .
\]
Moreover the norm of $E_{\gamma}$ is independent of the closed set $F.$
\end{theorem}

Note however that the extension operator and its norm does depend on
$\gamma\geq1.$ It is elementary to show that Lip-$\gamma$ functions defined on
an arbitrary Borel set $F$ extend to Lip-$\gamma$ functions on the closure
$\overline{F}$ without changing the Lipschitz constant$.$

\begin{lemma}
Let $\gamma\geq1,V$ and $W$ be Banach spaces, and $F\subseteq V$ a Borel set.
Then any $f\in$Lip-$\gamma\left(  F,%
%TCIMACRO{\U{211d} }%
%BeginExpansion
\mathbb{R}
%EndExpansion
^{k}\right)  $ extends to a function $\widetilde{f}\in$Lip-$\gamma\left(
\overline{F},%
%TCIMACRO{\U{211d} }%
%BeginExpansion
\mathbb{R}
%EndExpansion
^{k}\right)  $ such that
\[
\left\Vert f\right\Vert _{\text{Lip}-\gamma}=\left\Vert \widetilde{f}%
\right\Vert _{\text{Lip-}\gamma}.
\]

\end{lemma}

\begin{proof}
For any two points $x,y\in F$ we have
\begin{align*}
&  \left\vert f^{j}(x)(v)-f^{j}\left(  y\right)  (v)\right\vert \\
&  \leq\left\vert \sum_{l=1}^{k-j}\frac{1}{l!}f^{j+l}(x)(v\otimes
(y-x)^{\otimes l})\right\vert +\left\vert \sum_{l=1}^{k-j}\frac{1}{l!}%
f^{j+l}(y)(v\otimes(x-y)^{\otimes l})\right\vert \\
&  +\left\vert R_{j}(x,y)(v)\right\vert +\left\vert R_{j}(y,x)(v)\right\vert ,
\end{align*}
which implies
\begin{equation}
\left\vert f^{j}(x)(v)-f^{j}\left(  y\right)  (v)\right\vert \leq4L\max\left(
\left\vert x-y\right\vert ,\left\vert x-y\right\vert ^{\gamma-j}\right)
\label{Cauchy}%
\end{equation}
if $\left\vert x-y\right\vert \leq1.$ Given any $x\in\overline{F}$ there
exists a Cauchy sequence $x_{n}\in F$ converging to $x.$ By inequality
$\left(  \ref{Cauchy}\right)  $ for any sequence $y_{n}\in F$ converging to
$x$ the sequence $f\left(  y_{n}\right)  $ is Cauchy and converges to a limit.
Moreover this limit is independent of the particular sequence $y_{n}$ and we
may define $f^{j}\left(  x\right)  $ to be $\lim_{n\rightarrow\infty}%
f^{j}\left(  x_{n}\right)  .$ It is straightforward to check that the thus
defined extension is Lipschitz on $\overline{F}.$
\end{proof}

Another useful property Lip-$\gamma$ functions possess is that they are well
behaved under composition; namely, the composition of two Lip-$\gamma$
functions remains a Lip-$\gamma$ function and the norm of the composition is
easily related to the norm of the composed functions. More precisely, we have
the following result.

\begin{lemma}
\label{composition of lip functions}Let $U,$ $V$ and $W$ be Banach spaces and
$E$ be some subset of $U.$ Suppose that $f:E\rightarrow V$ and $g:V\rightarrow
W$ are Lip-$\gamma$ for $\gamma\geq1$ then the composition $g\circ
f:E\rightarrow W$ is Lip-$\gamma$ and moreover
\[
\left\vert \left\vert g\circ f\right\vert \right\vert _{\text{Lip-}\gamma}\leq
C\left(  \gamma\right)  \left\vert \left\vert g\right\vert \right\vert
_{\text{Lip-}\gamma}\max\left(  \left\vert \left\vert f\right\vert \right\vert
_{\text{Lip-}\gamma}^{k},1\right)  ,
\]
where $k$ is unique integer such that $\gamma\in(k,k+1]$ $\ $and $C\left(
\gamma\right)  $ does not depend on $f$ or $g.$
\end{lemma}

\begin{proof}
The chain rule (generalised to higher derivatives) provides candidates for the
auxiliary functions $\left(  g\circ f\right)  ^{i}$ for $i=0,1,...,k$.
\ Defining the remainder functions using these candidates and using the fact
that $f$ and $g$ are Lip$-\gamma$ it is straight-forward to verify that
$g\circ f$ is Lip$-\gamma$ together with the stated estimate on $\left\vert
\left\vert g\circ f\right\vert \right\vert _{\text{Lip-}\gamma}.$
\end{proof}

Lip-$\gamma$ is not a local property it turns out however under certain
conditions it can be deduce from uniform local estimates. We proceed in two steps.

\begin{lemma}
\label{convex domains}Let $\gamma\geq1$, $k$ be the unique integer such that
$\gamma\in(k,k+1]$ and suppose $U$ is a convex open subset of a finite
dimensional normed vector space. Then a Banach space valued function $f$ is
Lip-$\gamma$ on $U\ $with $\left\Vert f\right\Vert _{Lip-\gamma}\leq C$ if and
only if $f$ has $k$ derivatives bounded by $C$ on $U$ and the $k^{th}$
derivative is $\left(  \gamma-k\right)  -$ H\"{o}lder with constant at most
$C$, i.e. it satisfies%
\[
\sup_{x,y\in U}\frac{\left\vert f^{k}\left(  x\right)  -f^{k}\left(  y\right)
\right\vert }{\left\vert x-y\right\vert ^{\gamma-k}}\leq C.
\]

\end{lemma}

\begin{proof}
Suppose $f$ is $Lip-\gamma,$ then it is self-evident from the convexity of $U$
and the definition of $Lip-\gamma$ that the auxiliary functions $f^{1}%
,...,f^{k}$ are precisely the first $k$ derivatives of $f$ and that they are
bounded by $\left\vert \left\vert f\right\vert \right\vert _{\text{Lip-}%
\gamma}.$ The H\"{o}lder continuity follows from this observation too together
with the fact that
\[
\left\vert f^{k}(y)-f^{k}(x)\right\vert \leq\left\vert R_{j}(x,y)\right\vert
\leq C%
%TCIMACRO{\U{a6}}%
%BeginExpansion
\vert
%EndExpansion
x-y%
%TCIMACRO{\U{a6}}%
%BeginExpansion
\vert
%EndExpansion
^{\gamma-k}.
\]
Conversely, if $f$ has $k$ derivatives bounded by $C$ and the $k^{th}$
derivative is $\left(  \gamma-k\right)  -$ H\"{o}lder with constant $C$ then a
straight-forward application of Taylor's Theorem together with elementary
remainder estimates is sufficient to see that $f$ is $Lip-\gamma$ with
$\left\vert \left\vert f\right\vert \right\vert _{\text{Lip-}\gamma}\leq C$
\end{proof}

\begin{lemma}
\label{local to global}Let $\delta>0$ and $U\ $be a convex open subset of a
finite dimensional normed vector space. Let $f$ be a Banach space valued
function on $U$ such that $f|_{B\left(  x,\delta\right)  \cap U}$ is
Lip-$\gamma$ with constant at most $C$ for all $x\in U.$ Then $f$ is
Lip-$\gamma$ on $U$ and $\left\Vert f\right\Vert _{Lip-\gamma}\leq\max\left(
C,\frac{2C}{\delta^{\gamma-\lfloor\gamma\rfloor}}\right)  .$
\end{lemma}

\begin{proof}
Let $k$ be the unique integer such that $\gamma\in(k,k+1].$ First observe that
the intersection $B\left(  x,\delta\right)  \cap U$ is convex ($B\left(
x,\delta\right)  $ and $U$ themselves being convex). It is immediate from the
assumption that $f|_{B\left(  x,\delta\right)  \cap U}$ is $Lip-\gamma$ with
Lipschitz constant at most $C$ and that $f$ has $k$ derivatives bounded by $C$
on $U.$ It remains to check that the $k^{th}$ derivative $f^{k}$ is H\"{o}lder
$\gamma-k,$i.e.%
\[
\sup_{x,y\in U}\frac{\left\vert f^{k}\left(  x\right)  -f^{k}\left(  y\right)
\right\vert }{\left\vert x-y\right\vert ^{\gamma-k}}\leq\max\left(
C,\frac{2C}{\delta^{\gamma-k}}\right)  .
\]
The bound is immediate if $y\in B\left(  x,\delta\right)  \cap U$ as $f$ is by
assumption Lip-$\gamma$ on $B\left(  x,\delta\right)  \cap U.$ If $y\notin
B\left(  x,\delta\right)  \cap U$ we have $\left\vert x-y\right\vert
^{\gamma-k}\geq\delta^{\gamma-k},$ which together with $f^{k}$ being bounded
by $C$ implies the claim.
\end{proof}

The concept of an almost $p-$rough path plays an important role in the
definition of the integral of a rough path against sufficiently regular one forms

\begin{definition}
Let $p>1$ and $\omega$ a control over $\left[  0,T\right]  .$ A continuous map
$X:\Delta_{T}\rightarrow T^{(\lfloor p\rfloor)}(W)$ is an almost $p-$ rough
path if it is controlled by $\omega$ in the sense of $\left(
\ref{finite-variation}\right)  $ and is "almost multiplicative", i.e. there
exists $\theta>1$ such that%
\[
\left\Vert \left(  X_{s,u}\otimes X_{u,t}\right)  ^{i}-X_{s,t}^{i}\text{
}\right\Vert \leq\omega\left(  s,t\right)  ^{\theta}\text{,}%
\]
for all $0\leq s$ $\leq u\leq t\leq T,$ $i=1,\ldots,\lfloor p\rfloor.$
\end{definition}

A crucial feature of the almost $p-$rough path that is exploited when defining
the integral is that there is exactly one $p$ rough path that is "close". More
precisely we have the following theorem (see \cite{lyons-stflour}, Theorem 4.3).

\begin{theorem}
Let $p\geq1,$ $\theta>1$ and $\omega$ a control. Let $\widetilde{X}$ be a $p$
-almost rough path controlled by $\omega.$ Then there exists a unique $p-$
rough path $X$ such that%
\[
\max_{i=0,\ldots,\lfloor p\rfloor}\sup_{0\leq s<t\leq T}\frac{\left\Vert
\widetilde{X}_{s,t}^{i}-X_{s,t}^{i}\text{ }\right\Vert }{\omega\left(
s,t\right)  ^{\theta}}<\infty.
\]

\end{theorem}

Let $\gamma>p$, $\alpha=\left(  \alpha^{0},\ldots,\alpha^{k}\right)  $ a
Lip-$\left(  \gamma-1\right)  $ one form in the sense of Stein and $X$ a
geometric $p$ -rough path. The next step in the development is to construct
(using..) an almost $p-$ rough path $Y$ corresponding to the integral of
$\alpha$ against $X.$ We refer to see Lyons, Caruana, Levy
\cite{lyons-stflour} Section 4.3.1 for the definition of $Y$ and the details
of the construction.

\begin{definition}
\label{rough-integral}Let $\gamma>p$, $\alpha$ be a Lip-$\left(
\gamma-1\right)  $ one form taking values in some Banach space $V$ and $X$ a
geometric $p$ -rough path. We define the rough integral of $\alpha$ against
the path $X$ to be the unique $p-$rough path associated to $Y$ and denote it
by $\int\alpha\left(  X_{t}\right)  dX.$
\end{definition}

The rough integral is continuous in the $p-$variation topology and maps
geometric $p-$ rough paths on $W$ to geometric $p-$rough paths on $V.$ The
following theorem may be found in \cite{lyons-stflour} (Theorem 4.12).

\begin{theorem}
With the notation of Definition \ref{rough-integral} the mapping
\[
X\rightarrow\int\alpha\left(  X_{t}\right)  dX
\]
is continuous from $\Omega G_{p}(W)$ to $\Omega G_{p}(V).$ If $\omega$
controls $X$ the rough path $\int\alpha\left(  X_{t}\right)  dX$ is controlled
by $C\omega,$ where $C$ is a constant that only depends on $\left\Vert
\alpha\right\Vert _{Lip},\gamma,p$ and $\omega\left(  0,T\right)  .$
\end{theorem}

For bounded variation paths the rough integral agrees with the usual
Riemann-Stieltjes integral. As a consequence many integration identities such
as the fundamental theorem of calculus are satisfied by the projection onto
level one of the rough integral of a geometric rough path. Uniqueness of the
lift of bounded variation paths together with the continuity are a powerful
tool to obtain such identities and will frequently be used in the remainder of
the paper.

We note that rough paths and rough integrals defined over an interval posses a
natural restriction to rough paths over subintervals. In more detail, we have

\begin{definition}
Let $X\in G\Omega_{p}\left(  V\right)  $ be defined over an interval $\left[
0,T\right]  $ and let $\left[  s,t\right]  \subseteq\left[  0,T\right]  .$Then
we define the restriction $X|_{\left[  s,t\right]  }:\Delta\left[  s,t\right]
\rightarrow T^{\left\lfloor p\right\rfloor }\left(  V\right)  $ by
\[
X|_{\left[  s,t\right]  }\left(  u,v\right)  =X_{u,v}%
\]
for $\left(  u,v\right)  $ in $\Delta\left[  s,t\right]  .$
\end{definition}

\begin{remark}
It follows immediately from this definition that the restriction $X|_{\left[
s,t\right]  }$ of an element of $G\Omega_{p}\left(  V\right)  $ is itself a
geometric $p-$rough path over the interval $\left[  s,t\right]  .$
\end{remark}

\begin{definition}
Let $X\in G\Omega_{p}\left(  V\right)  $ be defined over the interval $\left[
s,t\right]  $ with a starting point $x\in V$. Then for any Lip-$\left(
\gamma-1\right)  $ one form on $V$ we define the rough integral
\[
\int\alpha\left(  X\right)  dX
\]
to be the geometric $p-$rough path over $\left[  s,t\right]  $ defined by
\[
\int_{u}^{v}\alpha\left(  X\right)  dX=\int_{u-s}^{v-s}\alpha\left(  Y\right)
dY
\]
where $Y\in G\Omega_{p}\left(  V\right)  $ is defined over $\left[
0,t-s\right]  $ with increments $Y_{u,v}=X_{s+u,s+v}$ and starting point $x.$
\end{definition}

Let $V$ and $W$ be Banach spaces and $X\in G\Omega_{p}\left(  V\right)  $. Let
$g$ be a Lip-$\gamma$ map from $W\rightarrow$L$\left(  V,W\right)  .$ Then $g$
may be interpreted as a one form on $V$ taking values in the vector fields on
$W.$ If $X$ is a bounded variation path we consider the differential equation%
\begin{equation}
dY_{t}=g\left(  Y_{t}\right)  dX_{t},\text{ \ }Y\left(  0\right)  =Y_{0}.
\label{rde-banach-space1}%
\end{equation}
For $p\geq2$ the increments at level one do not determine a rough path. Thus,
a rough differential equation will be defined for the pair $\left(
X,Y\right)  ,$ living on $V\oplus W,$ that takes into account the cross
integrals between $X$ and $Y$, which do not appear in either projection of
$\left(  X,Y\right)  $ onto $V$ or $W$. Let $h:V\oplus W\rightarrow
$End$\left(  V\oplus W\right)  $ be defined by%
\begin{equation}
h(x,y)=\left(
\begin{array}
[c]{cc}%
\text{Id}_{V} & 0\\
f\left(  y+Y_{0}\right)  & 0
\end{array}
\right)  , \label{def-h}%
\end{equation}
If $X$ is a proper geometric rough path the equation $\left(
\ref{rde-banach-space1}\right)  $ may be reformulated as
\begin{equation}
dZ_{t}=h\left(  Z_{t}\right)  dZ_{t},\text{ }Z_{0}=0,\text{ }\pi_{V}\left(
Z\right)  =X,
\end{equation}
where $\pi_{V}$ denotes the projection onto $T^{(\lfloor p\rfloor)}(V).$ We
define a solution in terms of a fixed point of the rough integral.

\begin{definition}
Let $g:W\rightarrow$L$\left(  V,W\right)  $ be a Lip$\left(  \gamma-1\right)
$ function and $X\in G\Omega_{p}\left(  V\right)  .$ We say $Z\in G\Omega
_{p}\left(  V\oplus W\right)  $ solves the rough differential equation
$\left(  \text{RDE}\right)  $%
\begin{equation}
dY_{t}=g\left(  Y_{t}\right)  dX_{t},\text{ \ }Y\left(  0\right)  =Y_{0}
\label{rde-banach-space}%
\end{equation}
if $\pi_{V}\left(  Z\right)  =X$ and%
\[
Z=\int h\left(  Z\right)  dZ,
\]
where $h$ is defined as in $\left(  \ref{def-h}\right)  .$
\end{definition}

If $g$ is sufficiently regular the rough differential equations $\left(
\ref{rde-banach-space}\right)  $ satisfies a universal limit theorem, namely
it has a unique solution and the It\^{o} map taking the driving signal to the
response is continuous in the $p-$variation topology. The following theorem
may be found in \cite{lyons-stflour} (Theorem 5.3)

\begin{theorem}
[Universal Limit Theorem]\label{universal limit theorem}Let $p\geq1,$
$\gamma>p$ and
\[
g:W\rightarrow L\left(  V,W\right)
\]
be a Lip $\gamma$ function. For all $X\in G\Omega_{p}\left(  V\right)  $ and
$Y_{0}\in V$ the equation%
\[
dY_{t}=g\left(  Y_{t}\right)  dX_{t},\ Y\left(  0\right)  =Y_{0}%
\]
admits a unique solution. The map from $G\Omega_{p}\left(  V\right)  \times
W\rightarrow G\Omega_{p}\left(  W\right)  $ that takes $\left(  X,Y_{0}%
\right)  $ to $Y$ is continuous in the $p-$variation topology. Moreover if $X$
is controlled by $\omega$ the rough path $Y$ is controlled by $C\omega,$ where
$C$ is a constant that only depends on $\left\Vert g\right\Vert _{\text{Lip-}%
\gamma},$ $p,$ $\gamma$ and $\omega.$
\end{theorem}

The proof of the theorem proceeds by considering a sequence of Picard
iterations and uses a rescaling of the driving signal $X_{t}$ together with
the division property of Lipschitz functions. The constant $C$ in the theorem
can be made more explicit leading to Gronwall type inequalities, see e.g.
Cass, Lyons \cite{cass} for some new sharp estimates. If $g$ is only
Lip-$\left(  \gamma-1\right)  $ the equation $\left(  \ref{rde-banach-space}%
\right)  $ no longer has a unique solution, existence however is still
guaranteed $\left(  \text{see e.g. Davie \cite{davie}}\right)  .$

\section{Lip-$\gamma$ manifolds}

We begin by introducing the concept of Lipschitz-$\gamma$ (or Lip-$\gamma$ for
short) manifolds, based on the notion of Lip-$\gamma$ functions introduced by
Stein \cite{stein} for functions on Banach spaces (see Definition
\ref{lipschitz functions}). We do not follow the usual approach to introducing
a differential structure which typically begins with a pseudo group of
functions. Instead we assume the existence of a particularly nice Lip-$\gamma$
atlas that is amenable to the analysis involved in developing a theory of
rough paths on a Lip-$\gamma$ manifold. We show that any finite dimensional
normed vector space and any compact $C^{\lfloor\gamma\rfloor+1}$ manifold are
Lip-$\gamma$ manifolds. The current paper does not address some important
geometric issues yet. For example, we do not develop a notion of submanifolds.
The main focus is to establish sufficient conditions a particular atlas has to
satisfy to guarantee the existence of solutions of rough differential
equations defined by sufficiently regular vector fields.

Throughout we work with a topological manifold $M$ (Hausdorff, second
countable and locally homeomorphic to $%
%TCIMACRO{\U{211d} }%
%BeginExpansion
\mathbb{R}
%EndExpansion
^{d},$ for some $d>0).$ As usual, a chart on $M$ is a pair $(\phi,U)$; $U$ is
an open subset of $M$ but, in contrast to usual practice, the co-orinate map
$\phi$ will be defined on globally on $M$ with its restriction $\phi|_{U}$ to
$U$ being homeomorphic onto an open subset of $%
%TCIMACRO{\U{211d} }%
%BeginExpansion
\mathbb{R}
%EndExpansion
^{d}.$ We begin by giving the definition of a Lip-$\gamma$ atlas. In the
following let for $x\in%
%TCIMACRO{\U{211d} }%
%BeginExpansion
\mathbb{R}
%EndExpansion
^{d},$ $u>0$ denote $B\left(  x,u\right)  $ the ball of radius $u,$ centred at
$x$ taken with respect to some norm $\left\Vert \cdot\right\Vert $ (which is
not necessarily the Euclidean norm). Later we will consider product manifolds
as the Cartesian product and it will then be convenient to specify a norm on
the coordinate space with certain desirable properties; it is for this reason
that we do not restrict ourselves to any particular norm at this stage.

\begin{definition}
\label{l-gamma-atlas}A Lipschitz-$\gamma$ atlas with Lipschitz constant $L$ on
$M$ is a countable collection of charts $\left\{  \left(  \phi_{i}%
,U_{i}\right)  :i\in I\right\}  $ with the following properties :

\begin{enumerate}
\item The sets $\left\{  U_{i}:i\in I\right\}  $ form a pre-compact and a
locally finite cover for $M.$

\item For all charts $\left(  \phi_{i},U_{i}\right)  $ we have $\phi_{i}%
(U_{i})=B_{\left\vert \left\vert \cdot\right\vert \right\vert }(0,1)$ and
there exists $R>0$ such that
\[
\text{ }\phi_{i}(M)\subseteq B_{\left\vert \left\vert \cdot\right\vert
\right\vert }\left(  0,R\right)  \subseteq\left(
%TCIMACRO{\U{211d} }%
%BeginExpansion
\mathbb{R}
%EndExpansion
^{d},\left\vert \left\vert \cdot\right\vert \right\vert \right)  .
\]

\item There is a $\delta$ in $\left(  0,1\right)  $ such that the open sets%
\begin{equation}
U_{i}^{\delta}=\phi_{i}|_{U_{i}}^{-1}\left(  B_{\left\vert \left\vert
\cdot\right\vert \right\vert }\left(  0,1-\delta\right)  \right)
\label{delta-sets}%
\end{equation}

cover $M.$

\item The functions $\phi:M\rightarrow%
%TCIMACRO{\U{211d} }%
%BeginExpansion
\mathbb{R}
%EndExpansion
^{d}$ are compactly supported and for any pair of charts $\left(  \phi
_{i},U_{i}\right)  $ and $\left(  \phi_{j},U_{j}\right)  $ the map $\phi
_{i}\circ\phi_{j}|_{U_{j}}^{-1}$ \bigskip$:\phi_{j}\left(  U_{j}\right)
\rightarrow\left(
%TCIMACRO{\U{211d} }%
%BeginExpansion
\mathbb{R}
%EndExpansion
^{d},\left\vert \left\vert \cdot\right\vert \right\vert \right)  $ is
Lip-$\gamma$ with
\[
\left\vert \left\vert \phi_{i}\circ\phi_{j}|_{U_{j}}^{-1}\right\vert
\right\vert _{Lip-\gamma}\leq L.
\]

\end{enumerate}
\end{definition}

The following elementary reformulation of the condition 3 in this definition
will be useful to us at times.

\begin{lemma}
Condition 3 of Definition \ref{l-gamma-atlas} is equivalent to the statement
that for all $m\in M$ there exists a chart $\left(  \phi,U\,\right)  $
containing $m$ with the property that%
\[
\overline{B\left(  \phi\left(  x\right)  ,\delta\right)  }\subseteq\phi\left(
U\right)  =B_{\left\vert \left\vert \cdot\right\vert \right\vert }\left(
0,1\right)
\]
.
\end{lemma}

Henceforth, we will not make any reference to the norm on the coordinate space
unless it is specifically relevant to the situation under discussion.

\begin{remark}
It can be shown (see e.g. Lee \cite{lee}) that any smooth manifold has an
atlas satisfying conditions one to three.
\end{remark}

\begin{definition}
A Lip-$\gamma$ manifold with Lipschitz constant $L$ is a topological manifold
with a Lip-$\gamma$ atlas with constant $L$.
\end{definition}

Later when we discuss solutions to RDEs on manifolds it will be important that
we can combine a signal rough path living on one Lip-$\gamma$ manifold $N$ and
a response rough path on another Lip-$\gamma$ manifold $M.$ This will involve
imposing the structure of a Lip-$\gamma$ manifold on the topological space
$N\times M$ (taken with the product topology). For definiteness, suppose that
$\left\{  \left(  \phi_{i},U_{i}\right)  :i\in I\right\}  $ is a Lip$-\gamma$
atlas for $N$ and that the coordinate spaces are subsets of the normed space
$\left(
%TCIMACRO{\U{211d} }%
%BeginExpansion
\mathbb{R}
%EndExpansion
^{d_{1}},\left\vert \left\vert \cdot\right\vert \right\vert _{1}\right)  $ and
likewise that $\left\{  \left(  \psi_{j},V_{j}\right)  :j\in J\right\}  $ is a
Lip$-\gamma$ atlas for $M$ coordinatised through subsets of $\left(
%TCIMACRO{\U{211d} }%
%BeginExpansion
\mathbb{R}
%EndExpansion
^{d_{2}},\left\vert \left\vert \cdot\right\vert \right\vert _{2}\right)  .$ We
then introduce the candidate Lip-$\gamma$ atlas for the product as%
\begin{equation}
\left\{  \left(  \xi_{ij},W_{ij}\right)  :i\in I,\text{ }j\in J\right\}
\label{product atlas}%
\end{equation}
where $W_{ij}=U_{i}\times V_{j}\subseteq N\times M$ and $\xi_{ij}:N\times
M\rightarrow%
%TCIMACRO{\U{211d} }%
%BeginExpansion
\mathbb{R}
%EndExpansion
^{d_{1}}\oplus%
%TCIMACRO{\U{211d} }%
%BeginExpansion
\mathbb{R}
%EndExpansion
^{d_{2}}$ is defined by $\xi_{ij}\left(  n,m\right)  =\left(  \phi_{i}\left(
n\right)  ,\psi_{j}\left(  m\right)  \right)  .$ The following lemma shows
that, for a particular choice of direct sum norm on the coordinate space $%
%TCIMACRO{\U{211d} }%
%BeginExpansion
\mathbb{R}
%EndExpansion
^{d_{1}}\oplus%
%TCIMACRO{\U{211d} }%
%BeginExpansion
\mathbb{R}
%EndExpansion
^{d_{2}},$ the atlas (\ref{product atlas}) is a Lip-$\gamma$ atlas for
$N\times M$.

\begin{lemma}
Let $N$ and $M$ be two Lip$-\gamma$ manifolds with associated constants
$L_{1}$ and $L_{2}$, $\delta_{1}$ and $\delta_{2}$, $R_{1}$ and $R_{2},$
atlases $\left\{  \left(  \phi_{i},U_{i}\right)  :i\in I\right\}  $ and
$\left\{  \left(  \psi_{j},V_{j}\right)  :j\in I\right\}  $ and such that
their coordinates spaces are subsets of $\left(
%TCIMACRO{\U{211d} }%
%BeginExpansion
\mathbb{R}
%EndExpansion
^{d_{1}},\left\vert \left\vert \cdot\right\vert \right\vert _{1}\right)  $ and
$\left(
%TCIMACRO{\U{211d} }%
%BeginExpansion
\mathbb{R}
%EndExpansion
^{d_{2}},\left\vert \left\vert \cdot\right\vert \right\vert _{2}\right)  $
respectively. Then (\ref{product atlas}) is a Lip-$\gamma$ atlas with
Lipschitz constant bounded by $L:=\max\left(  L_{1},L_{2}\right)  $ when the
coordinate space $%
%TCIMACRO{\U{211d} }%
%BeginExpansion
\mathbb{R}
%EndExpansion
^{d_{1}}\oplus%
%TCIMACRO{\U{211d} }%
%BeginExpansion
\mathbb{R}
%EndExpansion
^{d_{2}}$ is endowed with the direct sum norm:%
\[
\left\vert \left\vert \left(  x,y\right)  \right\vert \right\vert =\max\left(
\left\vert \left\vert x\right\vert \right\vert _{1},\left\vert \left\vert
y\right\vert \right\vert _{2}\right)  .
\]

\end{lemma}

\begin{proof}
We need to verify the conditions of Definition \ref{l-gamma-atlas}. It is
elementary to check that $\left\{  W_{ij}:i\in I,\text{ }j\in J\right\}  $
\ forms a cover for the product space. Condition 2 follows more or less
immediately from the (crucial) observation that for any $r>0$ we have
\[
B_{\left\vert \left\vert \cdot\right\vert \right\vert }\left(  0,r\right)
=B_{\left\vert \left\vert \cdot\right\vert \right\vert _{1}}\left(
0,r\right)  \times B_{\left\vert \left\vert \cdot\right\vert \right\vert _{2}%
}\left(  0,r\right)  ,
\]
so that
\begin{align*}
\xi_{ij}\left(  W_{ij}\right)   &  =\xi_{ij}\left(  U_{i}\times V_{j}\right)
\\
&  =\phi_{i}\left(  U_{i}\right)  \times\psi_{j}\left(  V_{j}\right) \\
&  =B_{\left\vert \left\vert \cdot\right\vert \right\vert _{1}}\left(
0,1\right)  \times B_{\left\vert \left\vert \cdot\right\vert \right\vert _{2}%
}\left(  0,1\right) \\
&  =B_{\left\vert \left\vert \cdot\right\vert \right\vert }\left(  0,1\right)
.
\end{align*}
Easy calculations using the same idea can be used to show that $\xi
_{ij}\left(  N\times M\right)  \subseteq B_{\left\vert \left\vert
\cdot\right\vert \right\vert }\left(  0,R\right)  $ for $R:=\max\left(
R_{1},R_{2}\right)  $. \ To check condition 3 we need to observe that if
$\delta:=\min\left(  \delta_{1},\delta_{2}\right)  $ then
\begin{align*}
W_{ij}^{\delta}  &  =\xi_{ij}|_{W_{ij}}^{-1}\left(  B_{\left\vert \left\vert
\cdot\right\vert \right\vert }\left(  0,1-\delta\right)  \right) \\
&  =\xi_{ij}|_{W_{ij}}^{-1}\left(  B_{\left\vert \left\vert \cdot\right\vert
\right\vert _{1}}\left(  0,1-\delta\right)  \times B_{\left\vert \left\vert
\cdot\right\vert \right\vert _{2}}\left(  0,1-\delta\right)  \right) \\
&  =\phi_{i}|_{U_{i}}^{-1}\left(  B_{\left\vert \left\vert \cdot\right\vert
\right\vert _{1}}\left(  0,1-\delta\right)  \right)  \times\psi_{j}|_{V_{j}%
}^{-1}\left(  B_{\left\vert \left\vert \cdot\right\vert \right\vert _{2}%
}\left(  0,1-\delta\right)  \right) \\
&  =U_{i}^{\delta}\times V_{j}^{\delta}\\
&  \supseteq U_{i}^{\delta_{1}}\times V_{j}^{\delta_{2}}.
\end{align*}
It follows that $\left\{  W_{ij}^{\delta}:i\in I,\text{ }j\in J\right\}  $
covers $N\times M$. Finally, condition 4 follows because
\begin{equation}
\xi_{i_{1}j_{1}}\circ\xi_{i_{2}j_{2}}|_{W_{i_{2}j_{2}}}^{-1}=\left(
\phi_{i_{1}}\circ\phi_{i_{2}}|_{U_{i_{2}}}^{-1},\psi_{j_{1}}\circ\psi_{j_{2}%
}|_{V_{j_{2}}}^{-1}\right)  \label{product transition}%
\end{equation}
and the observation that the choice of norm on $%
%TCIMACRO{\U{211d} }%
%BeginExpansion
\mathbb{R}
%EndExpansion
^{d_{1}}\oplus%
%TCIMACRO{\U{211d} }%
%BeginExpansion
\mathbb{R}
%EndExpansion
^{d_{2}}$ implies that any $f:B_{\left\vert \left\vert \cdot\right\vert
\right\vert }\left(  0,1\right)  \rightarrow%
%TCIMACRO{\U{211d} }%
%BeginExpansion
\mathbb{R}
%EndExpansion
^{d_{1}}\oplus%
%TCIMACRO{\U{211d} }%
%BeginExpansion
\mathbb{R}
%EndExpansion
^{d_{2}}$ which is of the form $f=\left(  f_{1},f_{2}\right)  $ where
$f_{1}:B_{\left\vert \left\vert \cdot\right\vert \right\vert _{1}}\left(
0,1\right)  \rightarrow%
%TCIMACRO{\U{211d} }%
%BeginExpansion
\mathbb{R}
%EndExpansion
,$ $i=1,2$ and where the $f_{i}$ are Lip-$\gamma$ will itself be Lip-$\gamma$
and
\[
\left\vert \left\vert f\right\vert \right\vert _{Lip-\gamma}\leq\max\left(
\left\vert \left\vert f_{1}\right\vert \right\vert _{Lip-\gamma},\left\vert
\left\vert f_{2}\right\vert \right\vert _{Lip-\gamma}\right)  .
\]
The proof is completed by combining this observation with
(\ref{product transition}).
\end{proof}

It turns out that it is easy to show that finite dimensional normed vector
spaces and compact $C^{\lfloor\gamma\rfloor+1}$ manifolds are Lip-$\gamma$ manifolds.

\begin{example}
\label{banach-manifold}A finite dimensional normed vector space $W$ is a
Lip-$\gamma$ manifold. To see this we may cover $W$ with balls $B\left(
x_{i},1-\delta\right)  $. The charts $\left(  \phi_{i},B(x_{i},1\right)  ),$
where $\phi_{i}$ is any extension map that translates the ball $B(x_{i},1)$ to
the unit ball centred at the origin, are a Lip-$\gamma$ atlas for any
$\gamma\geq1.$
\end{example}

\begin{lemma}
Any compact $C^{\lfloor\gamma\rfloor+1}$-manifold is a Lip-$\gamma$ manifold.
\end{lemma}

\begin{proof}
It is elementary to show that any $C^{\lfloor\gamma\rfloor+1}$-manifold has an
atlas consisting of charts $\left(  \phi_{i},U_{i}\right)  $ that is a regular
refinement, i.e. satisfies conditions one and three of a Lip-$\gamma$ atlas
with $\delta=1$ and in addition satisfies $\phi_{i}\left(  U_{i}\right)
=B\left(  0,3\right)  $ (see e.g. Lee \cite{lee} page 53). Using the
compactness of $M$ we can find a finite number of sets $U_{i_{j}}^{\delta}$
defined as in (\ref{delta-sets}) that cover the manifold. The set
\[
V=\overline{\phi_{i_{j}}\left(  U_{i_{j}}^{\delta/2}\right)  }=\overline
{B\left(  0,5/2\right)  }%
\]
is compact and therefore $U_{i_{j}}^{\delta/2}\subseteq\phi_{i_{j}}%
^{-1}\left(  V\right)  $ is pre compact. The functions $\phi_{i_{j}%
}|_{U_{i_{k}}^{\delta/3}}$ may be extended to $C^{\lfloor\gamma\rfloor+1}$
functions $\widetilde{\phi}_{i_{j}}$ that vanish off $U_{i_{k}}^{\delta/4}$
and thus are compactly supported. Also note that the functions
$\widetilde{\phi}_{i_{j}}\circ\widetilde{\phi}_{i_{k}}^{-1}|_{B\left(
0,8/9\right)  }$ are $C^{\lfloor\gamma\rfloor+1}$ and therefore as it has
$\lfloor\gamma\rfloor+1$ bounded derivatives on $B\left(  0,5/2\right)  $ the
function $\widetilde{\phi}_{i_{j}}\circ\widetilde{\phi}_{i_{k}}^{-1}%
|_{B\left(  0,5/2\right)  }$ is Lip-$\gamma$ with some constant $L_{j,k}.$
Hence, the finite collection of charts $\left(  \widetilde{\phi}_{i_{j}%
},U_{i_{j}}^{\delta/2}\right)  $ also satisfies condition 4. Note that
$\frac{2}{5}\widetilde{\phi}_{i_{j}}\left(  U_{i_{j}}^{\delta/2}\right)
=B\left(  0,1\right)  $ and $\widetilde{\phi}_{i_{j}}\left(  M\right)
\subseteq B(0,3)$ so condition 2. is also satisfied and the collection
$\left(  \frac{2}{5}\widetilde{\phi}_{i_{j}},U_{i_{j}}^{\delta/2}\right)  $ is
a Lip-$\gamma$ atlas.
\end{proof}

The next step is to introduce Lip-$\gamma$ functions on Lip-$\gamma$ manifolds
taking values in Banach spaces.

\begin{definition}
\label{lipschitz function manifold}Let $f:M\rightarrow W$ be a function
between a Lip-$\gamma$ manifold $M$ and a Banach space $W.$ We say $f$ is
Lip-$\gamma$ if there exists a $C\geq0$ such that for every chart $(\phi
_{i},U_{i})$ on $M$ the functions $f\circ\phi_{i}|_{\phi\left(  U_{i}\right)
}^{-1}:\phi_{i}(U_{i})\rightarrow W$ are Lip-$\gamma$ with Lipschitz constant
at most $C$ in the sense of Definition \ref{lipschitz functions}. We say the
minimal $C$ satisfying this requirement is the Lipschitz constant of $f$ and
denoted it by $\left\Vert f\right\Vert _{Lip-\gamma}$.
\end{definition}

The Lipschitz constants define a norm on the Lip-$\gamma$ functions on $M$
taking values in $W$. We can similarly define Lip-$\gamma$ functions on any
open subset $V$ of the manifold by considering the intersection of all charts
with the set $V.$ We now introduce the important notion of the equivalence
between two Lip-$\gamma$ atlases on a manifold $M.$ The important feature that
we need to capture is that equivalent atlases should give rise to the same set
of Lip-$\gamma$ functions on the manifold. Thus, we make the following definition.

\begin{definition}
[Equivalence of Lip-$\gamma$ atlases]Suppose that the sets
\[
\mathcal{A}_{1}=\left\{  \left(  \phi_{i},U_{i}\right)  :i\in I\right\}
\text{ and }\mathcal{A}_{2}=\left\{  \left(  \psi_{j},V_{j}\right)  :j\in
J\right\}
\]
define two Lip$-\gamma$ atlases on $M$ with respective coordinate spaces
$\left(
%TCIMACRO{\U{211d} }%
%BeginExpansion
\mathbb{R}
%EndExpansion
^{d},\left\vert \cdot\right\vert _{1}\right)  $ and $\left(
%TCIMACRO{\U{211d} }%
%BeginExpansion
\mathbb{R}
%EndExpansion
^{d},\left\vert \cdot\right\vert _{2}\right)  $ $.$ Then we say that
$\mathcal{A}_{1}$ and $\mathcal{A}_{2}$ are equivalent with constant $C$ if
for every $i$ in $I$ and $j$ in $J$ we have that the functions
\begin{align*}
\phi_{i}\circ\psi_{j}|_{V_{j}}^{-1}  &  :B_{\left\vert \cdot\right\vert _{2}%
}\left(  0,1\right)  \rightarrow\left(
%TCIMACRO{\U{211d} }%
%BeginExpansion
\mathbb{R}
%EndExpansion
^{d},\left\vert \cdot\right\vert _{1}\right) \\
\psi_{j}\circ\phi_{i}|_{U_{i}}^{-1}  &  :B_{\left\vert \cdot\right\vert _{1}%
}\left(  0,1\right)  \rightarrow\left(
%TCIMACRO{\U{211d} }%
%BeginExpansion
\mathbb{R}
%EndExpansion
^{d},\left\vert \cdot\right\vert _{2}\right)
\end{align*}
are Lip-$\gamma$ in the sense of Stein and moreover
\[
\sup_{i\in I,j\in J}\left\vert \left\vert \phi_{i}\circ\psi_{j}|_{V_{j}}%
^{-1}\right\vert \right\vert _{Lip-\gamma}\leq C\text{ and }\sup_{i\in I,j\in
J}\left\vert \left\vert \psi_{j}\circ\phi_{i}|_{U_{i}}^{-1}\right\vert
\right\vert _{Lip-\gamma}\leq C.\text{ }%
\]

\end{definition}

The following technical lemma will be useful in characterising the notion of
equivalence atlases. It tells us that the uniform control on the Lip-$\gamma$
norm of the transition functions in the definition of equivalence allows us to
control the size of balls when passing between the two coordinate systems.
More precisely, we have:

\begin{lemma}
\label{connectedness}Suppose that $\mathcal{A}_{1}=\left\{  \left(  \phi
_{i},U_{i}\right)  :i\in I\right\}  $ and $\mathcal{A}_{2}=\left\{  \left(
\psi_{j},V_{j}\right)  :j\in J\right\}  $ define two Lip$-\gamma$ atlases on
$M$ with coordinate spaces $\left(
%TCIMACRO{\U{211d} }%
%BeginExpansion
\mathbb{R}
%EndExpansion
^{d},\left\vert \cdot\right\vert _{1}\right)  $ and $\left(
%TCIMACRO{\U{211d} }%
%BeginExpansion
\mathbb{R}
%EndExpansion
^{d},\left\vert \cdot\right\vert _{2}\right)  $ respectively. Suppose further
that $\mathcal{A}_{1}$and $\mathcal{A}_{2}$ are equivalent with finite
constant $C\geq1.$ Let $x$ be in $B_{\left\vert \cdot\right\vert _{1}}\left(
0,1\right)  $; for any chart $\left(  \phi,U\right)  $ in $\mathcal{A}_{1}$ we
have that $m:=\phi|_{U}^{-1}\left(  x\right)  $ is in $U.$ Moreover, if
$\left(  \psi,V\right)  $ is any other chart in the second atlas
$\mathcal{A}_{2}$ which also\ contains $m$ and is such that for some $r>0$%
\[
B_{\left\vert \cdot\right\vert _{2}}\left(  \psi\left(  m\right)  ,r\right)
\subseteq B_{\left\vert \cdot\right\vert _{2}}\left(  0,1\right)
\]
then we have that
\begin{equation}
\phi|_{U}^{-1}\left(  B_{\left\vert \cdot\right\vert _{1}}\left(  x,u\right)
\cap B_{\left\vert \cdot\right\vert _{1}}\left(  0,1\right)  \right)
\subseteq\psi|_{V}^{-1}\left(  B_{\left\vert \cdot\right\vert _{2}}\left(
\psi\left(  m\right)  ,r\right)  \right)  \label{contain}%
\end{equation}
provided $0\leq u<\frac{r}{C}.$
\end{lemma}

\begin{proof}
Suppose for a contradiction that (\ref{contain}) does not hold; then for some
$y$ in $B_{\left\vert \cdot\right\vert _{1}}\left(  0,1\right)  $ with
$\left\vert x-y\right\vert _{1}\leq u$ we have that $\phi|_{U}^{-1}\left(
y\right)  $ is in $M\setminus\psi|_{V}^{-1}\left(  B_{\left\vert
\cdot\right\vert _{2}}\left(  \psi\left(  m\right)  ,r\right)  \right)  $.
Define $w:\left[  0,1\right]  \rightarrow B_{\left\vert \cdot\right\vert _{1}%
}\left(  0,1\right)  $ \ by taking $w\left(  \lambda\right)  =\left(
1-\lambda\right)  x+\lambda y$; it follows that $g:\left[  0,1\right]
\rightarrow M$ defined by
\[
g\left(  \lambda\right)  =\phi|_{U}^{-1}\left(  w\left(  \lambda\right)
\right)
\]
is a continuous path with values in $U\subseteq M.$ We have $g\left(
0\right)  =\phi|_{U}^{-1}\left(  x\right)  $ is in $\psi|_{V}^{-1}\left(
B_{\left\vert \cdot\right\vert _{2}}\left(  \psi\left(  m\right)  ,r\right)
\right)  ,$ and it must happen that for some $\lambda$ in $(0,1]$ the path $g$
is in the set
\begin{equation}
\psi|_{V}^{-1}\left(  B_{\left\vert \cdot\right\vert _{2}}\left(  \psi\left(
m\right)  ,r\right)  \setminus B_{\left\vert \cdot\right\vert _{2}}\left(
\psi\left(  m\right)  ,Cu\right)  \right)  ; \label{crossing set}%
\end{equation}
if this is were not the case then for any continuous function $\Phi
:M\rightarrow%
%TCIMACRO{\U{211d} }%
%BeginExpansion
\mathbb{R}
%EndExpansion
$ such that
\[
\Phi\left(  m\right)  =\left\{
\begin{array}
[c]{cc}%
1 & \text{on }\psi|_{V}^{-1}\left(  B_{\left\vert \cdot\right\vert _{1}%
}\left(  \psi\left(  m\right)  ,Cu\right)  \right)  \text{ \ }\\
0 & \text{ \ \ on }M\setminus\psi|_{V}^{-1}\left(  B_{\left\vert
\cdot\right\vert _{2}}\left(  \psi\left(  m\right)  ,r\right)  \right)
\end{array}
\right.
\]
we would have that $\Phi\circ g:\left[  0,1\right]  \rightarrow%
%TCIMACRO{\U{211d} }%
%BeginExpansion
\mathbb{R}
%EndExpansion
$ is a continuous function with values in $\left\{  0,1\right\}  $ and hence
constant. It follows that $\left(  \Phi\circ g\right)  \left(  1\right)
=\left(  \Phi\circ g\right)  \left(  0\right)  =1,$ contradicting $\phi
|_{U}^{-1}\left(  y\right)  $ not being in $\psi|_{V}^{-1}\left(
B_{\left\vert \cdot\right\vert _{2}}\left(  \psi\left(  m\right)  ,r\right)
\right)  .$ Hence, by choosing $\lambda$ such that $g\left(  \lambda\right)  $
is in (\ref{crossing set}) we can deduce
\[
\left\vert \psi\left(  g\left(  \lambda\right)  \right)  -\psi\left(
m\right)  \right\vert _{2}=\left\vert \left(  \psi\circ\phi|_{U}^{-1}\right)
\left(  w\left(  \lambda\right)  \right)  -\left(  \psi\circ\phi|_{U}%
^{-1}\right)  \left(  x\right)  \right\vert _{2}\geq Cu.
\]
On the other hand, the fact that $\psi$ and $\phi$ are drawn from equivalent
atlases implies that
\[
\left\vert \left(  \psi\circ\phi|_{U}^{-1}\right)  \left(  w\left(
\lambda\right)  \right)  -\left(  \psi\circ\phi|_{U}^{-1}\right)  \left(
x\right)  \right\vert _{2}\leq C\left\vert w\left(  \lambda\right)
-x\right\vert _{1}\leq C\left\vert x-y\right\vert _{1}<Cu;
\]
the two inequalities stand in contradiction to one another and the proof is complete.
\end{proof}

We record the following corollary for later use.

\begin{corollary}
\label{last resort}\bigskip Let $\ \mathcal{A}$ be a Lip-$\gamma$ atlas on a
manifold $M$ with constant $L.$ Suppose that $\left(  \phi,U\right)  $ and
$\left(  \psi\,,V\right)  $ are two charts and $m$ in $M$ is such that $m\in
U\cap V$ . Suppose there exists a $r>0$ such that $B\left(  \phi\left(
m\right)  ,r\right)  \subseteq B\left(  0,1\right)  $ and $B\left(
\psi\left(  m\right)  ,r\right)  \subseteq B\left(  0,1\right)  .$ Then we
have
\[
\psi|_{U}^{-1}\left(  B\left(  \psi\left(  x\right)  ,u\right)  \right)
\subseteq\phi|_{V}^{-1}\left(  B\left(  \phi\left(  x\right)  ,r\right)
\right)  ,
\]
provided $0\leq u<\frac{r}{L}$.
\end{corollary}

\begin{proof}
The claim is immediate from the previous lemma upon noticing that
$\mathcal{A}$ is equivalent to itself with constant $L.$
\end{proof}

The following lemma is the key result on equivalent Lip-$\gamma$ atlases. It
tells us that two equivalent atlases yield the same set of Lip-$\gamma$
functions up to changes of constant.

\begin{lemma}
\label{equiv}Let $\gamma\geq\gamma^{\prime}\geq1$ and suppose that
$\mathcal{A}_{1}=\left\{  \left(  \phi_{i},U_{i}\right)  :i\in I\right\}  $
and $\mathcal{A}_{2}=\left\{  \left(  \psi_{j},V_{j}\right)  :j\in J\right\}
$ define two Lip$-\gamma$ atlases on $M$ with respective constants $\delta
_{1}$ and $\delta_{2}$ . Suppose further that $\mathcal{A}_{1}$and
$\mathcal{A}_{2}$ are equivalent with finite constant $C.$ Then a Banach space
valued function $f:M\rightarrow W$ on $M$ is Lip-$\gamma^{\prime}$ with
respect to $\mathcal{A}_{1}$ on $M$ if and only if it is Lip$-\gamma^{\prime}$
with respect to $\mathcal{A}_{2}$ on $M.$ Furthermore, denoting the Lipschitz
norms of $f$ on $\mathcal{A}_{1}$and $\mathcal{A}_{2}$ by $\left\vert
\left\vert f\right\vert \right\vert _{1,Lip-\gamma}$ and $\left\vert
\left\vert f\right\vert \right\vert _{2,Lip-\gamma}$ we have that these norms
are equivalent; more precisely, we have
\begin{equation}
c\left(  C,\gamma,\delta_{1}\right)  \left\vert \left\vert f\right\vert
\right\vert _{2,Lip-\gamma}\leq\left\vert \left\vert f\right\vert \right\vert
_{1,Lip-\gamma}\leq d\left(  C,\gamma,\delta_{2}\right)  \left\vert \left\vert
f\right\vert \right\vert _{2,Lip-\gamma} \label{equivalence}%
\end{equation}
for some finite constants $c\left(  C,\gamma,\delta_{1}\right)  $ and
$d\left(  C,\gamma,\delta_{2}\right)  $
\end{lemma}

\begin{proof}
Suppose that $f$ is Lip$-\gamma^{\prime}$ with respect to the atlas
$\mathcal{A}_{2}$, we need to show that for any chart $\left(  \phi_{i}%
,U_{i}\right)  $ in $\mathcal{A}_{1}$ the function
\[
f\circ\phi_{i}|_{U_{i}}^{-1}:B_{\left\vert \cdot\right\vert _{1}}\left(
0,1\right)  \rightarrow W
\]
is $Lip-\gamma^{\prime}$ and has a Lipschitz norm that can be bounded
uniformly in $i\in I$. \ Due to the convexity of $B_{\left\vert \cdot
\right\vert _{1}}\left(  0,1\right)  $ it suffices, by Lemma
\ref{convex domains} to prove that $f\circ\phi_{i}|_{U_{i}}^{-1}$ is $k$-times
differentiable (where $k$ is the unique integer such that $\gamma^{\prime}%
\in(k,k+1]$) and that the $k^{th}$ derivative $D^{k}\left(  f\circ\phi
_{i}|_{U_{i}}^{-1}\right)  $ is $\left(  \gamma^{\prime}-k\right)
-$H\"{o}lder continuous. To prove this let us first that for every $x$ in
$B_{\left\vert \cdot\right\vert _{1}}\left(  0,1\right)  $ there exists a
chart $\left(  \psi,V\right)  $ in $\mathcal{A}_{2}$ such that $\phi
_{i}|_{U_{i}}^{-1}\left(  x\right)  $ is in $V$ and hence on the neighbourhood
$\phi_{i}\left(  U_{i}\cap V\right)  $ of $x$ we have that
\[
f\circ\phi_{i}|_{U_{i}}^{-1}=f\circ\psi|_{V}^{-1}\circ\psi\circ\phi
_{i}|_{U_{i}}^{-1}.
\]
Because the right hand side is the composition of two Lip-$\gamma^{\prime}$
functions it is itself Lip-$\gamma^{\prime}$ and we can deduce that
$f\circ\phi_{i}|_{U_{i}}^{-1}$ is $k-$times differentiable on this
neighbourhood and that its derivatives of these two functions agree. Using the
bounds in Lemma \ref{composition of lip functions} we can deduce that for
$j=0,1....,k$%
\begin{align}
\sup_{x\in B_{\left\vert \cdot\right\vert _{1}}\left(  0,1\right)  }\left\vert
D^{j}\left(  f\circ\phi_{i}|_{U_{i}}^{-1}\right)  \left(  x\right)
\right\vert  &  \leq\sup_{x\in B_{\left\vert \cdot\right\vert _{1}}\left(
0,1\right)  }\left\vert D^{j}\left(  f\circ\psi|_{V}^{-1}\circ\psi\circ
\phi_{i}|_{U_{i}}^{-1}\right)  \left(  x\right)  \right\vert \nonumber\\
&  \leq c_{1}\left(  C,\gamma\right)  \left\vert \left\vert f\right\vert
\right\vert _{2,Lip-\gamma^{\prime}}. \label{uniform bound on derivatives}%
\end{align}

It remains to prove that $D^{k}\left(  f\circ\phi_{i}|_{U_{i}}^{-1}\right)  $
is $\left(  \gamma^{\prime}-k\right)  -$H\"{o}lder continuous. As a first step
towards proving this we notice that if $x$ is in $B_{\left\vert \cdot
\right\vert _{1}}\left(  0,1\right)  $ then using condition 3 of Definition
\ref{l-gamma-atlas} there exists a chart $\left(  \psi^{x},V^{x}\right)  $ in
$\mathcal{A}_{2}$ such that $\phi_{i}|_{U_{i}}^{-1}\left(  x\right)  $ is in
$V^{x}$ and moreover that the preimage under $\psi^{x}$ of the ball (in
$\left\vert \cdot\right\vert _{2}$) of radius $\delta_{2}$ centred on $\left(
\psi^{x}\circ\phi_{i}|_{U_{i}}^{-1}\right)  \left(  x\right)  $ is contained
in $V^{x};$ that is we have
\begin{equation}
\psi^{x}|_{V^{x}}^{-1}\left[  B_{\left\vert \cdot\right\vert _{2}}\left(
\left(  \psi^{x}\circ\phi_{i}|_{U_{i}}^{-1}\right)  \left(  x\right)
,\delta_{2}\right)  \right]  \subseteq V^{x}. \label{delta ball}%
\end{equation}
We now observe that since $\mathcal{A}_{1}$and $\mathcal{A}_{2}$ are
equivalent atlases Lemma \ref{connectedness} together with (\ref{delta ball})
shows that
\begin{equation}
\phi|_{U}^{-1}\left(  B_{\left\vert \cdot\right\vert _{1}}\left(  0,1\right)
\cap B_{\left\vert \cdot\right\vert _{1}}\left(  x,\frac{\delta_{2}}%
{2C}\right)  \right)  \subseteq\psi^{x}|_{V^{x}}^{-1}\left[  B_{\left\vert
\cdot\right\vert _{2}}\left(  \left(  \psi^{x}\circ\phi_{i}|_{U_{i}}%
^{-1}\right)  \left(  x\right)  ,\delta_{2}\right)  \right]  \subseteq V^{x}.
\label{chart delta condition}%
\end{equation}
An elementary calculation then gives the bound
\begin{align}
&  \sup_{x,y\in B_{\left\vert \cdot\right\vert _{1}}\left(  0,1\right)  }%
\frac{\left\vert D^{k}\left(  f\circ\phi_{i}|_{U_{i}}^{-1}\right)  \left(
x\right)  -D^{k}\left(  f\circ\phi_{i}|_{U_{i}}^{-1}\right)  \left(  y\right)
\right\vert }{\left\vert x-y\right\vert ^{\gamma^{\prime}-k}}\nonumber\\
&  \leq\sup_{\substack{x,y\in B_{\left\vert \cdot\right\vert _{1}}\left(
0,1\right)  ,\\\left\vert x-y\right\vert _{1}<\frac{\delta_{2}}{2C}}%
}\frac{\left\vert D^{k}\left(  f\circ\phi_{i}|_{U_{i}}^{-1}\right)  \left(
x\right)  -D^{k}\left(  f\circ\phi_{i}|_{U_{i}}^{-1}\right)  \left(  y\right)
\right\vert }{\left\vert x-y\right\vert ^{\gamma^{\prime}-k}}\nonumber\\
&  +2\left(  \frac{2C}{\delta_{2}}\right)  ^{\gamma^{\prime}-k}\left\vert
D^{k}\left(  f\circ\phi_{i}|_{U_{i}}^{-1}\right)  \right\vert _{\infty
;B_{\left\vert \cdot\right\vert _{1}}\left(  0,1\right)  }
\label{kth derivative}%
\end{align}
We can use (\ref{chart delta condition}) to deduce, again using Lemma
\ref{composition of lip functions}, that for any $y$ in $B_{\left\vert
\cdot\right\vert _{1}}\left(  x,\frac{\delta_{2}}{2C}\right)  \cap
B_{\left\vert \cdot\right\vert _{1}}\left(  0,1\right)  $ we have
\begin{align*}
&  \left\vert D^{k}\left(  f\circ\phi_{i}|_{U_{i}}^{-1}\right)  \left(
x\right)  -D^{k}\left(  f\circ\phi_{i}|_{U_{i}}^{-1}\right)  \left(  y\right)
\right\vert \\
&  =\left\vert D^{k}\left(  f\circ\psi^{x}|_{V^{x}}^{-1}\circ\psi^{x}\circ
\phi_{i}|_{U_{i}}^{-1}\right)  \left(  x\right)  -D^{k}\left(  f\circ\psi
^{x}|_{V^{x}}^{-1}\circ\psi^{x}\circ\phi_{i}|_{U_{i}}^{-1}\right)  \left(
y\right)  \right\vert \\
&  \leq c_{1}\left(  C,\gamma\right)  \left\vert \left\vert f\right\vert
\right\vert _{2,Lip-\gamma^{\prime}}\left\vert x-y\right\vert ^{\gamma
^{\prime}-k}.
\end{align*}
Using this together with (\ref{uniform bound on derivatives}) in
(\ref{kth derivative}) yields%
\[
\left\vert D^{k}\left(  f\circ\phi_{i}|_{U_{i}}^{-1}\right)  \right\vert
_{\left(  \gamma^{\prime}-k\right)  -\text{H\"{o}l}}\leq\left(  1+2\left(
\frac{2C}{\delta_{2}}\right)  ^{\gamma^{\prime}-k}\right)  c_{1}\left(
C,\gamma\right)  \left\vert \left\vert f\right\vert \right\vert _{2,Lip-\gamma
^{\prime}}%
\]
and hence that $f\circ\phi_{i}|_{U_{i}}^{-1}$ is Lip-$\gamma^{\prime}$ with
\[
\left\vert \left\vert f\circ\phi_{i}|_{U_{i}}^{-1}\right\vert \right\vert
_{Lip-\gamma^{\prime}}\leq d\left(  C,\gamma,\delta_{2}\right)  \left\vert
\left\vert f\right\vert \right\vert _{2,Lip-\gamma^{\prime}}.
\]
These bounds are uniform with respect to indexing set $I$ so we take the
supremum over $i$ in $i$ to deduce the result
\[
\left\vert \left\vert f\right\vert \right\vert _{1,Lip-\gamma^{\prime}}\leq
d\left(  C,\gamma,\delta_{2}\right)  \left\vert \left\vert f\right\vert
\right\vert _{2,Lip-\gamma^{\prime}}.
\]
The first half of inequality in (\ref{equivalence}) can be shown in the same
way just by re-running the proof and interchanging the roles of the two atlases.
\end{proof}

In the definition of a Lip-$\gamma$ function $f$ \ we require $f\circ\phi
|_{U}^{-1}$ to be Lip-$\gamma$ for all charts $\left(  \phi,U\right)  .$
Clearly to verify if a function on a smooth manifold is smooth it is
sufficient to check that the functions $f\circ\phi|_{\phi\left(  U\right)
}^{-1}$ are smooth for any cover by charts. In our setting coordinate
transformations amplify the Lipschitz constants; nonetheless, we are able to
obtain an analogous result and involving the constants of the Lip-$\gamma$
atlas. Elucidating slightly, the following lemma shows that to check whether a
function is Lip-$\gamma$ on a given Lipschitz atlas it suffices to check that
$f\circ\phi|_{U}^{-1}$ is Lip-$\gamma$ for a sub-collection of charts in the
atlas satisfying a certain covering property. The key observation of course,
is that if this sub-collection is chosen appropriately then it constitutes a
Lip-$\gamma$ atlas for $M\ $and is equivalent to the original atlas.

\begin{lemma}
Let $M$ be a Lip-$\gamma$ manifold with atlas $\mathcal{A=}\left\{  \left(
\phi_{i},U_{i}\right)  :i\in I\right\}  $ and constants $\delta$ and $L.$
Suppose $f$ is a Banach space valued function on $M$ and that $I^{\prime}$ is
a subset of the indexing set $I$ such that the collection $\left\{
U_{i}^{\delta}:i\in I^{\prime}\right\}  $ $\ $forms a cover for $M.$ If there
is a finite constant $C^{\prime}>0$ such that
\[
\sup_{i\in I^{\prime}}\left\Vert f\circ\phi_{i}|_{U_{i}}^{-1}\right\Vert
_{Lip-\gamma}\leq C^{\prime}\text{ }%
\]
then $f$ is a Lip-$\gamma$ function on $M$ (with respect to $\mathcal{A}$) and
$\left\Vert f\right\Vert _{Lip-\gamma}\leq C,$ where $C$ is a constant that
only depends on the constants $\delta$ and $L$ of the Lip-$\gamma$ atlas and
the constant $C^{\prime}.$
\end{lemma}

\begin{proof}
It is easy to see that $\mathcal{A}^{^{\prime}}\mathcal{=}\left\{  \left(
\phi_{i},U_{i}\right)  :i\in I^{\prime}\right\}  $ is a $Lip-\gamma$ atlas for
$M$ with constants $\delta$ and $L.$ Moreover, it follows from condition 4 of
\ref{l-gamma-atlas} that $\mathcal{A}$ and $\mathcal{A}^{^{\prime}}$ are
equivalent with constant $L$. The conclusion is now immediate from Lemma
\ref{equiv}.
\end{proof}

\begin{definition}
Let $\alpha$ be a one-form on a Lip-$\gamma_{0}$ manifold $M$ taking values in
a Banach space $V.$ We say $\alpha$ is Lip-$\gamma$ if the one-form
\[
\left(  \phi_{i}|_{U_{i}}^{-1}\right)  ^{\ast}\alpha:\phi_{i}\left(
U_{i}\right)  \rightarrow L\left(
%TCIMACRO{\U{211d} }%
%BeginExpansion
\mathbb{R}
%EndExpansion
^{d},W\right)
\]
is Lip-$\gamma$ for any chart $\left(  \phi_{i},U_{i}\right)  $ in the atlas
defining the manifold and for some finite constant $C$
\begin{equation}
\sup_{i\in I}\left\vert \left\vert \left(  \phi_{i}|_{U_{i}}^{-1}\right)
^{\ast}\alpha\right\vert \right\vert _{Lip-\gamma}\leq C.
\label{uniform chart bound}%
\end{equation}
The Lipschitz constant is the minimal positive number $C$ that bounds
(\ref{uniform chart bound})$.$
\end{definition}

The following lemma shows the preceding definition is consistent with the
definition made for Banach space valued functions.

\begin{lemma}
\label{pullback-oneform}Let $M$ be a Lip-$\gamma$ manifold and $W$ a Banach
space, $h:M\rightarrow W$ a Lip-$\gamma$ function from $M$ to $W$ and
$\alpha:E\rightarrow L\left(  W,V\right)  $ a Lip-$\left(  \gamma-1\right)  $
one form on a subset $E$ of $W$ taking values in another Banach space $V.$
Then, if the range of $h$ is contained in $E$ (i.e. $h\left(  M\right)
\subseteq E$) the pullback $h^{\ast}\alpha$ is a $V-$valued Lip-$\left(
\gamma-1\right)  $ one form on $M$. Moreover, we have
\[
\left\Vert h^{\ast}\alpha\right\Vert _{Lip-\left(  \gamma-1\right)  }\leq
C\left(  \gamma,d\right)  \left\Vert \alpha\right\Vert _{Lip-\left(
\gamma-1\right)  }\left\Vert h\right\Vert _{Lip-\gamma}\max\left(  \left\Vert
h\right\Vert _{Lip-\gamma}^{\left\lfloor \gamma\right\rfloor },1\right)  ,
\]
where $C$ is a constant independent of $h$ and $\alpha.$
\end{lemma}

\begin{proof}
Let $\left(  \phi,U\right)  $ be any chart then we need to verify that the
pullback $\left(  \phi|_{U}^{-1}\right)  ^{\ast}\left(  h^{\ast}\alpha\right)
=\left(  h\circ\phi^{-1}\right)  ^{\ast}\alpha:\phi\left(  U\right)
\rightarrow L\left(
%TCIMACRO{\U{211d} }%
%BeginExpansion
\mathbb{R}
%EndExpansion
^{d},V\right)  $ is Lip-$\left(  \gamma-1\right)  $ and that we can bound the
Lip-$\left(  \gamma-1\right)  $ norm uniformly over the charts. Recall that
the pullback $\left(  h\circ\phi|_{U}^{-1}\right)  ^{\ast}\alpha$ is defined
by
\[
\left[  \left(  h\circ\phi|_{U}^{-1}\right)  ^{\ast}\right]  \alpha\left(
x\right)  \left(  v_{x}\right)  =\alpha\circ\left(  h\circ\phi|_{U}%
^{-1}\right)  \left(  x\right)  \left(  h\circ\phi|_{U}^{-1}\right)  _{\ast
}\left(  v_{x}\right)
\]
and so can be expressed in terms of the two functions
\begin{align*}
f_{1}  &  :=\alpha\circ\left(  h\circ\phi|_{U}^{-1}\right)  :\phi\left(
U\right)  \rightarrow L\left(  W,V\right) \\
f_{2}  &  :=\left(  h\circ\phi|_{U}^{-1}\right)  _{\ast}:\phi\left(  U\right)
\rightarrow L\left(
%TCIMACRO{\U{211d} }%
%BeginExpansion
\mathbb{R}
%EndExpansion
^{d},W\right)
\end{align*}
by $\left(  h\circ\phi|_{U}^{-1}\right)  ^{\ast}\alpha\left(  x\right)
=f_{1}\left(  x\right)  \circ f_{2}\left(  x\right)  .$ Because $h$ is
Lip-$\gamma$ it follows by definition that $h\circ\phi|_{U}^{-1}:\phi\left(
U\right)  \rightarrow W$ is Lip-$\gamma$ hence the derivative $f_{2}=\left(
h\circ\phi|_{U}^{-1}\right)  _{\ast}$ is Lip-$\left(  \gamma-1\right)  $
too$;$ furthermore, since $\alpha$ is Lip-$\left(  \gamma-1\right)  $ Lemma
\ref{composition of lip functions} implies that $f_{1}$ is Lip-$\left(
\gamma-1\right)  $ since it is the composition of two Lip-$\left(
\gamma-1\right)  $ functions. Given these observations it is elementary to
verify that
\[
\left(  h\circ\phi|_{U}^{-1}\right)  ^{\ast}\alpha\left(  \cdot\right)
=f_{1}\left(  \cdot\right)  \circ f_{2}\left(  \cdot\right)  :\phi\left(
U\right)  \rightarrow L\left(
%TCIMACRO{\U{211d} }%
%BeginExpansion
\mathbb{R}
%EndExpansion
^{d},V\right)
\]
is Lip-$\left(  \gamma-1\right)  $ and that the norm can be bounded by
\begin{equation}
\left\vert \left\vert \left(  h\circ\phi|_{U}^{-1}\right)  ^{\ast}%
\alpha\right\vert \right\vert _{\text{Lip-}\left(  \gamma-1\right)  }\leq
C\left(  \gamma\right)  \left\vert \left\vert f_{1}\right\vert \right\vert
_{\text{Lip-}\left(  \gamma-1\right)  }\left\vert \left\vert f_{2}\right\vert
\right\vert _{\text{Lip-}\left(  \gamma-1\right)  }. \label{lip norm bound}%
\end{equation}
Definition \ref{lipschitz function manifold} implies that
\begin{align*}
\left\vert \left\vert f_{2}\right\vert \right\vert _{\text{Lip}-\left(
\gamma-1\right)  }  &  \leq C\left(  \gamma,d\right)  \left\vert \left\vert
\left(  h\circ\phi|_{U}^{-1}\right)  _{\ast}\right\vert \right\vert
_{\text{Lip}-\left(  \gamma-1\right)  }\\
&  \leq C\left(  \gamma,d\right)  \left\vert \left\vert h\circ\phi|_{U}%
^{-1}\right\vert \right\vert _{\text{Lip}-\gamma}\leq C\left(  \gamma
,d\right)  \left\vert \left\vert h\right\vert \right\vert _{\text{Lip}-\gamma
},
\end{align*}
moreover, Lemma \ref{composition of lip functions} \ yields the bound%
\begin{align*}
\left\vert \left\vert f_{1}\right\vert \right\vert _{\text{Lip}-\left(
\gamma-1\right)  }  &  \leq C\left(  \gamma\right)  \left\vert \left\vert
\alpha\right\vert \right\vert _{\text{Lip}-\left(  \gamma-1\right)  }%
\max\left(  \left\vert \left\vert h\circ\phi|_{U}^{-1}\right\vert \right\vert
_{\text{Lip}-\left(  \gamma-1\right)  }^{\left\lfloor \gamma-1\right\rfloor
},1\right) \\
&  \leq C\left(  \gamma\right)  \left\vert \left\vert \alpha\right\vert
\right\vert _{\text{Lip}-\left(  \gamma-1\right)  }\max\left(  \left\vert
\left\vert h\right\vert \right\vert _{\text{Lip}-\gamma}^{\left\lfloor
\gamma\right\rfloor },1\right)
\end{align*}
from which the desired conclusion follows.
\end{proof}

Recall that for a chart $\left(  U,\phi\right)  $ the set $U^{\delta}$ denotes
the pre image of the open subset of $\phi\left(  U\right)  $ consisting of
points with distance at most $\delta$ from the boundary.

\begin{lemma}
\label{bump-functions}Let $M$ be a Lip-$\gamma$ manifold, $\left(
\phi,U\right)  $ a chart and suppose $f$ is a Lip-$\gamma$ function supported
on $U^{\delta}.$ Then $f$ extends to a Lip-$\gamma$ function on $M$ that
vanishes outside $U^{\delta/2}.$ Moreover the Lipschitz constant of the
extension only depends on the Lipschitz constant of $f$, $\delta$ and $L$ in
the definition of the Lip-$\gamma$ atlas.
\end{lemma}

\begin{proof}
The function $g$ defined to be $f\circ\phi|_{U}^{-1}$ on $\phi(U^{\delta})$
and $0$ on $\phi\left(  U\right)  \setminus\phi(U^{\delta/2})$ is Lip-$\gamma$
on the union of the two sets (the two sets are separated by a distance
$\delta/2).$ Thus by the Stein-Whitney extension theorem (Theorem
\ref{whitney-extension}) $g$ extends to a Lip-$\gamma$ function on $%
%TCIMACRO{\U{211d} }%
%BeginExpansion
\mathbb{R}
%EndExpansion
^{d}.$ The required function is now obtained by considering $g\circ\phi$ on
$U$ extended to vanish outside $U.$ Note that for any chart $\left(
\psi,V\right)  $ intersecting $U$ the function $g\circ\psi|_{V\cap U}^{-1}$ is
Lip-$\gamma$ as
\[
g\circ\psi|_{V\cap U}^{-1}=g\circ\phi|_{U}^{-1}\circ\phi\circ\psi|_{V\cap
U}^{-1}%
\]
is the composition of two Lip-$\gamma$ functions is Lip-$\gamma$ and $g$
vanishes by construction outside the intersection.
\end{proof}

As in the theory of differentiable manifolds partitions of unity are an
indispensable tool that will allow us to extend Lip-$\gamma$ functions on a
manifold. The following proposition guarantees their existence.

\begin{lemma}
\label{partition-of unity}For a Lipschitz-$\gamma$ atlas $\left\{  (\phi
_{i},U_{i}):i\in I\right\}  $ there are functions $f_{i}:M\rightarrow%
%TCIMACRO{\U{211d} }%
%BeginExpansion
\mathbb{R}
%EndExpansion
$ such that each $f_{i}$ has support contained in $U_{i}$ and the following
conditions hold:
\end{lemma}

\begin{enumerate}
\item $f_{i}\geq0$

\item
\[
\sum_{i\in I}f_{i}=1
\]

\item Each $f_{i}$ is Lip-$\gamma.$
\end{enumerate}

\begin{proof}
By Lemma \ref{bump-functions} there exist Lip-$\gamma$ functions $c_{i}$ on
$M$ that are identically one on $U_{i}^{\delta/2}$ and vanish outside $U_{i}.$
Note that by definition of a Lip-$\gamma$ atlas the sets $U_{i}^{\delta/2}$
are a locally finite cover of $M.$ Let $J_{i}$ denote the (finite) set of
indices consisting of the indices of all sets $U_{j}$ intersecting $U_{i}.$
Note that $\sum_{j\in J_{i}}c_{j}$ is a Lip-$\gamma$ function and by
construction bounded below by one on $U_{i}.$ Let $g:%
%TCIMACRO{\U{211d} }%
%BeginExpansion
\mathbb{R}
%EndExpansion
\rightarrow%
%TCIMACRO{\U{211d} }%
%BeginExpansion
\mathbb{R}
%EndExpansion
$ be a smooth function with bounded derivatives of all orders that agrees with
the function $1/x$ on $[1/2,\infty).$ For any Lip-$\gamma$ function $h\geq1$
on $U_{i}$ the composition $g\circ h$ is Lip-$\gamma.$ We deduce that
\[
1/\sum_{j\in J_{i}}c_{j}%
\]
is Lip-$\gamma$ on $U_{i}$ and we may define the function $f_{i}$ by setting
\[
f_{i}=\frac{c_{i}}{\sum_{j\in J_{i}}c_{j}}%
\]
on $U_{i}$ and taking $f_{i}$ to be zero on the complement. As $c_{i}$
vanishes outside $U_{i}$ the function $f_{i}$ is Lip-$\gamma$ and by
construction $\sum f_{i}=1$ on $M.$
\end{proof}

\begin{lemma}
Let $g$ be a Lip-$\gamma$ function defined on a compact subset $K.$ Then $g$
has a Lip-$\gamma$ extension to $M$ such that the Lipschitz constant of the
extension only depends on $\left\Vert g\right\Vert _{Lip-\gamma}$ , the
constants $L,\delta$ and $R$ in the definition on the Lip-$\gamma$ atlas and
the compact set $K$.
\end{lemma}

\begin{proof}
Let $J$ be a finite set of indices such the charts $\left(  \phi_{i}%
,U_{i}\right)  _{i\in J}$ cover the set $K.$ Let $g_{i}=f_{i}\widetilde{g} $
on $U_{i},$ where the $f_{i}$ are a partition of unity for the Lip-$\gamma$
atlas and $\widetilde{g}$ a Lip-$\gamma$ extension of $g$ from $U_{i}\cap K$
to $U_{i}$ . As $f_{i}$ vanishes outside $U_{i}$ we may extend $g_{i}$ to $M$
by setting it to zero outside $U_{i}$ without increasing the Lipschitz
constant. The required function is now given by $\sum_{j\in J}g_{j}.$
\end{proof}

\section{Rough paths on manifolds}

\subsection{Introduction and motivation}

In the following we establish a notion of rough paths on Lip-$\gamma$
manifolds that is global and coordinate free. All our definitions are
consistent if the manifold is a finite dimensional normed vector space, in
particular we show that in this case there is a bijective correspondence
between classical rough paths augmented by a starting point and rough paths in
the manifold sense. In this context it is important to realise that on a
manifold we do not have the translation invariance of a Banach space, so a
rough path on a manifold necessarily comes with a specific starting point
$x\in M$. Conceptually we regard a $p-$rough path $Z$ on a manifold as a
non-linear functional mapping sufficiently regular one forms $\alpha$ taking
values in any Banach space $W$ to geometric $p-$rough paths on $W.$ This
functional corresponds to the rough integral of $Z$ and so can only be
expected to be linear at level one. In this spirit of the rough integral on
Banach spaces we impose two natural consistency conditions on the map. First
we require the integral map to be continuous in the $p-$variation topology.
Second the push forward of $Z$ under a smooth map to a Banach space $V$ must
behave like a classical rough path on $V.$

A classical $p-$rough path on a Banach space $V$ corresponds via its trace
(projections to level one) to a path of bounded $p-$variation on $V.$ We
introduce a notion of support for $Z$ and show that this support is a
continuous path with values in $M$. We furthermore demonstrate that for
sufficiently small times the support of $Z$ is contained in some coordinate
chart. This allows us to prove a theorem can rather loosely be described as a
non-linear analogy of the classical Riesz representation theorem: rough paths
on a manifold are precisely the push forwards of finitely many classical rough
paths from the coordinate charts.

The structure of the section is as follows: We first present the basic
definitions and prove the basic properties of rough paths manifold and, their
support, restrictions and concatenation. In doing so we defer the somewhat
technical proofs to two key theorems that characterise the support of the
rough path. We subsequently give these proofs in subsections \ref{sec-proof1}
and \ref{sec-proof2}.

\subsection{Definition and basic properties of rough paths on manifolds.}

Before we can attempt a definition of \ a rough path on a Lip-$\gamma$
manifold we verify that the pushforward of a classical rough path on a Banach
space is well behaved and satisfies a chain rule that will be central to our
definition in the manifold setting.

\begin{lemma}
\label{chain-rule}Let $\gamma>p\geq1$ , $M$ be a Banach space and suppose $X$
is a classical geometric $p-$rough path on $M$ with starting point $x$. Then,
if $\psi:M\rightarrow$ $W$ and $\alpha:W\rightarrow L\left(  W,V\right)  $ are
a Lip$-\gamma$ function and a Lip-$\left(  \gamma-1\right)  $ one form
respectively for two arbitrary Banach spaces $W$ and $V,$ we have that
\[
\int\psi^{\ast}\alpha\left(  X_{s}\right)  dX_{s}=\int\alpha\left(
Z_{s}\right)  dZ_{s},
\]
where
\[
Z=\int d\psi\left(  X_{s}\right)  dX_{s}%
\]
with starting point $\psi\left(  x\right)  .$
\end{lemma}

\begin{proof}
We observe that $\psi^{\ast}\alpha$ and $d\psi$ are Lip$-\left(
\gamma-1\right)  $ one forms on $M$. Let $\tilde{x}:\left[  0,T\right]
\rightarrow M$ be a path of finite variation, note that $\tilde{z}\equiv
\psi\left(  \tilde{x}\right)  $ still has finite variation by the regularity
assumptions on $\psi$. If we denote the lifts of these paths by $\tilde
{X}=S_{\lfloor p\rfloor}\left(  \tilde{x}\right)  $ and $\tilde{Z}=S_{\lfloor
p\rfloor}\left(  \tilde{z}\right)  $ then we claim that for any Lip-$\left(
\gamma-1\right)  $ one-form $f$ we have
\begin{equation}
S_{\lfloor p\rfloor}\left(  \int f\left(  \tilde{x}_{s}\right)  d\tilde{x}%
_{s}\right)  \equiv\int f\left(  \tilde{X}\right)  d\tilde{X}. \label{lift2}%
\end{equation}
To see this we recall that the the Lyons extension theorem tells us that
$S_{\lfloor p\rfloor}\left(  \int f\left(  \tilde{x}_{s}\right)  d\tilde
{x}_{s}\right)  $ is the unique multiplicative functional on $T^{_{\lfloor
p\rfloor}}\left(  V\right)  $ with finite $1-$variation such that the
projection $\pi_{1}S_{\lfloor p\rfloor}\left(  \int f\left(  \tilde{x}%
_{s}\right)  d\tilde{x}_{s}\right)  $ $=\int f\left(  \tilde{x}\right)
d\tilde{x}_{s}.$ Since the geometric $p$-rough path $\int f\left(  \tilde
{X}\right)  d\tilde{X}$ is also a multiplicative functional of finite
$1-$variation with this property the two must agree by uniqueness and
(\ref{lift2}) holds. We therefore have that
\begin{align}
\int\psi^{\ast}\alpha\left(  \tilde{X}_{s}\right)  d\tilde{X}_{s}  &
=S_{\lfloor p\rfloor}\left(  \int\psi^{\ast}\alpha\left(  \tilde{x}%
_{s}\right)  d\widetilde{x}_{s}\right) \nonumber\\
&  =S_{\lfloor p\rfloor}\left(  \int\alpha\left(  \psi\left(  \tilde{x}%
_{s}\right)  \right)  \psi_{\ast}\overset{\cdot}{\tilde{x}}_{s}ds\right)
\nonumber\\
&  =S_{\lfloor p\rfloor}\left(  \int\alpha\left(  \psi\left(  \tilde{x}%
_{s}\right)  \right)  d\psi\left(  \tilde{x}_{s}\right)  \right) \nonumber\\
&  =S_{\lfloor p\rfloor}\left(  \int\alpha\left(  \tilde{z}_{s}\right)
d\tilde{z}_{s}\right) \nonumber\\
&  =\int\alpha\left(  \tilde{Z}_{s}\right)  d\tilde{Z}_{s}.
\label{change of variable}%
\end{align}
Now suppose that $\left(  x_{n}\right)  _{n=1}^{\infty}$ is a sequence of
finite variation paths with $X_{n}:=$ $S_{\lfloor p\rfloor}\left(
x_{n}\right)  \rightarrow X$ in the rough path topology and give each of the
$X_{n}$s a common starting point $x.$ Then, if $Z_{n}=S_{\lfloor p\rfloor
}\left(  \psi\left(  x_{n}\right)  \right)  $ with starting point $\psi\left(
x\right)  $ the relation (\ref{change of variable}) gives on the one hand
\begin{equation}
\int\psi^{\ast}\alpha\left(  X_{n}\left(  s\right)  \right)  dX_{n}\left(
s\right)  =\int\alpha\left(  Z_{n}\left(  s\right)  \right)  dZ_{n}\left(
s\right)  , \label{approximation}%
\end{equation}
on the other hand, taking $\alpha$ to be the identity one-form in
(\ref{change of variable}) and using the continuity of the map $Z\rightarrow
\int d\psi\left(  Z\right)  dZ$ in the rough path topology shows that
\[
Z_{n}=S_{\lfloor p\rfloor}\left(  \int d\psi\left(  x_{n}\right)
dx_{n}\right)  =\int d\psi\left(  Z_{n}\left(  s\right)  \right)
dZ_{n}\left(  s\right)  \rightarrow\int d\psi\left(  X_{s}\right)  dX_{s}=Z.
\]
Letting $n$ tend to infinity in (\ref{approximation}) and making a final use
of continuity concludes the proof.
\end{proof}

We first fix $\gamma_{0}>1$ and let $M$ be a Lip-$\gamma_{0}$ manifold. On a
manifold with Lip-$\gamma_{0}$ regularity we will define $p$-rough paths any
$p<\gamma_{0}.$ As in the classical Banach space case we expect such paths to
integrate Lip-($\gamma-1)$ one forms if $\gamma>p.$ Intuitively a rough path
on a manifold is an object that has a starting point and which consistently
integrates sufficiently regular compactly supported one forms taking values in
any Banach space.

\begin{definition}
\label{rp-manifold-def}Let $\gamma_{0}>p\geq1$ and $\gamma>p.$ A geometric
$p$-rough path over $\left[  0,T\right]  $ on a Lip-$\gamma_{0}$ manifold $M$
starting at a point $x\in M\ $is a map $Z$ which, for any Banach space $V,$
maps every $V$valued Lip-$\gamma$ one form $\alpha$ on $M $ to a continuous
$V-$valued geometric $p$-rough path over $\left[  0,T\right]  $ (in the
classical sense) such that :

\begin{enumerate}
\item Given any compactly supported Lip-$\gamma$ function $\psi:M\rightarrow
W$ mapping $M$ to a Banach space $W$ and any Banach space valued Lip-$\left(
\gamma-1\right)  $ one form $\alpha$ on $W$ we have
\[
Z(\psi^{\ast}\alpha)=\int\alpha(Y_{t})dY_{t},
\]
where $Y$ is the rough path starting at $\psi(x)$ with increments given by
$Z(d\psi)$ $.$

\item There exists a finite control $\omega$ such that for any $V-$ valued
Lip$-\left(  \gamma-1\right)  $ one form $\alpha$ on $M$ the geometric
$p-$rough path $Z\left(  \alpha\right)  $ is controlled by $\left\Vert
\alpha\right\Vert _{\text{Lip-}\left(  \gamma-1\right)  }\omega(s,t),$ i.e. we
have:
\[
\left\Vert Z\left(  \alpha\right)  ^{i}\right\Vert \leq\frac{\left(
\left\Vert \alpha\right\Vert _{\text{Lip-}\left(  \gamma-1\right)  }%
\omega(s,t)\right)  ^{i/p}}{\beta\left(  \frac{i}{p}\right)  !},\text{
}i=1,\ldots,\lfloor p\rfloor,~\left(  s,t\right)  \in\Delta_{T}.
\]

\end{enumerate}
\end{definition}

As we only define geometric rough paths on manifolds we will in the following
frequently drop the qualifier "geometric" and just refer to these objects as
rough paths. The map $Z$ can be thought of as functional on sufficiently
regular Banach space valued one forms, taking values in the geometric
$p-$rough paths. Condition one intuitively forces the pushforward of the rough
path onto a Banach space to behave like a classical rough path. The second
condition captures the continuity of the integral in $p-$variation another
feature that is expected of a classical rough path.

The behaviour of a classical rough path on a Banach space is determined by its
integral again the compactly supported one forms. To obtain a bijective
correspondence with the classical rough paths on finite dimensional normed
vector spaces we define an equivalence relation $\sim$ on the \thinspace$p-
$rough paths on a Lip-$\gamma_{0}$ manifold $M$ by saying two rough paths
$Z,Z^{\prime}$ satisfy $Z\sim Z^{\prime}$ if their starting points agree and
\begin{equation}
Z\left(  \alpha\right)  =Z^{\prime}\left(  \alpha\right)
\end{equation}
for any compactly supported Banach space valued Lip-$\left(  \gamma-1\right)
$ one form on $M.$ Before we show the consistency of Definition
\ref{rp-manifold-def} with the definition of classical rough paths on finite
dimensional normed vector space we require two elementary lemmas for preparation.

\begin{lemma}
\label{fix}Let $g:M\rightarrow V$ be a compactly supported Lip-$\gamma$
function and $Z$ a rough path on Lip-$\gamma$ manifold $M$. Then
\[
\text{supp}\left(  Z\left(  dg\right)  \right)  \subseteq\overline{g\left(
M\right)  }.
\]

\end{lemma}

\begin{proof}
As $dg=g^{\ast}$Id$_{V}$ we see that%
\[
Z\left(  dg\right)  =Z\left(  g^{\ast}\text{Id}_{V}\right)  =\int
d\text{Id}_{V}dY_{t},
\]
where Id$_{V}$ is the identity map on $V$ and $Y_{t}$ the classical rough path
with starting point $g\left(  x\right)  $ and increments $Z\left(  dg\right)
.$ If $Y$ had support off $\overline{g\left(  M\right)  }$ we could change the
one form on the right hand side off $\overline{g\left(  M\right)  }$ but on
supp$\left(  Y_{t}\right)  $ without affecting the left hand side, leading to
a contradiction.
\end{proof}

\begin{lemma}
\label{prep-2}Let $V$ be a finite dimensional normed vector space. There
exists $C>0$ such that for all $u\geq1,$ $x\in V$ there exist compactly
supported Lip-$\gamma$ functions $f_{u}:V\rightarrow V$ such that
\[
f_{u}|_{B\left(  x,u\right)  }=\text{Id}_{V}|_{B\left(  x,u\right)  }%
\]
and the one forms $df_{u}$ are Lip-$\gamma$ with Lipschitz constant at most
$C.$
\end{lemma}

\begin{proof}
To see the existence of the functions $f_{u}$ one may start with $f_{1}$, a
compactly supported smooth extension of the identity function on the unit ball
and define $f_{u}\left(  x\right)  $ for $x\in V$ by $f_{u}\left(  x\right)
=uf_{1}\left(  x/u\right)  $ and observe that the derivatives of $f_{u}$ are
bounded by the derivatives of $f_{1}.$
\end{proof}

Clearly only the one forms $df_{u}$ can be expected to have a uniformly
bounded Lipschitz constant in $u$ as only $d$Id$_{V,}$ but not Id$_{V}$ itself
is Lip-$\gamma$ on the whole of $V.$

The following theorem constructs on finite dimensional normed vector spaces an
explicit bijection between the quotient of the manifold rough paths under
$\sim$ and classical geometric rough paths with a starting point. The proof
that the map is onto is somewhat complex but instructive as many of the
techniques and strategies introduced there will continue to be used when we
proof our main results on the localisation of rough paths on manifolds. In
particular we express compactly supported one forms $\alpha$ as the pullback
of another one form $\widehat{\alpha}$ under a compactly supported
Lip-$\gamma$ function $g$. The first condition in the definition of a rough
path allows us then to re-express $Z\left(  \alpha\right)  $ as a classical
rough integral of the one form $\widehat{\alpha}$ against classical path with
increments given by $Z\left(  dg\right)  .$

\begin{theorem}
\label{bijection-bspace}Let $V$ be a finite dimensional normed vector space
with the Lip-$\gamma$ atlas of Example \ref{banach-manifold}. Then there
exists a bijection between classical geometric rough paths with a starting
point specified in $V$ and the quotient of the rough paths in the manifold
sense by the equivalence relation $\sim.$
\end{theorem}

\begin{proof}
Let $X\in G\Omega_{p}\left(  V\right)  $ be a classical geometric $p$ -rough
path on $V$ and $x\in V$ a starting point$.$ Note that there are finitely many
open sets $U_{i}$ (corresponding to charts in the Lip-$\gamma$ atlas on $V)$
that cover supp$\left(  X\right)  .$ If $\alpha$ is a Lip-$\left(
\gamma-1\right)  $ one form in the manifold sense it is immediate that
$\alpha$ is Lip-$\left(  \gamma-1\right)  $ in the classical sense on each set
$U_{i}$ with uniformly bounded Lipschitz constants$.$ Therefore the classical
rough integral $\int\alpha\left(  X\right)  dX$ exists and the continuity
theorem for the classical rough integral tells us that there exists a control
$\omega$ such that $\int\alpha\left(  X\right)  dX$ is controlled by
$C\omega,$ where $C$ is a constant depending on $\alpha$ only via its
Lip-$\gamma$ norm. Note that for this we have used that the collection
$\phi_{i}\left(  U_{i}\right)  $ covering supp$\left(  X\right)  $ is finite.
Thus we may set for any vector valued Lip-$\left(  \gamma-1\right)  $ one form
$\alpha$
\[
Z\left(  \alpha\right)  =\int\alpha\left(  X\right)  dX
\]
and $Z$ satisfies condition 2. of a rough path on a manifold$.$ By Lemma
\ref{chain-rule} $Z$ with starting point $x$ also satisfies condition 1. and
hence, $Z$ is a rough path in the manifold sense. We define a map $h$ from
classical geometric rough paths to manifold rough paths by setting $h\left(
X,x\right)  =Z$ for any $X\in G\Omega_{p}\left(  V\right)  $ with starting
point $x\in V.$ We confirm that that the map $h$ is injective by considering
\[
h\left(  X,x\right)  \left(  d\text{Id}_{V}\right)  =\int d\text{Id}%
_{V}\left(  X\right)  dX.
\]
To finish the proof we show the map $h$ is onto. Let $Z$ be a rough path in
the manifold sense on $V$ with starting point $x$, then to prove this claim we
have to demonstrate that there exists a classical rough path $X_{t}$ such that
$h\left(  X,x\right)  =Z;$ i.e. we have to show that
\begin{equation}
Z\left(  \alpha_{1}\right)  =\int\alpha_{1}\left(  X_{t}\right)  dX_{t}
\label{bij-id1}%
\end{equation}
for any compactly supported Lip$\left(  \gamma-1\right)  $ one form
$\alpha_{1}.$

Let $f_{u}$ be the Lip-$\gamma$ extensions of the identity obtained in Lemma
\ref{prep-2} and $c_{u}:V\rightarrow%
%TCIMACRO{\U{211d} }%
%BeginExpansion
\mathbb{R}
%EndExpansion
$ be a family of Lip-$\gamma$ bump functions with Lip-$\gamma$ norm
independent of $u$ that are identically one on $B\left(  x,u-1/2\right)  $ and
have support contained in $B\left(  x,u\right)  .$ Then the functions
$g_{u}:V\rightarrow V\oplus%
%TCIMACRO{\U{211d} }%
%BeginExpansion
\mathbb{R}
%EndExpansion
$ defined by $g_{u}=\left(  f_{u},c_{u}\right)  $ are Lip-$\gamma$ with
Lipschitz constants uniformly bounded in $u.$ We deduce by the continuity
condition 2. in Definition $\ref{rp-manifold-def}$ that there exists
$u^{\prime}>0$ such that
\begin{equation}
\text{supp}\left(  Z\left(  dg_{u}\right)  \right)  \subseteq B\left(
g_{u}\left(  x\right)  ,u^{\prime}-1\right)  \label{bij-id0}%
\end{equation}
for all $u\geq1.$

Suppose $\alpha_{1}$ is supported in the ball $B(x,r-1)$ for some
$r>u^{\prime}+1.$ Let $W:=V\oplus%
%TCIMACRO{\U{211d} }%
%BeginExpansion
\mathbb{R}
%EndExpansion
$ and define projections $\pi_{%
%TCIMACRO{\U{211d} }%
%BeginExpansion
\mathbb{R}
%EndExpansion
}$ and $\pi_{V}$ onto $%
%TCIMACRO{\U{211d} }%
%BeginExpansion
\mathbb{R}
%EndExpansion
$ and $V$ respectively. We define a one form $\widehat{\alpha_{1}}$ on by
letting $\widehat{\alpha_{1}}=x_{%
%TCIMACRO{\U{211d} }%
%BeginExpansion
\mathbb{R}
%EndExpansion
}\pi_{V}^{\ast}\alpha,$ where $x_{%
%TCIMACRO{\U{211d} }%
%BeginExpansion
\mathbb{R}
%EndExpansion
}$ is a Lip-$\gamma$ extension of coordinate function $\pi_{%
%TCIMACRO{\U{211d} }%
%BeginExpansion
\mathbb{R}
%EndExpansion
}$ from the subspace $V\times\lbrack0,1]$ to $W.$ \ Clearly $\widehat{\alpha
_{1}}$ is Lip-$\gamma$ and by construction of $g$ and $\widehat{\alpha_{1}}$
we have $\alpha=c_{r}\alpha=g_{r}^{\ast}\widehat{\alpha_{1}}.$ Hence, it
follows that%
\begin{equation}
Z\left(  \alpha\right)  =Z\left(  c_{r}\alpha\right)  =Z\left(  g_{r}^{\ast
}\widehat{\alpha_{1}}\right)  =\int\widehat{\alpha_{1}}Z\left(  dg_{r}\right)
. \label{bij-id2}%
\end{equation}
Let $\alpha_{2}$ be a Lip-$\left(  \gamma-1\right)  $ extension of $\alpha
_{1}|_{B\left(  x,u^{\prime}-1\right)  }$ that vanishes off $B\left(
x,u^{\prime}\right)  $ and $\widehat{\alpha_{2}}$ be defined by
$\widehat{\alpha_{2}}=x_{%
%TCIMACRO{\U{211d} }%
%BeginExpansion
\mathbb{R}
%EndExpansion
}\pi_{V}^{\ast}\alpha_{2}.$ We note that $g_{r}^{\ast}\widehat{\alpha_{2}%
}=\alpha_{2}$ and by identity $\left(  \ref{bij-id0}\right)  $ the rough path
$Z\left(  dg_{r}\right)  $ has support contained inside the ball $B\left(
g\left(  x\right)  ,u^{\prime}-1\right)  .$ Hence, as both one forms agree on
the support of the path $Z\left(  dg_{r}\right)  ,$ we deduce that
\begin{equation}
\int\widehat{\alpha_{1}}Z\left(  dg_{r}\right)  =\int\widehat{\alpha_{2}%
}Z\left(  dg_{r}\right)  =Z\left(  g_{r}^{\ast}\widehat{\alpha_{2}}\right)
=Z\left(  \alpha_{2}\right)  . \label{bij-id3}%
\end{equation}
With the same arguments used to derive $\left(  \ref{bij-id2}\right)  $ (note
the smaller support of the one form $\alpha_{2}$) it is straightforward to
show that
\begin{equation}
Z\left(  \alpha_{2}\right)  =\int\widehat{\alpha_{2}}Z\left(  dg_{u^{\prime}%
}\right)  . \label{bij-id4}%
\end{equation}
Using that $\alpha_{1}$ and $\alpha_{2}$ agree on an open subset containing
the support of the classical rough path $Z\left(  dg_{u^{\prime}}\right)  $
based at $g(x)$ we conclude that
\begin{equation}
\int\widehat{\alpha_{2}}Z\left(  dg_{u^{\prime}}\right)  =\int\widehat{\alpha
_{1}}Z\left(  dg_{u^{\prime}}\right)  . \label{bij-id5}%
\end{equation}
From left to right the chain of equalities $\left(  \ref{bij-id2}\right)  $ to
$\left(  \ref{bij-id5}\right)  $ yield%
\begin{equation}
Z\left(  \alpha_{1}\right)  =\int\widehat{\alpha_{1}}Z\left(  dg_{u^{\prime}%
}\right)  . \label{bij-id6}%
\end{equation}
Let $V_{a,b}$ denote the affine subspace $\left\{  w\in W:a\leq\pi_{%
%TCIMACRO{\U{211d} }%
%BeginExpansion
\mathbb{R}
%EndExpansion
}w\leq b\right\}  .$ Combining $\left(  \ref{bij-id0}\right)  $ and Lemma
\ref{fix} it follows that
\begin{align}
\text{supp}Z\left(  dg_{u^{\prime}}\right)   &  \subseteq\overline
{g_{u^{\prime}}\left(  M\right)  }\cap B\left(  g_{u^{\prime}}\left(
x\right)  ,u^{\prime}-1\right) \label{bij-id7}\\
&  =\overline{g_{u^{\prime}}\left(  M\right)  }\cap B\left(  g_{u^{\prime}%
}\left(  x\right)  ,u^{\prime}-1\right)  \cap V_{0,1}.\nonumber
\end{align}
In fact, as the bump function $c_{u^{\prime}}$ only takes values in $\left(
0,1\right)  $ when $f_{u^{\prime}}$ is the identity on the closed annulus of
radii $u^{\prime}-1/2$ and $u^{\prime},$ centred at $g_{u^{\prime}}\left(
x\right)  $ we have
\begin{equation}
\left(  \overline{g_{u^{\prime}}\left(  M\right)  }\cap B\left(  g_{u^{\prime
}}\left(  x\right)  ,u^{\prime}-1\right)  \right)  \cap\left\{  w\in
W:0<\pi_{R}w<1\right\}  =\emptyset, \label{bij-id9}%
\end{equation}
and combining this identity with $\left(  \ref{bij-id7}\right)  $ we get%
\[
\text{supp}Z\left(  dg_{u^{\prime}}\right)  \subseteq V_{0,0}\cup V_{1,1}.
\]
But $g_{u^{\prime}}\left(  x\right)  \in V_{1\,,1}$ and the classical rough
path $Z\left(  dg_{u^{\prime}}\right)  $ is continuous. Therefore as $Z\left(
dg_{u^{\prime}}\right)  $ starts in $V_{1,1}$ it remains in $V_{1,1}$ and we
deduce that supp$Z\left(  dg_{u^{\prime}}\right)  $ is contained in the affine
subspace $V_{1,1}.$ Using that the function $x_{%
%TCIMACRO{\U{211d} }%
%BeginExpansion
\mathbb{R}
%EndExpansion
}$ is identically one on supp$Z\left(  dg_{u^{\prime}}\right)  \subseteq
V_{1,1}$ we deduce that%
\begin{align}
\int\widehat{\alpha_{1}}dZ\left(  dg_{u^{\prime}}\right)   &  =\int x_{%
%TCIMACRO{\U{211d} }%
%BeginExpansion
\mathbb{R}
%EndExpansion
}\pi_{V}^{\ast}\alpha_{1}Z\left(  dg_{u^{\prime}}\right) \nonumber\\
&  =\int\pi_{V}^{\ast}\alpha_{1}Z\left(  dg_{u^{\prime}}\right) \nonumber\\
&  =Z\left(  \left(  \pi_{V}\circ g_{u^{\prime}}\right)  ^{\ast}\alpha
_{1}\right) \nonumber\\
&  =Z\left(  f_{u^{\prime}}^{\ast}\alpha_{1}\right)  =\int\alpha_{1}Z\left(
df_{u^{\prime}}\right)  \label{bij-id8}%
\end{align}
Setting $X_{t}=Z\left(  df_{u^{\prime}}\right)  $ equality $\left(
\ref{bij-id1}\right)  $ follows from equations $\left(  \ref{bij-id6}\right)
$ and $\left(  \ref{bij-id8}\right)  .$
\end{proof}

A map $f$ between two Lip-$\gamma$ manifolds $M$ and $N$ that pulls
Lip-$\left(  \gamma-1\right)  $ one forms on $M$ back to Lip-$\left(
\gamma-1\right)  $ on $N$ allows us to consider pushforwards of rough paths
between manifolds. More precisely the property we care about is that there is
a constant $C_{f}$ independent of $\alpha$ such that%
\begin{equation}
\left\Vert f^{\ast}\alpha\right\Vert _{Lip-\left(  \gamma-1\right)  }\leq
C_{f}\left\Vert \alpha\right\Vert _{Lip-\left(  \gamma-1\right)  }
\label{l-gamma}%
\end{equation}
for any Banach space valued Lip-$\left(  \gamma-1\right)  $ one form $\alpha$
on $M.$ Note that by Lemma \ref{pullback-oneform} a Lip-$\gamma$ map from a
Lip-$\gamma$ manifold to a Banach space satisfies $\left(  \ref{l-gamma}%
\right)  .$ For maps between manifolds we have so far not been able to
identify a similar, natural condition. A natural way to develop paths from one
manifold into another will be to consider differential equations.

\begin{lemma}
\label{pushforward-rp}Let $\gamma_{0},\gamma>p\geq1,$ $Z$ be a $p$-rough path
on a Lip-$\gamma_{0}$ manifold $M$ and $g:M\rightarrow N$ between $M$ and a
Lip-$\gamma_{0}$ manifold $N$ satisfying $\left(  \ref{l-gamma}\right)  $.
Then $g$ induces a pushforward $g_{\ast}$ from $p-$rough paths on $M$ to
$p-$rough paths on $N$ given by%
\[
g_{\ast}Z\left(  \alpha\right)  =Z\left(  g^{\ast}\alpha\right)
\]
for any Lip-$\gamma$ one form $\alpha$ defined on $N$ and new starting point
$g\left(  x\right)  .$
\end{lemma}

\begin{proof}
We need to verify that $g_{\ast}Z$ is a rough path on $N$. Let $\alpha$ be a
Lip-($\gamma-1)$ one form on $N.$ By definition%
\[
g_{\ast}Z(\psi^{\ast}\alpha)=Z(g^{\ast}\psi^{\ast}\alpha)=Z\left(  \left(
\psi\circ g\right)  ^{\ast}\alpha\right)  =\int\alpha dZ(d\left(  \psi\circ
g\right)  ).
\]
Using the identity $g^{\ast}d\psi=d\left(  \psi\circ g\right)  $ we deduce
that
\[
\int\alpha dZ(d\left(  \psi\circ g\right)  )=\int\alpha dZ(g^{\ast}d\psi
)=\int\alpha dg_{\ast}Z(d\psi).
\]
By assumption we have $\left\Vert g^{\ast}\alpha\right\Vert _{Lip-\left(
\gamma-1\right)  }\leq C_{g}\left\Vert \alpha\right\Vert _{Lip-\left(
\gamma-1\right)  }.$
\end{proof}

\begin{lemma}
{Pushforwards are morphisms on rough paths.}
\end{lemma}

\begin{proof}
Let $M,N,S$ be Lip-$\gamma_{0}$ manifolds, $g:M\rightarrow N$, $h:N\rightarrow
S$ maps satisfying $\left(  \ref{l-gamma}\right)  $ and $Z$ a rough path on
$M$. Let $\alpha$ be a $W$-valued Lip-$\left(  \gamma-1\right)  $ one form on
$S$. We need to show that $h_{\ast}g_{\ast}Z=(hg)_{\ast}Z.$ To see this note
that
\[
Z\left(  (h\circ g)^{\ast}\alpha\right)  =Z\left(  g^{\ast}(h^{\ast}%
\alpha)\right)  ,
\]
which follows directly from
\[
(h\circ g)^{\ast}(\alpha)=g^{\ast}(h^{\ast}\alpha).
\]

\end{proof}

The consistency of the composition of pushforwards is crucial as it ensures
invariance under change of coordinates.

\begin{lemma}
Let $\gamma_{0},\gamma>p\geq1,$ $(\phi,U),(\psi,V)$ be charts on a
Lip-$\gamma_{0}$ manifold $M$ and let $Z$ be a $p-$rough path on $M$ such that
$Z\left(  \alpha\right)  =Z\left(  \beta\right)  $ for any two Lip-$\left(
\gamma-1\right)  $ one forms that agree on $U\cap V$. Let $\Phi$ and $\Psi$
denote the pushforwards of $Z$ under $\phi$ and $\psi$. Then if
$\widetilde{\psi\circ\phi|_{U}^{-1}}$ denotes any Lip-$\gamma_{0}$ extension
of $\psi\circ\phi|_{U}^{-1}$ to a function from $%
%TCIMACRO{\U{211d} }%
%BeginExpansion
\mathbb{R}
%EndExpansion
^{d}$ to $%
%TCIMACRO{\U{211d} }%
%BeginExpansion
\mathbb{R}
%EndExpansion
^{d}.$ Then
\[
\Psi=(\widetilde{\psi\circ\phi|_{U}^{-1}})_{\ast}(\Phi).
\]

\end{lemma}

\begin{proof}
Let $\alpha$ be a $W-$valued Lip-$(\gamma-1)$ one form on $%
%TCIMACRO{\U{211d} }%
%BeginExpansion
\mathbb{R}
%EndExpansion
^{d}$. Then on the one hand $\Psi(\alpha)=Z(\psi^{\ast}\alpha)$, and on the
other hand
\begin{equation}
\left(  (\widetilde{\psi\circ\phi|_{U}^{-1}})_{\ast}\Phi\right)  (\alpha
)=\Phi((\widetilde{\psi\circ\phi|_{U}^{-1}})^{\ast}\alpha)=Z(\phi^{\ast
}(\widetilde{\psi\circ\phi|_{U}^{-1}})^{\ast}\alpha)=Z(\psi^{\ast}\alpha),
\end{equation}
where the last equality follows from the fact that $\phi^{\ast}%
(\widetilde{\psi\circ\phi|_{U}^{-1}})^{\ast}\alpha=\psi^{\ast}\alpha$ on
$U\cap V.$
\end{proof}

\subsection{Restriction and concatenations of rough paths on manifolds}

To reconstruct the entire rough path on $M$ from its local pushforwards onto
the coordinate neighbourhoods we need a way of concatenating rough paths on
the manifold defined on consecutive intervals. For classical rough paths this
is elementary (as the following lemma shows) however, because a rough path on
the manifold comes equipped with a starting point a little more care is needed
when performing the concatenation. Let us first observe the following.

\begin{lemma}
\label{classical concat}Let $p\geq1$ and suppose $Z$ and $Y$ are two
(classical) geometric $p-$rough paths on a Banach space $V$ over the two
intervals $\left[  s,t\right]  $ and $\left[  t,u\right]  $ respectively. Then
the concatenation
\[
Z\ast Y:\Delta\left[  s,u\right]  :=\left\{  \left(  r,v\right)  :s\leq r\leq
v\leq u\right\}  \rightarrow T^{\left\lfloor p\right\rfloor }\left(  V\right)
\]
defined by
\[
\left(  Z\ast Y\right)  _{r,v}=\left\{
\begin{array}
[c]{cc}%
Z_{r,v} & \text{for }s\leq r\leq v\leq t\text{ \ \ \ \ \ \ \ \ \ \ \ \ }\\
Z_{r,t}\otimes Y_{t,v} & \text{for }s\leq r\leq t\text{ and }t\leq v\leq u\\
Y_{r,v} & \text{for }t\leq r\leq v\leq u\text{ \ \ \ \ \ \ \ \ \ \ \ \ }%
\end{array}
\right.
\]
is also a geometric $p-$rough paths on $V$ over $\left[  s,u\right]  .$
\end{lemma}

\begin{proof}
This follows from the observation that if $\left(  z_{n}\right)
_{n=1}^{\infty}$ and $\left(  y_{n}\right)  _{n=1}^{\infty}$ are two sequences
of finite-varitation paths such that
\[
d_{p}\left(  S_{\left\lfloor p\right\rfloor }\left(  z_{n}\right)  ,Z\right)
\rightarrow0\text{ and }d_{p}\left(  S_{\left\lfloor p\right\rfloor }\left(
y_{n}\right)  ,Y\right)  \rightarrow0\text{ as }n\rightarrow\infty
\]
then it is easy to establish that $d_{p}\left(  S_{\left\lfloor p\right\rfloor
}\left(  z_{n}\ast y_{n}\right)  ,Z\ast Y\right)  \rightarrow0$ as
$n\rightarrow\infty.$
\end{proof}

\begin{definition}
If $\gamma_{0}>p\geq1$ and $Z$ and $Y$ are two geometric $p$-rough paths on a
Lip$-\gamma_{0}$ manifold $M$ defined over two intervals $\left[  s,t\right]
$ and $\left[  t,u\right]  $ respectively. Then the concatenation $Z\ast Y$ is
defined to be the map which for any $\gamma>p$ takes any Lip-$\gamma$ one form
$\alpha$ on $M$ with values in any Banach space $V$ to the function
\[
\left(  Z\ast Y\right)  \left(  \alpha\right)  :\Delta\left[  s,u\right]
\rightarrow T^{\left\lfloor p\right\rfloor }\left(  V\right)
\]
defined by
\[
\left(  Z\ast Y\right)  \left(  \alpha\right)  _{r,v}=\left\{
\begin{array}
[c]{cc}%
Z\left(  \alpha\right)  _{r,v} & \text{for }s\leq r\leq v\leq t\text{
\ \ \ \ \ \ \ \ \ \ \ \ }\\
Z\left(  \alpha\right)  _{r,t}\otimes Y\left(  \alpha\right)  _{t,v} &
\text{for }s\leq r\leq t\text{ and }t\leq v\leq u\\
Y\left(  \alpha\right)  _{r,v} & \text{for }t\leq r\leq v\leq u\text{
\ \ \ \ \ \ \ \ \ \ \ \ }%
\end{array}
\right.  .
\]

\end{definition}

Suppose $X$ is in $G\Omega_{p}\left(  V\right)  $ and is defined over the
interval $\left[  0,T\right]  .$ If $\left[  s,t\right]  \subseteq\left[
0,T\right]  $ is a subinterval and we let $Y|_{\left[  s,t\right]  }$ and $Y$
denote respectively the geometric $p$-rough paths over $\left[  s,t\right]  $
and $\left[  0,T\right]  $ with starting points $y_{1}$ and $y_{2}.$ Then for
any Lip-$\left(  \gamma-1\right)  $ one form $\alpha$ we have the consistency
relation:
\[
\int\alpha\left(  Y|_{\left[  s,t\right]  }\right)  dY|_{\left[  s,t\right]
}=\left.  \left(  \int\alpha\left(  Y\right)  dY\right)  \right\vert _{\left[
s,t\right]  }%
\]
provided the starting points are consistent, i.e. provided
\begin{equation}
y_{1}+\pi_{1}Y_{0,s}=y_{2}. \label{classical consistency}%
\end{equation}
Furthermore, if $Y$ is in $G\Omega_{p}\left(  V\right)  $ defined over
$\left[  s,t\right]  $ with starting point $y_{1}$ and $\bar{Y}$ is in
$G\Omega_{p}\left(  V\right)  $ defined over $\left[  t,u\right]  $ with
starting point $y_{2}$ the same condition (\ref{classical consistency})
ensures that
\[
\left(  \int\alpha\left(  Y\right)  dY\right)  \ast\left(  \int\alpha\left(
\bar{Y}\right)  d\bar{Y}\right)  =\int\alpha\left(  Y\ast\bar{Y}\right)
d\left(  Y\ast\bar{Y}\right)
\]
for any Lip-$\left(  \gamma-1\right)  $ one form, when $Y\ast\bar{Y}$ is
endowed with the starting point $y_{1}.$ In other words, the condition
(\ref{classical consistency}) ensures that the two operations of concatenating
rough paths and taking integrals along one-forms commute. The following
definition attempts to capture the same notion as this for a rough path on a
manifold, once again the lack of translation-invariance of the underlying
space means that this property needs to be slightly reformulated.

\begin{definition}
If $\gamma_{0}>p\geq1$ and $Z$ and $Y$ are two geometric $p$-rough paths on a
Lip$-\gamma_{0}$ manifold $M$, defined over the intervals $\left[  s,t\right]
$ and $\left[  t,u\right]  $ and having starting points $x_{1}$ and $x_{2}$
respectively. Then we say that $Z$\textbf{\ has an end point consistent with
the starting point of }$Y$ if for every $\gamma>p$ and every compactly
supported Lip-$\gamma$ function $g:M\rightarrow V,$ taking values in any
Banach space $V,$ we have
\[
g\left(  x_{1}\right)  +\pi_{1}Z\left(  dg\right)  _{s,t}=g\left(
x_{2}\right)  .
\]

\end{definition}

\begin{lemma}
\label{concat is a grp}If $\gamma_{0}>p\geq1$ and $Z$ and $Y$ are two
geometric $p$-rough paths on a Lip$-\gamma_{0}$ manifold $M$, defined over
$\left[  s,t\right]  $ and $\left[  t,u\right]  $, with starting points
$x_{1}$ and $x_{2}$ respectively and such that $Z$ has an end point consistent
with the starting point of $Y.$ Then the concatenation $Z\ast Y$ of $Z$ and
$Y$ is a geometric $p-$rough path on $M$ at the starting point $x_{1}$.
\end{lemma}

\begin{proof}
We need to verify that $Z\ast Y$ satisfies the conditions of Definition
\ref{rp-manifold-def}. We first observe that for any Banach space $V$ and any
$V- $valued Lip-$\left(  \gamma-1\right)  $ one form $\alpha$ on $M,$ Lemma
\ref{classical concat} guarantees that $\left(  Z\ast Y\right)  \left(
\alpha\right)  $ is a geometric $p$-rough path on $V$.

Let us now assume that $\alpha$ is a $V-$valued Lip-$\left(  \gamma-1\right)
$ one form on another Banach space $W$ and that $g:M\rightarrow W$ is a
compactly supported Lip-$\gamma$ function$;$ to check the first condition of
Definition \ref{rp-manifold-def} we need to show that
\[
\left(  Z\ast Y\right)  \left(  g^{\ast}\alpha\right)  =\int\alpha\left(
X\right)  dX,
\]
where $X$ is the $p-$rough path on $W$ with increments $\left(  Z\ast
Y\right)  \left(  dg\right)  $ and starting point $g\left(  x_{1}\right)  .$
Consider $r\leq v$ in $\left[  s,u\right]  $ such that $s\leq r\leq t\leq
v\leq u$ (the other cases necessitating only slight adaptation of the
following argument), then using the multiplicative property of rough paths we
see that
\begin{equation}
\int_{r}^{v}\alpha\left(  X\right)  dX=\int_{r}^{t}\alpha\left(  X\right)
dX\otimes\int_{t}^{v}\alpha\left(  X\right)  dX. \label{multi}%
\end{equation}
Using the definition of the concatenation we see that
\[
\int_{r}^{t}\alpha\left(  X\right)  dX=\int_{r}^{t}\alpha\left(  X_{1}\right)
dX_{1}%
\]
where $X_{1}=Z\left(  dg\right)  $ over $\left[  s,t\right]  $ with starting
point $g\left(  x_{1}\right)  .$ Furthermore, the consistency between the end
points and starting points of $Z$ and $Y$ shows that
\[
g\left(  x_{1}\right)  +\pi_{1}\left(  Z\ast Y\right)  \left(  dg\right)
_{s,t}=g\left(  x_{1}\right)  +\pi_{1}Z\left(  dg\right)  _{s,t}=g\left(
x_{2}\right)
\]
and hence
\[
\int_{t}^{v}\alpha\left(  X\right)  dX=\int_{t}^{v}\alpha\left(  X_{2}\right)
dX_{2}%
\]
where $X_{2}=Y\left(  dg\right)  $ with starting point $g\left(  x_{2}\right)
.$ Putting everything together in (\ref{multi}) we may conclude:%
\begin{align*}
\int_{r}^{v}\alpha\left(  X\right)  dX  &  =\int_{r}^{t}\alpha\left(
X_{1}\right)  dX_{1}\otimes\int_{r}^{t}\alpha\left(  X_{2}\right)  dX_{2}\\
&  =Z\left(  g^{\ast}\alpha\right)  _{r,t}\otimes Y\left(  g^{\ast}%
\alpha\right)  _{t,v}\\
&  =\left(  Z\ast Y\right)  \left(  g^{\ast}\alpha\right)  .
\end{align*}

To check the second condition of Definition \ref{rp-manifold-def} we observe
the following elementary estimate for $s\leq r\leq t\leq v\leq u$ and
$i=1,...,\left\lfloor p\right\rfloor $
\begin{align*}
\left\vert \left\vert \left(  Z\ast Y\right)  \left(  \alpha\right)
_{r,v}^{i}\right\vert \right\vert  &  =\left\vert \left\vert \sum_{k=0}%
^{i}Z\left(  \alpha\right)  _{r,t}^{k}\otimes Y\left(  \alpha\right)
_{t,v}^{k-i}\right\vert \right\vert \\
&  \leq\sum_{k=0}^{i}\left\vert \left\vert Z\left(  \alpha\right)  _{r,t}%
^{k}\right\vert \right\vert \left\vert \left\vert Y\left(  \alpha\right)
_{t,v}^{k-i}\right\vert \right\vert \\
&  \leq C\left(  p\right)  \frac{\left(  \left\vert \left\vert \alpha
\right\vert \right\vert _{\text{Lip-}\left(  \gamma-1\right)  }\omega\left(
r,v\right)  \right)  ^{i/p}}{\beta\left(  \frac{i}{p}\right)  !}.
\end{align*}
Similar but more straight-forward calculations pertain to the cases $s\leq
r\leq v\leq t$ and $t\leq r\leq v\leq u.$
\end{proof}

\begin{lemma}
[Associativity of concatenation]\label{associativity of concat}If $\gamma
_{0}>p\geq1$ and $Z$ , $Y$ and $W$ are three geometric $p$-rough paths on a
Lip$-\gamma_{0}$ manifold $M$, defined over $\left[  s,t\right]  $, $\left[
t,u\right]  $ and $\left[  u,v\right]  $ with starting points $x_{1}$, $x_{2}$
and $x_{3}$ respectively. Suppose that $Z$ has an end point consistent with
the starting point of $Y$ and $Y$ has an end point consistent with the
starting point of $W.$ Then $\left(  Z\ast Y\right)  \ast W$ and $Z\ast\left(
Y\ast W\right)  $ are geometric $p$-rough paths over $\left[  s,v\right]  $
and we have
\[
\left(  Z\ast Y\right)  \ast W=Z\ast\left(  Y\ast W\right)  .
\]

\end{lemma}

\begin{proof}
It is self-evident that if $Z$ has an end point consistent with the starting
point of $Y$ then $Z$ has an end point consistent with the starting point of
the concatenation $Y\ast W.$ By Lemma \ref{concat is a grp} $Z\ast\left(
Y\ast W\right)  $ is a geometric $p$-rough path. In a similar way for any
compactly supported Lip-$\gamma$ function $g:M\rightarrow V$ \ we have that
\begin{align*}
g\left(  x_{1}\right)  +\pi_{1}\left(  Z\ast Y\right)  \left(  dg\right)
_{s,u}  &  =g\left(  x_{1}\right)  +\pi_{1}Z\left(  dg\right)  _{s,t}+\pi
_{1}Y\left(  dg\right)  _{t,u}\\
&  =g\left(  x_{2}\right)  +\pi_{1}Y\left(  dg\right)  _{t,u}\\
&  =g\left(  x_{3}\right)  .
\end{align*}
Hence, Lemma \ref{concat is a grp} implies that $\left(  Z\ast Y\right)  \ast
W $ is a geometric $p$-rough path. \ A simple exercise in unravelling the
definitions then shows that
\[
\left(  \left(  Z\ast Y\right)  \ast W\right)  \left(  \alpha\right)  =\left(
Z\ast\left(  Y\ast W\right)  \right)  \left(  \alpha\right)
\]
for any Banach space valued Lip-$\left(  \gamma-1\right)  $ one form $\alpha$
on $M,$ which completes the proof.
\end{proof}

\begin{remark}
Henceforth we can write $Z\ast Y\ast W$ to mean either $\left(  Z\ast
Y\right)  \ast W$ or $Z\ast\left(  Y\ast W\right)  $ free of ambiguity.
\end{remark}

\subsection{The support of a rough path on a manifold}

Built into the definition of a rough path on a manifold is a notion of an
non-linear, rough integral against sufficiently regular one forms. We show
that this integral map considered on a core of compactly supported one forms
eventually allows us to show that the non-linear functional that is the rough
path corresponds to a continuous path on the manifold. In a first step we
prove that for sufficiently small but uniform times $t_{0}$ a rough path on a
Lip-$\gamma$ manifold restricts to a chart. The pushforward of the restricted
path under the chart map corresponds (up to equivalence under $\sim)\ $to a
classical geometric rough path and we can identify its support with the pre
image of this classical rough path. Once this localisation is achieved it is
not difficult to demonstrate that the localised rough path has an endpoint and
we can define its restriction to an arbitrary time interval $\left[
s,t\right]  \subseteq\left[  0,t_{0}\right]  .$ Finally by considering the
concatenation of localised paths we can prove that rough paths on a manifold
are (up to equivalence under $\sim$ ) nothing but the pushforward of finitely
many classical geometric rough paths on the images of the chart maps.

\begin{definition}
If $\gamma_{0}>p\geq1$ we say that a geometric $p$-rough path $Z$ on a
Lip-$\gamma_{0}$ manifold $M$ with starting point $x$ in $M$ misses an open
set $U\subseteq M$ if $Z(\alpha)=0$ for any compactly supported Banach space
valued Lip-($\gamma-1)$ one form whose support is contained in $U.$
\end{definition}

\begin{definition}
\label{support}If $\gamma_{0}>p\geq1,$ $Z$ is a geometric $p$-rough on a
Lip-$\gamma_{0}$ manifold $M$ with starting point $x$ in $M,$ then we define
the support of $Z$ by
\[
\text{supp }Z=\left\{  x\right\}  \cup\left(  M\setminus\bigcup
_{U~\text{misses}~Z}U\right)  .
\]

\end{definition}

\begin{remark}
It is important to note that our notion of support only considers compactly
supported one forms.
\end{remark}

It is immediate from the definition that the support is a closed set. We have
already seen that every finite dimensional normed vector space $V$ may be
equipped with the structure of a Lip$-\gamma_{0}$ manifold; moreover, we have
shown a bijective correspondence between the set of classical rough paths on
$V$ together with a starting point and the quotient of the set of rough paths
on $V$ (where $V$ is regarded as a Lip-$\gamma_{0}$ manifold) by the
equivalence relation $\sim$ . A classical rough path $X$ in $V$ with starting
point $x$ already has a natural notion of support, namely the support of the
underlying base path $x_{t}:=x+\pi_{1}X_{0,t}$ and we need to first check that
these two concepts are the same.

\begin{proposition}
Let $X$ be a classical rough path on a finite dimensional normed vector space
$V$ with starting point $x_{0}\in V.$ Let $Z$ denote a rough path on $V $
(with the same starting point) corresponding to $X$ under the bijection of
Lemma \ref{bijection-bspace} when $V$ is viewed as a $Lip-\gamma$ manifold.
Then if $x:\left[  0,T\right]  \rightarrow V$ is the continuous path defined
by $x_{t}=x_{0}+\pi_{1}X_{0,s}$ we have%
\[
\text{supp }Z|_{\left[  0,T\right]  }=\left[  x\right]  :=\left\{  x_{t}%
:t\in\left[  0,T\right]  \right\}  .
\]

\end{proposition}

\begin{proof}
If $x$ is constant then the result is trivial because $\left[  x\right]
=\left\{  x_{0}\right\}  $ and $Z\left(  \alpha\right)  $ vanishes for every
one form $\alpha$ so that supp $Z|_{\left[  0,T\right]  }=\left\{
x_{0}\right\}  $. Suppose $x$ is not constant. Recall that $X$ and $Z$ are
related in the sense that the integral of $X_{s,t}$ against any compactly
supported Lip-$\left(  \gamma-1\right)  $ coincides with the integral of
$\alpha$ against $Z\left(  d\text{Id}_{B\left(  x,u\right)  }\right)  _{s,t}$
for all sufficiently large $u>0$ and where Id$_{B\left(  x,u\right)
}:V\rightarrow V$ is a suitable extension of the identity map from $B\left(
x,u\right)  $. Moreover arguing as in the proof of Lemma
\ref{bijection-bspace} we see that for any compactly supported $Lip-\left(
\gamma-1\right)  $ one form $\alpha$ on $V$%
\[
Z\left(  \alpha\right)  =\int\alpha\left(  X_{s}\right)  dX_{s}%
\]
where $X$ is considered together with its starting point $x_{0}\in V$. Thus,
if $\alpha$ is supported on $V\setminus\left[  x\right]  $ it follows from the
definition of the rough integral (see \cite{lyons-stflour}) that $Z\left(
\alpha\right)  =0$ and hence supp $Z|_{\left[  0,T\right]  }\subseteq\left[
x\right]  $.

Conversely, consider $w=x_{t}\in\left[  x\right]  $ and let $N$ be some
neighbourhood of $w$ in $V.$ Let $r>0$ be such that the closure of ball or
radius $r$ centred at $w$ is contained in $N,$ i.e.%
\[
\overline{B\left(  w,r\right)  }\subset N\text{.}%
\]
Because $x$ is continuous and non-constant there exists $s\neq t$ in $\left[
0,T\right]  $ such that between $s$ and $t$ the path $x$ stays in the ball
$B\left(  w,r\right)  $ while $x_{s}\neq x_{t}$. Suppose $s<t$ (the argument
for $t<s$ proceeds \emph{mutatis mutandis}) then by considering a $Lip-\left(
\gamma-1\right)  $ one form $\alpha$ which agrees with the one form $d$Id on
$B\left(  w,r\right)  $ and vanishes off $N$ we can deduce that
\[
\pi_{1}Z\left(  \alpha\right)  _{s,t}=\pi_{1}\int_{s}^{t}\alpha\left(
X_{u}\right)  dX_{u}=x_{t}-x_{s}\neq0\text{.}%
\]
Thus supp $Z|_{\left[  0,T\right]  }\supseteq\left[  x\right]  $ and the
result follows.
\end{proof}

The main of results of this subsection are following two theorems that
demonstrate that for sufficiently small times a rough path on a manifold
restricts to a chart and on this chart the support can be identified with a
continuous path.

\begin{theorem}
\label{one chart localisation}Let $\gamma_{0}>p\geq1$ and suppose $Z$ is a
geometric $p-$rough path on a Lip-$\gamma_{0}$ manifold $M$ with starting
point $x$ and is controlled by some control function $\omega.$ There exists a
strictly positive time $t_{0}$ independent of $x$ such that%
\[
\text{supp}(Z%
%TCIMACRO{\U{a6}}%
%BeginExpansion
\vert
%EndExpansion
_{[0,t_{0}]})\subseteq U
\]
for some chart $\left(  \phi,U\right)  $ containing $x.$
\end{theorem}

In the next step we identify the pushforward of the restricted rough path
$Z|_{\left[  0,t_{0}\right]  }$ under the chart map $\phi$ with its
corresponding classical rough path and use this relation to show that the
support of $Z|_{\left[  0,t_{0}\right]  }$ is a continuous path in the manifold.

\begin{theorem}
\label{push forward localisation}Suppose $\gamma_{0}>p\geq1.$ Let $Z$ be a
geometric $p$-rough path on a Lip$-\gamma_{0}$ manifold $M$ with starting
point $x\in M$ and $t_{0}$, $\left(  \phi,U\right)  $ as in Theorem
\ref{one chart localisation}. Then the support of the pushforward of $Z$ under
$\phi$ is contained in the coordinate neighbourhood $\phi\left(  U\right)  ,$
i.e.
\[
\text{supp }\left(  \phi_{\ast}Z\right)  |_{\left[  0,t_{0}\right]  }%
\subseteq\phi\left(  U\right)  .
\]
Moreover, we can recover the support of $Z$ restricted to $\left[
0,t_{0}\right]  $ as the preimage under $\phi$ of the support of the pushed
forward rough path; that is, we have
\[
\text{supp }Z|_{\left[  0,t_{0}\right]  }=\left(  \phi|_{U}\right)
^{-1}\left(  \text{supp }\left(  \phi_{\ast}Z\right)  |_{\left[
0,t_{0}\right]  }\right)  .
\]

\end{theorem}

As a corollary we obtain the following characterisation of the localised rough path.

\begin{corollary}
\label{last-in-line}Suppose $\gamma_{0}>p\geq1.$ Let $Z$ be a geometric $p
$-rough path on a Lip$-\gamma_{0}$ manifold $M$ with starting point $x\in M$
and $t_{0}$, $\left(  \phi,U\right)  $ as in Theorem
\ref{one chart localisation}. Then

\begin{enumerate}
\item $\left(  \phi_{\ast}Z\right)  |_{\left[  0,t_{0}\right]  }\left(
\xi\right)  =\left(  \phi_{\ast}Z\right)  |_{\left[  0,t_{0}\right]  }\left(
\eta\right)  $ for any compactly supported Lip-$\left(  \gamma-1\right)  $ one
forms which agree on $\phi\left(  U\right)  .$

\item The restriction $Z|_{\left[  0,t_{0}\right]  }$ is characterised by the
push forward $\left(  \phi_{\ast}Z\right)  |_{\left[  0,t_{0}\right]  }$ in
the sense that, for any compactly supported Banach space valued Lip$-\left(
\gamma-1\right)  $ one form $\alpha$ on $M,$ \
\[
Z|_{\left[  0,t_{0}\right]  }\left(  \alpha\right)  =\left(  \phi_{\ast
}Z\right)  |_{\left[  0,t_{0}\right]  }\left(  \xi_{\alpha}\right)  ,
\]
where $\xi_{\alpha}$ is any compactly supported Lip$-\left(  \gamma-1\right)
$ one form on $%
%TCIMACRO{\U{211d} }%
%BeginExpansion
\mathbb{R}
%EndExpansion
^{d}$ which agrees with $\left(  \phi|_{U}^{-1}\right)  ^{\ast}\alpha$ on
$\phi\left(  U\right)  .$
\end{enumerate}
\end{corollary}

The proofs of these theorems are slightly delicate. They will be carried out
in the following subsections which we will precede by proving a number of
small lemmas.

We emphasise once more that the value of $t_{0}$ in Theorem
\ref{one chart localisation} (and hence also in Theorem
\ref{push forward localisation}) does not depend on the starting point of the
rough path. This observation, which is made possible by the definition we took
for a Lipschitz manifold, allows us to decompose any rough path using the
rough paths obtained by pushing the restricted path forward under appropriate
charts. The uniformity of $t_{0}$ with respect to the starting point is
important because it ensures that we need only consider restrictions of $Z$
over a finite number of intervals. Because the locally pushed forward rough
paths are rough paths on a finite dimensional normed vector space (the
coordinate space is a subset of $%
%TCIMACRO{\U{211d} }%
%BeginExpansion
\mathbb{R}
%EndExpansion
^{d}$), they may be identified with classical rough paths on this space. By
showing that these different pieces have consistent starting and end points we
may concatenate the decomposed path to recover the original path. This has a
number of consequences: it allows us to construct a base path on the manifold
underlying our rough path; it implies that locally a rough path on $M$ is no
more difficult to understand than a classical rough path; and it implies that
the support of $Z$ is compact.

\bigskip We close this section by proving a small lemma which shows that the
pushforwards of rough paths under compactly supported maps preserve the support.

\begin{lemma}
\label{lemma-support-image}Suppose $\gamma_{0}>p\geq1.$ Let $Z$ be a geometric
$p$-rough path on a Lip-$\gamma_{0}$ manifold $M.$ Suppose $U$ is an open,
pre-compact subset of $M$ and $f$ is a compactly supported map from $M$ to a
Lip-$\gamma_{0}$ manifold $N$ satisfying $\left(  \ref{l-gamma}\right)  .$
Moreover suppose that supp$\left(  Z\right)  \subseteq U$, then%
\[
\text{supp}\left(  f_{\ast}Z\right)  \subseteq f(U).
\]

\end{lemma}

\begin{proof}
As supp$\left(  Z\right)  $ is compact we can find an open set $V$ such that
\[
\text{supp}\left(  Z\right)  \subseteq V\subseteq\overline{V}\subseteq U.
\]
To see this note that $M$ is completely metrisable and the continuous function
$d\left(  x\right)  $ defined by the distance of $x$ to $U^{c}$ for all $x\in
M$ is strictly positive on $U$ and therefore $\inf_{x\in\text{supp}\left(
Z\right)  }d(x)>0.$ Let $\alpha$ be a Lip-$\gamma-1$ one form on $N$ that
vanishes on $f\left(  V\right)  $. By definition $f_{\ast}Z\left(
\alpha\right)  =Z\left(  f^{\ast}\alpha\right)  $ and by construction
$f^{\ast}\alpha$ vanishes on $V.$ Therefore, as supp$\left(  Z\right)
\subseteq U$ we deduce $f_{\ast}Z\left(  \alpha\right)  =0.$ Hence, $f_{\ast
}Z\left(  \alpha\right)  $ misses the open set $f\left(  \overline{V}\right)
^{c}$ and by definition of the support we have supp$\left(  f_{\ast}Z\right)
\subseteq f\left(  \overline{V}\right)  \subseteq f\left(  U\right)  .$
\end{proof}

\subsection{Rough paths on manifolds are concatenations of localisable rough
paths.}

Suppose that we are given a geometric $p-$rough path $Z$ on a Lip-$\gamma_{0}$
manifold $M$ with some starting point $x\in M,$ $\gamma_{0}>p\geq1$ and
suppose $Z$ defined over the interval $\left[  0,T\right]  $. We now describe
how to construct a finite partition $D=\left\{  s_{i}:i=0,1,....,N\right\}  $
of $\left[  0,T\right]  $ such that $Z$ is the concatenation of rough paths
over the intervals $\left[  s_{i},s_{i+1}\right]  $ which are localisable in
the sense that they may be related to the pushforward (under the coordinate
map) of a restriction of $Z$. More precisely, we will provide an inductive
definition of a sequence $\left\{  \left(  Z^{n},z^{n},\phi^{n},U^{n}\right)
\right\}  _{n=1}^{N}$ whose elements have components consisting of the
following: a geometric $p-$rough $Z^{n}$ on $M\ $over the interval $\left[
s_{n-1},s_{n}\right]  ;$ a\ starting point $z^{n}$ associated with $Z^{n};$
and a chart $\left(  \phi^{n},U^{n}\right)  $ such that $Z^{n}$ and the
(restriction of the) pushforward $\left(  \phi_{\ast}^{n}Z\right)  |_{\left[
s_{n-1},s_{n}\right]  }$ are related by
\begin{equation}
Z^{n}\left(  \alpha\right)  =\left(  \phi_{\ast}^{n}Z\right)  |_{\left[
s_{n-1},s_{n}\right]  }\left(  \xi_{\alpha}\right)  ,
\label{localisation relation}%
\end{equation}
for any (compactly supported) Lip$-\left(  \gamma-1\right)  $ one form
$\xi_{\alpha}$ which agrees with $\left(  \phi^{n}|_{U^{n}}^{-1}\right)
^{\ast}\alpha$ on $\phi^{n}\left(  U^{n}\right)  .$ The idea is that the
original rough path $Z$ can be reconstructed by concatenating this finite
sequence; and because the rough paths in this sequence have rather straight
forward representation using the coordinate maps they are not much more
difficult to understand than classical rough paths. In this way (global)
properties of the original rough path can be understood in terms of the
properties of these localisations.

To carry out this procedure we first initialise by taking $t_{0}>0$ to be as
in Theorem \ref{one chart localisation}, $N=\left\lceil \frac{T}{t_{0}%
}\right\rceil $ \ and define $D=\left\{  s_{k}:k=0,1,...,N\right\}  $ by
setting $s_{k}=\left(  kt_{0}\right)  \wedge T.$ We then perform the following steps:

\begin{description}
\item[Base case] Define $z^{1}=x,$ $Z^{1}=Z|_{\left[  0,s_{1}\right]  }$ and
let $\left(  \phi^{1},U^{1}\right)  $ be any chart (the existence of which is
guaranteed by Corollary \ref{last-in-line}) so that $Z^{1}$ and $\left(
\phi_{\ast}^{1}Z\right)  |_{\left[  0,s_{1}\right]  }$ are related by
(\ref{localisation relation}).

\item[Induction step] Suppose $\left(  Z^{n},z^{n},\phi^{n},U^{n}\right)  $
have been constructed for $k=1,....,n.$ If $n=N$ then we stop. Otherwise, we
define
\[
z^{n+1}=\phi^{n}|_{U^{n}}^{-1}\left(  \phi^{n}\left(  z^{n}\right)  +\pi
_{1}Z\left(  d\phi^{n}\right)  _{s_{n-1},s_{n}}\right)
\]
and note $\left(  Z|_{\left[  s_{n},T\right]  }\right)  |_{\left[
s_{n},s_{n+1}\right]  }=Z|_{\left[  s_{n},s_{n+1}\right]  }.$ We define
$Z^{n+1}$ to be $Z|_{\left[  s_{n},s_{n+1}\right]  }$ with starting point
$z^{n+1},$ then Corollary \ref{last-in-line} guarantees the existence of a
chart $\left(  \phi^{n+1},U^{n+1}\right)  $ such that $Z^{n+1}$ and
\[
\left(  \phi_{\ast}^{n+1}Z|_{\left[  s_{n},s_{n+1}\right]  }\right)
|_{\left[  s_{n},s_{n+1}\right]  }=\left(  \phi_{\ast}^{n+1}Z\right)
|_{\left[  s_{n},s_{n+1}\right]  }%
\]
are related by (\ref{localisation relation}).
\end{description}

\begin{definition}
Let $Z$ be a geometric $p-$rough path on a Lip-$\gamma_{0}$ manifold $M$ with
starting point $x\in M,$ $\gamma_{0}>p\geq1$ defined over the interval
$\left[  0,T\right]  .$ Then any sequence $\left\{  \left(  Z^{n},z^{n}%
,\phi^{n},U^{n}\right)  \right\}  _{n=1}^{N}$ of be a geometric $p- $rough
paths, starting points and charts defined according to the above construction
relative to some partition $D$ of $\left[  0,T\right]  $ is a called a
localising sequence for $Z$.
\end{definition}

\begin{lemma}
\label{consistency of localisations}Let $Z$ be a geometric $p-$rough path on a
Lip-$\gamma_{0}$ manifold $M$ with starting point $x$ and having a localising
sequence $\left\{  \left(  Z^{n},z^{n},\phi^{n},U^{n}\right)  \right\}
_{n=1}^{N}.$ Then for any $k=1,....,N-1,$ $Z^{k}$\textbf{\ }has an end point
consistent with $z^{k+1},$ the starting point of $Z^{k+1}.$
\end{lemma}

\begin{proof}
Suppose that $g:M\rightarrow V$ is a compactly supported Lip-$\gamma$ function
in a Banach space $V.$ Let $f_{n}:%
%TCIMACRO{\U{211d} }%
%BeginExpansion
\mathbb{R}
%EndExpansion
^{d}\rightarrow V$ be any compactly supported Lip-$\gamma$ function which
agrees with $g\circ\left(  \phi^{n}|_{U^{n}}^{-1}\right)  $ on $\phi
^{n}\left(  U^{n}\right)  .$ Hence, $df_{n}$ is a compactly supported
$V-$valued one-form which on $\phi^{n}\left(  U^{n}\right)  $ agrees with%
\[
d\left(  g\circ\phi^{n}|_{U^{n}}^{-1}\right)  =\left(  \phi^{n}|_{U^{n}}%
^{-1}\right)  ^{\ast}\left(  dg\right)
\]
and it follows that
\[
\pi_{1}Z^{n+1}\left(  dg\right)  _{s_{n},s_{n+1}}=\pi_{1}\left(  \phi_{\ast
}^{n}Z\right)  |_{\left[  s_{n},s_{n+1}\right]  }\left(  df_{n}\right)
_{s_{n},s_{n+1}}=\pi_{1}\int_{s_{n}}^{s_{n+1}}df_{n}\left(  Y_{u}\right)
dY_{u},
\]
where $Y$ has starting point $z^{n}$ and increments $Z\left(  d\phi
^{n}\right)  |_{\left[  s_{n},s_{n+1}\right]  }.$ This implies
\begin{equation}
\pi_{1}Z^{n}\left(  dg\right)  _{s_{n},s_{n+1}}=f_{n}\left(  \phi^{n}\left(
z^{n}\right)  +\pi_{1}Z\left(  d\phi^{n}\right)  _{s_{n},s_{n+1}}\right)
-f_{n}\left(  \phi^{n}\left(  z^{n}\right)  \right)  , \label{consist}%
\end{equation}
and because the support of $\left(  \phi_{\ast}^{n}Z\right)  |_{\left[
s_{n},s_{n+1}\right]  }$ is contained in the coordinate neighbourhood we have
\[
\left\{  z^{n}+\pi_{1}Z\left(  d\phi^{n}\right)  _{s_{n},t}:t\in\left[
s_{n},s_{n+1}\right]  \right\}  \subseteq\phi^{n}\left(  U^{n}\right)  .
\]
Thus, (\ref{consist}) yields%
\begin{align*}
g\left(  z^{n}\right)  +\pi_{1}Z^{n}\left(  dg\right)  _{s_{n},s_{n+1}}  &  =
\left[  g_{n}\circ\left(  \phi^{n}|_{U^{n}}^{-1}\right)  \right]  \left(
\phi^{n}\left(  z^{n}\right)  +\pi_{1}Z\left(  d\phi^{n}\right)
_{s_{n},s_{n+1}}\right) \\
&  =g\left(  z^{n+1}\right)  .
\end{align*}

\end{proof}

\begin{remark}
Invoking Lemma \ref{associativity of concat} we can use Lemma
\ref{consistency of localisations} to show that we can concatenate the
elements of any localising sequence $\left\{  \left(  Z^{n},z^{n},\phi
^{n},U^{n}\right)  \right\}  _{n=1}^{N}$ to generate a new geometric $p-$rough
path on $M$ over $\left[  0,T\right]  $%
\[
Z^{1}\ast Z^{2}\ast....\ast Z^{N}:=(...(\left(  Z^{1}\ast Z^{2}\right)  \ast
Z^{3})\ast....)\ast Z^{N}%
\]
Associativity again guarantees that the order in which the brackets are
expanded in the concatenation is irrelevant.
\end{remark}

\begin{theorem}
Let $Z$ be a geometric $p-$rough path on a Lip-$\gamma_{0}$ manifold $M$ with
starting point $x.$ Then there exists a localising sequence $\left\{  \left(
Z^{n},z^{n},\phi^{n},U^{n}\right)  \right\}  _{n=1}^{N}$ for $Z.$ Moreover,
$Z$ equals (up to equivalence) the concatenation of its localisations; that
is
\begin{equation}
Z=Z^{1}\ast Z^{2}\ast....\ast Z^{N}. \label{concat equality}%
\end{equation}

\end{theorem}

\begin{proof}
The existence of the localising sequence is just the content of the above
recursive procedure. The fact that the concatenation in (\ref{concat equality}%
) is a geometric $p$-rough path follows from Lemmas \ref{concat is a grp} and
\ref{consistency of localisations}. \ To check that the equality in
(\ref{concat equality}) holds is straight forward from the construction of the
sequence; more precisely, if $D=\left\{  s_{n}:n=0,1,...,N\right\}  $ again
denotes the partition underlying the sequence then, using the multiplicative
property of $Z\left(  \alpha\right)  $ we have for any $\left(  s,t\right)  $
in $\Delta\left[  0,T\right]  $ with $s_{k-1}<s<s_{k}<...<s_{l}<t<s_{l+1}$ and
any Lip-$\left(  \gamma-1\right)  $ one-form:%
\begin{align*}
Z\left(  \alpha\right)  _{s,t}  &  =Z\left(  \alpha\right)  _{s,s_{k}}\otimes
Z\left(  \alpha\right)  _{s_{k},s_{k+1}}\otimes...\otimes Z\left(
\alpha\right)  _{s_{l,}t}\\
&  =\left(  Z\left(  \alpha\right)  |_{\left[  s_{k-1},s_{k}\right]  }\right)
_{s,s_{k}}\otimes\left(  Z\left(  \alpha\right)  |_{\left[  s_{k}%
,s_{k+1}\right]  }\right)  _{s_{k},s_{k+1}}\otimes...\otimes\left(  Z\left(
\alpha\right)  |_{\left[  s_{l},s_{l+1}\right]  }\right)  _{s_{l},t}\\
&  =Z^{k}\left(  \alpha\right)  _{s,s_{k}}\otimes Z^{k+1}\left(
\alpha\right)  _{s_{k},s_{k+1}}\otimes...\otimes Z^{l}\left(  \alpha\right)
_{s_{l},t}\\
&  =\left(  Z^{1}\ast Z^{2}\ast....\ast Z^{N}\right)  \left(  \alpha\right)
_{s,t}.
\end{align*}

\end{proof}

\begin{corollary}
Let $\gamma_{0}>p\geq1.$ Then the support of any geometric $p-$rough path on a
Lip-$\gamma_{0}$ manifold is compact.
\end{corollary}

\begin{proof}
Identifying any localising sequence it is easy to show that
\[
\text{supp }Z=\bigcup\nolimits_{i=1}^{N}\text{supp }Z^{i}=\bigcup
\nolimits_{i=1}^{N}\left(  \phi^{i}|_{U^{i}}^{-1}\right)  \left(  \text{supp
}\left(  \phi_{\ast}^{i}Z\right)  |_{\left[  s_{i-1},s_{i}\right]  }\right)
.
\]

\end{proof}

\subsection{Proof of Theorem \ref{one chart localisation}\label{sec-proof1}}

The embedding map constructed in the following lemma is in spirit similar to
the usual construction of Whitney embedding maps for compact manifolds.

\begin{lemma}
\label{whitney}Let $\gamma_{0},\gamma$ be such that $1<\gamma\leq\gamma_{0}$
and suppose $\alpha_{1},...,\alpha_{l}$ is a finite collection of compactly
supported Lip$\left(  \gamma-1\right)  $ one-forms on a Lipschitz-$\gamma_{0}$
manifold $M,$ which take values in some Banach space $W$. Then, for some $n\in%
%TCIMACRO{\U{2115} }%
%BeginExpansion
\mathbb{N}
%EndExpansion
$ we can find a compactly supported Lip$-\gamma$ function $g:M\rightarrow%
%TCIMACRO{\U{211d} }%
%BeginExpansion
\mathbb{R}
%EndExpansion
^{n},$ and a collection of $W-$valued Lip$-(\gamma-1)$ one-forms on $%
%TCIMACRO{\U{211d} }%
%BeginExpansion
\mathbb{R}
%EndExpansion
^{n},$ $\beta_{1},...,\beta_{l}$ $,$ such that the pull-back of $\beta_{i} $
under $g$ coincides with $\alpha_{i},$ i.e.%
\[
\alpha_{i}=g^{\ast}\beta_{i}\text{ \ for }i=1,....,l.
\]

\end{lemma}

\begin{proof}
Suppose $\mathcal{A=}\left\{  \left(  \phi_{i},U_{i}\right)  :i\in I\right\}
$ is the Lip-$\gamma_{0}$ atlas associated with $M$. Consider $\left\{
\left(  \tilde{U}_{i},\tilde{\phi}_{i}\right)  :i\in I\right\}  ,$ where
$\tilde{U}_{i}=U_{i}^{\delta/2}$ and $\tilde{\phi}_{i}$ are Lip-$\gamma_{0}$
maps satisfying $\tilde{\phi}_{i}\equiv\phi_{i}$ on $\tilde{U}_{i}$ and
vanishing off $U_{i}$. Using Lemma \ref{partition-of unity} we can find a
Lip-$\gamma_{0}$ partition of unity $\left\{  f_{i}:i\in I\right\}  $
subordinate to the atlas $\left\{  \left(  \tilde{U}_{i},\tilde{\phi}%
_{i}\right)  :i\in I\right\}  $ with the same index set; in particular $f_{i}$
is (compactly) supported in $U_{i}$ for each $i\in I$ and is strictly positive
on $\hat{U}_{i}:=\tilde{U}_{i}^{\delta/2}.$

\underline{Step 1}: We first identify the function $g$. Before doing so we
have to establish a target Euclidean space which is large enough to contain
enough information about each of the $l$ one forms $\alpha_{1},...,\alpha_{l}%
$. To this end let $J_{i}=\left\{  j_{1}^{i},....,j_{k_{i}}^{i}\right\}  $
denote the set of all indices $j\in I$ such that $U_{j}$ has non-empty
intersection with supp$\alpha_{i};$ because supp$\alpha_{i}$ is compact and
$M$ is locally finite it is guaranteed that $J_{i}$ is a finite set. Using the
partition of unity we may decompose each compactly supported one-form
$\alpha_{i}$ as the sum
\begin{equation}
\alpha_{i}=\sum_{j\in J_{i}}f_{j}\alpha_{i},\text{ }i\in\left\{
1,....,l\right\}  . \label{decomposition}%
\end{equation}
On each chart $U_{j}$ it is rudimentary to represent $\alpha_{i}\left(
\cdot\right)  $ as $\gamma_{j}\left(  \cdot\right)  \circ d\phi_{j}\left(
\cdot\right)  $ for some $\gamma_{j}:U_{j}\rightarrow L\left(
%TCIMACRO{\U{211d} }%
%BeginExpansion
\mathbb{R}
%EndExpansion
^{d},W\right)  $ which is Lip$\left(  \gamma-1\right)  $; in fact, we may
identify $\gamma_{j}$ explicitly as%
\[
\gamma_{j}\left(  m\right)  \left(  v\right)  =\alpha_{i}\left(  m\right)
\left(  \phi_{j}^{-1}\right)  _{\ast}\left(  v\right)
\]
where $v$ is, as usual, identified with $v_{\phi_{i}^{j}\left(  m\right)  }\in
T_{\phi_{i}^{j}\left(  m\right)  }%
%TCIMACRO{\U{211d} }%
%BeginExpansion
\mathbb{R}
%EndExpansion
^{d}\cong%
%TCIMACRO{\U{211d} }%
%BeginExpansion
\mathbb{R}
%EndExpansion
^{d}$. It follows that we can write $\alpha_{i}$ as
\[
\alpha_{i}\left(  \cdot\right)  =\sum_{j\in I_{i}}f_{j}\left(  \cdot\right)
\gamma_{j}\left(  \cdot\right)  \circ d\phi_{j}\left(  \cdot\right)  .
\]
We define the natural number $n$ (which will be the dimension of our target
space) by taking $n:=\left(  d+1\right)  \left(  k_{1}+....+k_{l}\right)  $
and let the functions $f:M\rightarrow%
%TCIMACRO{\U{211d} }%
%BeginExpansion
\mathbb{R}
%EndExpansion
^{\left(  k_{1}+....+k_{l}\right)  }$ and $\phi:M\rightarrow%
%TCIMACRO{\U{211d} }%
%BeginExpansion
\mathbb{R}
%EndExpansion
^{\left(  k_{1}+....+k_{l}\right)  d}$ be given by
\begin{align*}
f\left(  m\right)   &  =\left(  f_{j_{1}^{1}}\left(  m\right)
,....,f_{j_{k_{1}}^{1}}\left(  m\right)  ,....,f_{j_{1}^{l}}\left(  m\right)
,....,f_{j_{k_{l}}^{l}}\left(  m\right)  \right) \\
\phi\left(  m\right)   &  =\left(  \phi_{j_{1}^{1}}\left(  m\right)
,....,\phi_{j_{k_{1}}^{1}}\left(  m\right)  ,....,\phi_{j_{1}^{l}}\left(
m\right)  ,....,\phi_{j_{k_{l}}^{l}}\left(  m\right)  \right)  .
\end{align*}
Finally, we can identify the required function $g:M\rightarrow%
%TCIMACRO{\U{211d} }%
%BeginExpansion
\mathbb{R}
%EndExpansion
^{n}$ by setting%
\[
g\left(  m\right)  =\left(  f\left(  m\right)  ,\phi\left(  m\right)  \right)
\in%
%TCIMACRO{\U{211d} }%
%BeginExpansion
\mathbb{R}
%EndExpansion
^{\left(  k_{1}+....+k_{l}\right)  }\times%
%TCIMACRO{\U{211d} }%
%BeginExpansion
\mathbb{R}
%EndExpansion
^{\left(  k_{1}+....+k_{l}\right)  d}\cong%
%TCIMACRO{\U{211d} }%
%BeginExpansion
\mathbb{R}
%EndExpansion
^{n}%
\]
\qquad\qquad\underline{Step 2}: We now construct the one-forms $\beta
_{1},...,\beta_{l}$ . To do this we need to introduce for $i\in I$ the
functions $h_{i}:%
%TCIMACRO{\U{211d} }%
%BeginExpansion
\mathbb{R}
%EndExpansion
^{d}\rightarrow L\left(
%TCIMACRO{\U{211d} }%
%BeginExpansion
\mathbb{R}
%EndExpansion
^{d},W\right)  $ defined by
\begin{equation}
h_{i}\left(  x\right)  =\left\{
\begin{array}
[c]{c}%
\left(  \gamma_{i}\circ\phi_{i}|_{U_{i}}^{-1}\right)  \left(  x\right)
\chi_{i}\left(  x\right)  \text{ \ \ \ \ on }\phi_{i}\left(  U_{i}\right)
\text{ \ \ \ \ \ \ \ \ }\\
\text{ \ \ \ \ \ \ \ \ \ \ \ \ \ \ \ }0\text{\ \ \ \ \ \ \ \ }%
\ \text{\ \ \ \ \ \ \ \ \ \ }\ \ \ \ \ \text{ on }%
%TCIMACRO{\U{211d} }%
%BeginExpansion
\mathbb{R}
%EndExpansion
^{d}\setminus\phi_{i}\left(  U_{i}\right)
\end{array}
\right.  \text{\ ,} \label{chart one form}%
\end{equation}
where $\chi_{i}$ $\in C^{\infty}\left(
%TCIMACRO{\U{211d} }%
%BeginExpansion
\mathbb{R}
%EndExpansion
^{d}\right)  $ is constructed so that supp $\chi_{i}\subseteq\phi_{i}\left(
U_{i}\right)  $ and moreover
\[
\chi_{i}\left(  x\right)  =\left\{
\begin{array}
[c]{c}%
\text{ }1\text{\ \ on }\phi_{i}\left(  \tilde{U}_{i}\right)
\text{\ \ \ \ \ \ \ \ \ \ \ }\\
0\text{\ \ on }%
%TCIMACRO{\U{211d} }%
%BeginExpansion
\mathbb{R}
%EndExpansion
^{d}\setminus\phi_{i}\left(  U_{i}\right)  \text{ \ \ \ \ }%
\end{array}
\right.  .
\]
It is trivial to see that $h_{i}$ is Lip-$\left(  \gamma-1\right)  $. We
uniformly bound the (finite subset of) partition of unity functions by letting%
\[
C=\max_{i\in\left\{  1,...,l\right\}  }\max_{k\in\left\{  1,....,k_{i}%
\right\}  }\sup_{m\in M}\left\vert f_{j_{k}^{i}}\left(  m\right)  \right\vert
<\infty
\]
and thus define the $W-$valued Lip-$\left(  \gamma-1\right)  $ one-forms
\[
\beta_{i}:\left[  -C,C\right]  ^{k_{1}+...+k_{l}}\times%
%TCIMACRO{\U{211d} }%
%BeginExpansion
\mathbb{R}
%EndExpansion
^{\left(  k_{1}+...+k_{l}\right)  d}\rightarrow L\left(
%TCIMACRO{\U{211d} }%
%BeginExpansion
\mathbb{R}
%EndExpansion
^{n},W\right)
\]
by
\[
\beta_{i}\left(  y_{1},...,y_{l},x_{1},....,x_{l}\right)  \left(
v_{1},...,v_{l},w_{1},...,w_{l}\right)  =\sum_{k=1}^{k_{i}}y_{i}^{k}%
h_{j_{k}^{i}}\left(  x_{i}^{k}\right)  w_{i}^{k}%
\]
where for $i=1,...,l$%
\begin{align*}
y_{i}  &  =\left(  y_{i}^{1},....,y_{i}^{k_{i}}\right)  \in\left[
-C,C\right]  ^{k_{i}}\\
x_{i}  &  =\left(  x_{i}^{1},...,x_{i}^{k_{i}}\right)  \text{, }w_{i}=\left(
w_{i}^{1},....,w_{i}^{k_{i}}\right)  \in%
%TCIMACRO{\U{211d} }%
%BeginExpansion
\mathbb{R}
%EndExpansion
^{k_{i}d}%
\end{align*}
As usual we can (and do) extend $\beta_{i}$ to a Lip-$\left(  \gamma-1\right)
$ one-form on $%
%TCIMACRO{\U{211d} }%
%BeginExpansion
\mathbb{R}
%EndExpansion
^{n},$ which we persist in calling $\beta_{i}.$

\underline{Step 3}: We verify the conclusion of the theorem. Firstly, the
definition of the pull-back gives that for any $m\in M,$ $v_{m}\in T_{m}M$ and
all $i=1,....,l$
\begin{equation}
\left(  g^{\ast}\beta_{i}\right)  \left(  m\right)  \left(  v_{m}\right)
=\beta_{i}\left(  g\left(  m\right)  \right)  \left(  g_{\ast}v_{m}\right)
=\sum_{k=1}^{k_{i}}f_{_{j_{k}^{i}}}\left(  m\right)  h_{j_{k}^{i}}\left(
\phi_{j_{k}^{i}}\left(  m\right)  \right)  \left(  \phi_{j_{k}^{i}}\right)
_{\ast}v_{m} \label{pull-back1}%
\end{equation}
and because supp $f_{i}\subset\tilde{U}_{i}$ we can observe that $f_{i}\left(
m\right)  h_{i}\left(  \phi_{i}\left(  m\right)  \right)  =f_{_{i}}\left(
m\right)  \gamma_{i}\left(  m\right)  $ for all $m\in M$ and $i\in I$. This
and (\ref{pull-back1}) immediately give that%
\[
\left(  g^{\ast}\beta_{i}\right)  \left(  m\right)  \left(  v_{m}\right)
=\sum_{k=1}^{k_{i}}f_{_{j_{k}^{i}}}\left(  m\right)  \gamma_{j_{k}^{i}}\left(
m\right)  \circ d\phi_{j_{k}^{i}}\left(  m\right)  v_{m}=\alpha_{i}\left(
m\right)  \left(  v_{m}\right)
\]
and hence $g^{\ast}\beta_{i}=\alpha_{i}$ for $i=1,....,l.$
\end{proof}

The following proposition is a key step in localising the support of the rough
path to a chart. Its proof builds on the techniques developed in the previous lemma.

\begin{proposition}
\label{localisation proposition}Let $1\leq p<\gamma\leq\gamma_{0}$ and suppose
that $Z$ is a geometric $p$-rough path on a Lipschitz-$\gamma_{0}$ manifold
$M$ with starting point $z$ and that $\alpha$ and $\beta$ are two compactly
supported Lip-$\left(  \gamma-1\right)  $ one-forms on $M$ taking values in
the same Banach space $W$. Suppose further that for some open subset of the
manifold $U\subseteq M$ we have%
\[
\text{supp }\beta\subset U\text{ \ and }\alpha\equiv\beta\text{ on }U
\]
then $Z\left(  \alpha\right)  \equiv0$ implies $Z\left(  \beta\right)
\equiv0.$
\end{proposition}

\begin{proof}
Let us note immediately that $U$ may be taken, without loss of generality, to
be precompact. This follows from the compactness of supp $\beta$ which allows
us to cover supp $\beta$ with finitely many precompact basis sets
$S_{1},...,S_{l}$ (this always exists for any topological - see Lee
\cite{lee}); hence, we may replace $U$ with the precompact set $\cup_{i=1}%
^{l}U\cap S_{i}$ without affecting the hypotheses. We will now finesse the
proof of the previous theorem somewhat. Recall that we constructed a
Lip-$\gamma$ function $g:M\rightarrow%
%TCIMACRO{\U{211d} }%
%BeginExpansion
\mathbb{R}
%EndExpansion
^{n},$ for some $n,$ and a Lip-$\left(  \gamma-1\right)  $ one form
$\xi_{\alpha}$ on $%
%TCIMACRO{\U{211d} }%
%BeginExpansion
\mathbb{R}
%EndExpansion
^{n}$ such that $\alpha=g^{\ast}\xi_{\alpha}.$ For our current purposes we
will need to modify our definition of $g$ and also expand the dimension of the
target space to capture more information about the support of the one-forms.

\underline{Step 1}:First we detail the modifications needed to the earlier
argument. Let $J_{\alpha}=\left\{  j_{1},...,j_{k_{\alpha}}\right\}  $ be the
set of all indices such that supp $\alpha$ has non-empty intersection with
$U_{i},$ \ and then recall from Lemma \ref{bump-functions} that we can find a
collection of real-valued Lip-$\gamma_{0}$ functions on $M,$ $\left\{
c_{i}:i\in I\right\}  ,$ such that $c_{i}$ vanishes outside $U_{i}$ and is
identically one on $\tilde{U}_{i}.$ We define $n=\left(  2+d\right)
k_{\alpha}$ and a Lip-$\left(  \gamma-1\right)  $ one-form on $%
%TCIMACRO{\U{211d} }%
%BeginExpansion
\mathbb{R}
%EndExpansion
^{n}$ by
\begin{gather}
\xi_{\alpha}\left(  x,y,z\right)  =\sum_{j\in J_{\alpha}}y^{j}h_{j}^{\alpha
}\left(  x^{j},z^{j}\right) \label{one form definition}\\
x=\left(  x^{j_{1}},...,x^{j_{k_{\alpha}}}\right)  \in%
%TCIMACRO{\U{211d} }%
%BeginExpansion
\mathbb{R}
%EndExpansion
^{dk_{\alpha}},\text{ }y=\left(  y^{j_{1}},...,y^{j_{k_{\alpha}}}\right)  \in%
%TCIMACRO{\U{211d} }%
%BeginExpansion
\mathbb{R}
%EndExpansion
^{k_{\alpha}},\text{ }z=\left(  z^{j_{1}},...,z^{j_{k_{\alpha}}}\right)  \in%
%TCIMACRO{\U{211d} }%
%BeginExpansion
\mathbb{R}
%EndExpansion
^{k_{\alpha}};\nonumber
\end{gather}
where $h_{j}^{\alpha}$\ is the coordinate representation of the one form
$\alpha$ defined in a similar way to (\ref{chart one form}), namely:
\[
h_{i}^{\alpha}\left(  x,z\right)  =\left\{
\begin{array}
[c]{cc}%
\left(  \gamma_{i}^{\alpha}\circ\phi_{i}^{-1}\right)  \left(  x\right)
\chi_{i}\left(  x\right)  \zeta\left(  z\right)  & \text{ \ on }\phi
_{i}\left(  U_{i}\right)  \text{ \ \ \ \ \ \ \ \ }\\
0\text{\ \ \ \ \ \ \ } & \text{on }%
%TCIMACRO{\U{211d} }%
%BeginExpansion
\mathbb{R}
%EndExpansion
^{d}\setminus\phi_{i}\left(  U_{i}\right)
\end{array}
\right.  ;
\]
$\gamma_{j}^{\alpha}:U_{j}\rightarrow L\left(
%TCIMACRO{\U{211d} }%
%BeginExpansion
\mathbb{R}
%EndExpansion
^{d},W\right)  $ \ is again given by
\[
\gamma_{j}^{\alpha}\left(  m\right)  \left(  v\right)  =\alpha_{i}\left(
m\right)  \left(  \phi_{j}|_{U_{j}}^{-1}\right)  _{\ast}\left(  v\right)  ;
\]
and where the additional dependence in the $z-$variable determined by the
smooth function $\zeta:%
%TCIMACRO{\U{211d} }%
%BeginExpansion
\mathbb{R}
%EndExpansion
\rightarrow%
%TCIMACRO{\U{211d} }%
%BeginExpansion
\mathbb{R}
%EndExpansion
$ which satisfies%
\[
\zeta\left(  z\right)  =\left\{
\begin{array}
[c]{c}%
\text{ \ }1\text{ on }\left\vert z\right\vert \leq\frac{1}{4}\text{\ \ }\\
0\text{\ on }\left\vert z\right\vert \geq\frac{1}{2}%
\end{array}
\right.  .
\]

Similarly to before we now define three Lip-$\gamma$ functions $\phi
:M\rightarrow$ $%
%TCIMACRO{\U{211d} }%
%BeginExpansion
\mathbb{R}
%EndExpansion
^{dk_{\alpha}},$ $\bar{c}:M\rightarrow$ $%
%TCIMACRO{\U{211d} }%
%BeginExpansion
\mathbb{R}
%EndExpansion
^{k_{\alpha}}$ and $f:M\rightarrow$ $%
%TCIMACRO{\U{211d} }%
%BeginExpansion
\mathbb{R}
%EndExpansion
^{k_{\alpha}}$ by
\begin{align*}
\phi\left(  m\right)   &  =\left(  \phi_{j_{1}}(m),....,\phi_{j_{k_{\alpha}}%
}\left(  m\right)  \right) \\
f\left(  m\right)   &  =\left(  f_{j_{1}}(m),....,f_{j_{k_{\alpha}}}\left(
m\right)  \right) \\
\bar{c}\left(  m\right)   &  =\left(  1-c_{j_{1}}(m),....,1-c_{j_{k_{\alpha}}%
}\left(  m\right)  \right)
\end{align*}
and then we augment the function considered in Theorem \ref{whitney} by
defining $q:M\rightarrow$ $%
%TCIMACRO{\U{211d} }%
%BeginExpansion
\mathbb{R}
%EndExpansion
^{n}$ via
\[
q\left(  m\right)  =\left(  \phi\left(  m\right)  ,f\left(  m\right)  ,\bar
{c}\left(  m\right)  \right)  .
\]

It is now easy to see using the same argument as before that $q^{\ast}%
\xi_{\alpha}=\alpha;$ we only need to make sure that the additional dependence
given to $h_{j}$ does not affect the analysis. But this follows from the
observation that if $m$ is in the support of $f_{j}$ then, since this implies
that $m$ is in $\tilde{U}_{i}$ (on which $c_{j}$ is identically one), we must
have $\zeta\left(  1-c_{j}\left(  m\right)  \right)  =1$ and hence for all
$m\in M$
\[
f_{j}\left(  m\right)  \zeta\left(  1-c_{j}\left(  m\right)  \right)
=f_{j}\left(  m\right)  .
\]
Because supp $\beta\subseteq$ supp $\alpha$ the same calculation shows that
the one-form on $%
%TCIMACRO{\U{211d} }%
%BeginExpansion
\mathbb{R}
%EndExpansion
^{n}$ defined by
\[
\xi_{\beta}\left(  x,y,z\right)  :=\sum_{j\in J_{\alpha}}y^{j}h_{j}^{\beta
}\left(  z^{j},x^{j}\right)
\]
satisfies $q^{\ast}\xi_{\beta}=\beta$.

\underline{Step 2:} We now claim that for any $w\in q\left(  M\right)
\subseteq%
%TCIMACRO{\U{211d} }%
%BeginExpansion
\mathbb{R}
%EndExpansion
^{n}$ there exists a neighbourhood $N\left(  w\right)  $ of $w$ (in $%
%TCIMACRO{\U{211d} }%
%BeginExpansion
\mathbb{R}
%EndExpansion
^{n}$) on which one of the following possibilities holds: either $\xi_{\beta
}\ $vanishes identically or $\xi_{\beta}$ and $\xi_{\alpha}$ are identically
equal. With this objective in mind we suppose $w=q\left(  m\right)  $ and
define $\tilde{J}_{\alpha}\subseteq J_{\alpha}$ to be the set of all indices
$j$ in $J_{\alpha}$ such that $m$ is contained in the chart $U_{j}$. \ Denote
the components of $w=\left(  x,y,z\right)  $ by
\[
x=\left(  x^{j_{1}},...,x^{j_{k_{\alpha}}}\right)  ,\text{ \ }y=\left(
y^{j_{1}},...,y^{j_{k_{\alpha}}}\right)  \text{ and }z=\left(  z^{j_{1}%
},...,z^{j_{k_{\alpha}}}\right)  .
\]
We will divide the analysis into two cases:

\underline{Case A:} \ Suppose $j$ is in $\tilde{J}_{\alpha}.$ At least one of
the following must be true: either $m$ is in $U,$ the open set in $M$ on which
$\alpha$ and $\beta$ agree, or $m$ is in $M\setminus$supp $\beta$. In the
first situation we may define $N\left(  m\right)  $ to be some neighbourhood
of $m$ in $M$ such that
\[
N\left(  m\right)  \subseteq U_{j}\cap U.
\]
In the second, we choose $N\left(  m\right)  $ such that $\beta\equiv0$ on
$N\left(  m\right)  ,$ so that $N\left(  m\right)  $ is again contained in
$U_{j}\cap U.$ We then use $N\left(  m\right)  $ to define a neighbourhood of
$w $ in $%
%TCIMACRO{\U{211d} }%
%BeginExpansion
\mathbb{R}
%EndExpansion
^{n}$ by setting
\[
N_{j}\left(  w\right)  =\underset{\phi_{j}\left(  N\left(  m\right)  \right)
\text{ in position }l\text{ where }j=i_{l}}{\text{ }\underbrace{%
%TCIMACRO{\U{211d} }%
%BeginExpansion
\mathbb{R}
%EndExpansion
^{d}\times...\times\phi_{j}\left(  N\left(  m\right)  \right)  \times...\times%
%TCIMACRO{\U{211d} }%
%BeginExpansion
\mathbb{R}
%EndExpansion
^{d}}}\times%
%TCIMACRO{\U{211d} }%
%BeginExpansion
\mathbb{R}
%EndExpansion
^{k_{\alpha}}\times%
%TCIMACRO{\U{211d} }%
%BeginExpansion
\mathbb{R}
%EndExpansion
^{k_{\alpha}}.
\]
It is easy to see that on $N_{j}\left(  w\right)  $ we have $h_{j}^{\alpha
}\equiv h_{j}^{\beta}\equiv0$ or $h_{j}^{\beta}\equiv0$ (or both) according to
whether $m$ is in $U\ $or $M\setminus$supp $\beta$ (or both).

\underline{Case B:} Suppose that $j$ is in $J_{\alpha}\setminus\tilde
{J}_{\alpha}.$ We can decompose our considerations further in two sub-cases.
If $j=j_{r}$ where $r$ is in $\left\{  1,...,k_{\alpha}\right\}  $ then

\underline{Case Bi}: $x^{j_{r}}$ is in $%
%TCIMACRO{\U{211d} }%
%BeginExpansion
\mathbb{R}
%EndExpansion
^{d}\setminus\phi_{j_{r}}\left(  U_{j_{r}}\right)  .$ Under this assumption we
have that $x^{j_{r}}$ is outside the support of the bump function $\chi
_{j_{r}}$ and thus $\chi_{j_{r}}$ \ must vanish on some neighbourhood
$N\left(  x_{j_{r}}\right)  $ in $%
%TCIMACRO{\U{211d} }%
%BeginExpansion
\mathbb{R}
%EndExpansion
^{d}.$ It follows that $h_{j_{r}}^{\alpha}\equiv h_{j_{r}}^{\beta}\equiv0$ on
the neighbourhood%
\[
N_{j_{r}}\left(  w\right)  =\underset{N\left(  x_{j_{r}}\right)  \text{ in
}j_{r}\text{th position}}{\underbrace{%
%TCIMACRO{\U{211d} }%
%BeginExpansion
\mathbb{R}
%EndExpansion
^{d}\times...\times N\left(  x_{j_{r}}\right)  \times...\times%
%TCIMACRO{\U{211d} }%
%BeginExpansion
\mathbb{R}
%EndExpansion
^{d}}}\times%
%TCIMACRO{\U{211d} }%
%BeginExpansion
\mathbb{R}
%EndExpansion
^{k_{\alpha}}\times%
%TCIMACRO{\U{211d} }%
%BeginExpansion
\mathbb{R}
%EndExpansion
^{k_{\alpha}}%
\]
of $w$ in $%
%TCIMACRO{\U{211d} }%
%BeginExpansion
\mathbb{R}
%EndExpansion
^{m}.$

\underline{Case Bii:} $x^{i_{r}}$ is in $\phi_{j_{r}}\left(  U_{j_{r}}\right)
.$ We then have $x_{j_{r}}=\phi_{j_{r}}\left(  m\right)  $ and since $j_{r}$
is in $J_{\alpha}\setminus\tilde{J}_{\alpha}$ we know that $m$ is not in the
chart $U_{i_{r}}$. Thus $z_{i_{r}}=\bar{c}_{i_{r}}\left(  m\right)
=1-c_{i_{r}}\left(  m\right)  =1$ so that $\zeta$ must vanish on some
neighbourhood $N\left(  z_{i_{r}}\right)  $ of $z_{i_{r}}$ in $%
%TCIMACRO{\U{211d} }%
%BeginExpansion
\mathbb{R}
%EndExpansion
$ and $h_{i_{r}}^{\alpha}\equiv h_{i_{r}}^{\beta}\equiv0$ on
\[
N_{i_{r}}\left(  w\right)  =%
%TCIMACRO{\U{211d} }%
%BeginExpansion
\mathbb{R}
%EndExpansion
^{dk_{\alpha}}\times%
%TCIMACRO{\U{211d} }%
%BeginExpansion
\mathbb{R}
%EndExpansion
^{k_{\alpha}}\times\underset{N\left(  z_{i_{r}}\right)  \text{ in }%
i_{r}\text{th position}}{\underbrace{%
%TCIMACRO{\U{211d} }%
%BeginExpansion
\mathbb{R}
%EndExpansion
\times...\times N\left(  z_{i_{r}}\right)  \times...\times%
%TCIMACRO{\U{211d} }%
%BeginExpansion
\mathbb{R}
%EndExpansion
}}.
\]

Putting everything together we can now define a neighbourhood of $w$ by
\[
N\left(  w\right)  =\bigcap\nolimits_{j\in J_{\alpha}}N_{j}\left(  w\right)
\]
on which we have $\xi_{\beta}\equiv0$ or $\xi_{\beta}$ $\equiv$ $\xi_{\alpha
}.$

\underline{Step 3:} We now prove the desired conclusion. Suppose for a
contradiction that $Z\left(  \beta\right)  _{s,t}\neq0$ for some $s<t$ in
$\left[  0,T\right]  ,$ then we have that
\[
Z\left(  \beta\right)  _{s,t}=\int_{s}^{t}\xi_{\beta}\left(  Y\right)
dY\neq0
\]
where $Y$ $\in G\Omega_{p}\left(
%TCIMACRO{\U{211d} }%
%BeginExpansion
\mathbb{R}
%EndExpansion
^{n}\right)  $ \ has starting point $q\left(  z\right)  $ and increments
$Z\left(  dq\right)  .$ We introduce the time $t_{0}$ in $[s,t]$ by
\[
t_{0}=\sup\left\{  u\geq s:Z\left(  \beta\right)  _{s,v}=0\text{ for all }%
v\in\left[  s,u\right]  \right\}  ,
\]
and note from the definition of $t_{0}$ and the continuity of the path
$y:\left[  0,T\right]  \rightarrow%
%TCIMACRO{\U{211d} }%
%BeginExpansion
\mathbb{R}
%EndExpansion
^{n}$ given by $y_{u}=q\left(  z\right)  +\pi_{1}Z\left(  dq\right)  _{0,u}$
that we have

\begin{enumerate}
\item $t_{0}<t;$

\item $Z\left(  \beta\right)  _{t_{0},u}\neq0$ for some $u$ in $\left(
t_{0},t\right)  $ which may be taken arbitrarily close to $t_{0};$

\item $y_{t_{0}}$ belongs to the closed set $q\left(  \bar{U}\right)
\subseteq q\left(  M\right)  .$
\end{enumerate}

Only the third of these needs any justification. To see this, notice that
Lemma \ref{fix} and Theorem \ref{bijection-bspace} together imply that
\[
\text{supp }Y=\left\{  y_{u}:u\in\left[  0,T\right]  \right\}  \subseteq
\overline{q\left(  M\right)  }%
\]
and hence if $y_{t_{0}}\in\overline{q\left(  M\right)  }\setminus q\left(
\bar{U}\right)  $ the continuity of $y$ implies that that for some time
$t_{1}>t_{0}$ the path
\[
\left\{  y_{u}:u\in\left[  t_{0},t_{1}\right]  \right\}
\]
stays inside the set $\overline{q\left(  M\right)  }\setminus q\left(  \bar
{U}\right)  $ which is open in the relative topology on $\overline{q\left(
M\right)  }.$ On the other hand, if $x\in\overline{q\left(  M\right)
}\setminus q\left(  \bar{U}\right)  $ then it follows from Steps 1 and 2 that
$x $ is the limit of the sequence $\left(  x_{l}\right)  _{l=1}^{\infty}$ in
$q\left(  M\right)  \setminus q\left(  \bar{U}\right)  .$ Because each
$x_{l}=q\left(  m_{l}\right)  $ with $m_{l}$ in $M\setminus$supp $\beta,$ the
proof of Step 2 shows that $\xi_{\beta}$ vanishes identically on some
neighbourhood of $x_{l};$ from which it is not difficult to deduce that for
all $t$ in $\left[  t_{0},t_{1}\right]  $ we have
\[
\int_{t_{0}}^{t}\xi_{\beta}\left(  Y\right)  dY=0
\]
contradicting the definition of $t_{0}$.

Equipped with these facts we can use the conclusion of Step 2 of the proof to
find a neighbourhood $N\left(  y_{t_{o}}\right)  $ of $y_{t_{o}}$ in $%
%TCIMACRO{\U{211d} }%
%BeginExpansion
\mathbb{R}
%EndExpansion
^{n}$ on which $\xi_{\alpha}$ and $\xi_{\beta}$ are identically equal; the
only alternative would be that $\xi_{\beta}$ must vanish on a neighbourhood
which would contradict the definition of $t_{0}$. Moreover, we can find
$u>t_{0}$ such that both $Z\left(  \beta\right)  _{t_{0},u}$ is non-zero and
such that the base path stays inside this neighbourhood $\left\{  y_{r}%
:r\in\left[  t_{0},u\right]  \right\}  \subseteq N\left(  y_{t_{o}}\right)  $
$.$ Together, these immediately yield the contradiction that
\[
0\neq Z\left(  \beta\right)  _{t_{0},u}=\int_{t_{0}}^{u}\xi_{\beta}\left(
Y\right)  dY=\int_{t_{0}}^{u}\xi_{\alpha}\left(  Y\right)  dY=Z\left(
\alpha\right)  _{t_{0},u}.
\]

\end{proof}

As a corollary to the preceding proposition we deduce that any two one forms
that agree on an open neighbourhood containing the support of a rough path
have the same integral against that path.

\begin{corollary}
\label{localisation 2}Let $1\leq p<\gamma\leq\gamma_{0}$. Suppose that $Z$ is
a geometric $p$-rough path on a Lipschitz-$\gamma_{0}$ manifold $M$ with
starting point $z$ and that $\alpha$ and $\beta$ are two compactly supported
Lip-$\left(  \gamma-1\right)  $ one-forms on $M$ taking values in the same
Banach space $W$. Suppose further that for some open subset of the manifold
$U\subseteq M$ we have%
\[
\alpha\equiv\beta\text{ on }U
\]
and also that the support of $Z$ is contained in $U,$ then we have $Z\left(
\alpha\right)  =Z\left(  \beta\right)  .$
\end{corollary}

\begin{proof}
Once again define $q:M\rightarrow%
%TCIMACRO{\U{211d} }%
%BeginExpansion
\mathbb{R}
%EndExpansion
^{n},$ $\xi_{\alpha}$ and $\xi_{\beta}$ as in (\ref{one form definition}). We
will first show that
\[
\text{supp }q_{\ast}Z\subseteq\overline{q\left(  U\right)  }.
\]
We notice that if $\xi$ is any compactly supported Lip-$\left(  \gamma
-1\right)  $ one form on $%
%TCIMACRO{\U{211d} }%
%BeginExpansion
\mathbb{R}
%EndExpansion
^{n}$ with support in the open set $N\setminus\overline{q\left(  U\right)  }$
then the definition of the pull-back shows immediately that $q^{\ast}\xi
\equiv0$ on $U$ and hence, by continuity, $q^{\ast}\xi\equiv0$ on $\bar{U}$ .
It follows because $U$ is open that
\[
\text{supp }q^{\ast}\xi:=\overline{\left\{  m\in M:q^{\ast}\xi\left(
m\right)  \neq0\right\}  }\subseteq\overline{M\setminus\bar{U}}\subseteq
M\setminus U;
\]
hence using that $Z$ is supported in $U$ and $q$ is compactly supported (so
that $q^{\ast}\xi$ is also compactly supported) we have that $0=Z\left(
q^{\ast}\xi\right)  =\left(  q_{\ast}Z\right)  \left(  \xi\right)  .$
Consequently, the support of the pushed forward rough path $q_{\ast}Z$ is
contained in $\overline{q\left(  U\right)  }$.

Using the proof of Proposition \ref{localisation proposition} we have that for
any $y$ in $q\left(  U\right)  $ there exists a neighbourhood of $y$ in $%
%TCIMACRO{\U{211d} }%
%BeginExpansion
\mathbb{R}
%EndExpansion
^{n}$ on which $\xi_{\alpha}\equiv\xi_{\beta}.$ It follows that all the
derivatives of $\xi_{\alpha}$ and $\xi_{\beta}$ up to order $\lfloor
\gamma\rfloor-1$ must agree on $q\left(  U\right)  $ and thus by continuity
they must also agree on the closure $\overline{q\left(  U\right)  }.$ Because
$q_{\ast}Z$ is also a classical rough path by Lemma \ref{bijection-bspace} we
know that supp $q_{\ast}Z$ is the support of the rough path $Y$ having
increments $Z\left(  dq\right)  $ and starting point $q\left(  z\right)  $.
Using the definition of the classical rough integral it is then immediate
that
\[
Z\left(  \alpha\right)  =Z\left(  q^{\ast}\xi_{\alpha}\right)  \text{ }%
=\int\xi_{\alpha}\left(  Y\right)  dY=\int\xi_{\beta}\left(  Y\right)
dY=Z\left(  \beta\right)  .
\]

\end{proof}

\begin{remark}
\label{one-forms-open}Any Lip-$\gamma$ one form defined $\alpha$ on an open
precompact subset $U$ of a Lip-$\gamma_{0}$ manifold $M$ may be extended to a
Lip-$\gamma$ one form on $M.$ In particular if $\left(  \phi,U\right)  $ is a
chart and $\alpha$ is a Lip-$\gamma$ form on $U^{\delta/2}$ (see $\left(
\ref{delta-sets}\right)  $ for the definition) we can find an extension
$\tilde{\alpha}$ to $M$ such that $\left\Vert \tilde{\alpha}\right\Vert
_{Lip-\gamma}\leq C\left\Vert \alpha\right\Vert _{Lip-\gamma}, $ where $C$ is
a constant that only depends on the constants in the Lip-$\gamma$ atlas but is
independent of $\alpha$ and $U.$ Corollary \ref{localisation 2} demonstrates
that if supp$\left(  Z\right)  \subseteq U$ the integral of the path $Z$
against the one form $\alpha$ is independent of any particular extension.
Hence we will in the following (slightly abusing notation) sometimes consider
$Z\left(  \alpha\right)  $ with one forms $\alpha$ only defined on $U$ without
making explicit reference to the extension of $\alpha.$ Notice that if $f$ is
a sufficiently regular map we have by Lemma \ref{lemma-support-image}
supp$\left(  f_{\ast}Z\right)  \subseteq f(U)$ and we may arguing as before
consider $f_{\ast}Z\left(  \widetilde{\alpha}\right)  $ for one forms
$\widetilde{\alpha}$ defined only on $f\left(  U\right)  .$\bigskip
\end{remark}

We are finally ready to show that the support of a rough path remains for
sufficiently small times in a chart. For the convenience of the reader we
restate the theorem first.

\textbf{Theorem \ref{one chart localisation}}\textit{\ Let }$Z$\textit{\ be a
rough path on a Lip-}$\gamma$\textit{\ manifold }$M$\textit{\ with starting
point }$x$\textit{\ controlled by some }$\omega.$\textit{\ There exists a
strictly positive time }$t_{0}$\textit{\ independent of }$x$\textit{\ such
that}%
\[
\text{supp}(Z%
%TCIMACRO{\U{a6}}%
%BeginExpansion
\vert
%EndExpansion
_{[0,t_{0}]})\subseteq U
\]
\textit{for some chart }$\left(  \phi,U\right)  .$

\begin{proof}
[Proof of Theorem \ref{one chart localisation}]Given the starting point $x$
there exists by definition of a Lipschitz atlas a chart $(\phi,U)$ such that
$B(\phi\left(  x\right)  ,\delta)$ $\subseteq\phi\left(  U\right)  =B\left(
0,1\right)  .$ Let $L$ be the uniform bound for the Lip-constant of the
composition of chart functions in the definition of a Lip-$\gamma_{0}$ atlas.
Note that to show that a point $z\in M$ is not in supp$\left(  Z%
%TCIMACRO{\U{a6}}%
%BeginExpansion
\vert
%EndExpansion
_{\left[  0,t_{0}\right]  }\right)  $ it is sufficient to exhibit an open
neighbourhood of $z$ that misses the path $Z%
%TCIMACRO{\U{a6}}%
%BeginExpansion
\vert
%EndExpansion
_{\left[  0,t_{0}\right]  }.$ Let $z\in M$ be an arbitrary point in the
manifold such that
\begin{equation}
z\notin\phi^{-1}\left(  B(\phi\left(  x\right)  ,\delta)\right)  .
\label{for contradiction}%
\end{equation}
We will later choose $t_{0}>0$ in a way that is independent of $z$ and
demonstrate that an open neighbourhood of $z$ misses $Z%
%TCIMACRO{\U{a6}}%
%BeginExpansion
\vert
%EndExpansion
_{\left[  0,t_{0}\right]  }.$ In fact the time $t_{0}$ will only depend on the
Lip-$\gamma_{0}$ atlas.

There exists a chart $(\psi,V)$ such that $B(\psi\left(  z\right)  ,\delta)$
$\subseteq\psi\left(  V\right)  =B\left(  0,1\right)  $. Given $\delta\geq
u>0$ and $y\in V$ define sets $B_{\psi}(y,u)$ by
\[
B_{\psi}(y,u)=\psi^{-1}\left(  B(\psi(y),u)\right)  .
\]
Let $\delta^{\prime}<\delta/3L.$ We will first demonstrate that%
\begin{equation}
B_{\psi}(z,\delta^{\prime})\cap\left\{  x\right\}  =\emptyset.
\label{the-final-line}%
\end{equation}
Suppose for a contradiction that $x\in B_{\psi}(z,\delta^{\prime}).$ Then
$x\in U\cap V$ and as $B(\psi\left(  z\right)  ,\delta)$ $\subseteq B\left(
0,1\right)  $ we have $B(\psi\left(  x\right)  ,\delta-\delta^{\prime
})\subseteq B\left(  0,1\right)  $ and we may consider the set $B_{\psi
}\left(  x,u\right)  $ for any $0<u\leq\delta-\delta^{\prime}.$ As $x\in
B_{\psi}(z,\delta^{\prime})$ it is immediate from the definition that $z\in$
$B_{\psi}\left(  x,\delta^{\prime}\right)  .$ By assumuption we have
$B(\phi\left(  x\right)  ,\delta)$ $\subseteq\phi\left(  U\right)  $ and as
$\delta-\delta^{\prime}>\delta/2$ we it follows that $B(\psi\left(  x\right)
,\delta/2)$ $\subseteq\psi\left(  U\right)  $. Thus we may apply Corollary
\ref{last resort} and deduce that $B_{\psi}\left(  x,\frac{\delta-\varepsilon
}{2L}\right)  \subseteq B_{\phi}\left(  x,\frac{\delta}{2}\right)  $ for all
$\varepsilon>0.$ Now for sufficiently small $\varepsilon$ we have
\[
z\in B_{\psi}\left(  x,\delta^{\prime}\right)  \subseteq B_{\psi}\left(
x,\frac{\delta-\varepsilon}{2L}\right)  \subseteq B_{\phi}\left(
x,\frac{\delta}{2}\right)
\]
contradicting $\left(  \ref{for contradiction}\right)  .$ We deduce that
$x\notin$ $B_{\psi}(z,\delta^{\prime})$ and $\left(  \ref{the-final-line}%
\right)  $ holds. By familiar arguments we may find Lipschitz functions
$g_{1}$ and $g_{2}$ on $M$ such that such that
\[
g_{1}%
%TCIMACRO{\U{a6}}%
%BeginExpansion
\vert
%EndExpansion
_{B_{\psi}(z,\delta^{\prime}/2)}=\psi%
%TCIMACRO{\U{a6}}%
%BeginExpansion
\vert
%EndExpansion
_{B_{\psi}(z,\delta^{\prime}/2)}%
\]
and $g_{1}$ vanishes outside $B_{\psi}(z,\delta^{\prime}).$ Similarly we can
find $g_{2}$ such that $g_{2}(x)=y_{0}$ for some $y_{0}\notin\overline
{B\left(  \psi\left(  z\right)  ,\delta^{\prime}+1\right)  }$ and $g_{2}$
vanishes on $B_{\psi}(z,\delta^{\prime}/2)$ . Note that we can choose $g_{1}$
and $g_{2}$ in a way that the Lipschitz constants depends on the Lip-$\gamma$
atlas (via the constants $\delta,$ $L$) as by definition of a Lip-$\gamma$
atlas the Lipschitz constant of the chart function $\psi$ on $M$ is itself
bounded by $L$. Hence, if we let $g=g_{1}+g_{2}$ we may assume
\[
Lip(g)\leq C(\delta,L,\gamma).
\]

Recall that $g_{1}$ vanishes off $B_{\psi}(z,\delta^{\prime})$ and $B_{\psi
}(z,\delta^{\prime})\cap\left\{  x\right\}  =\emptyset$ so we deduce that
$g\left(  x\right)  =y_{0}.$ Let $Y_{t}$ be the classical rough path with
starting point $g\left(  x\right)  $ and increments $Z\left(  dg\right)  .$
Recall that by definition of $Z,$ the rough path $Y_{t}$ is controlled by
$Lip(g)\omega.$ Thus there exist a strictly positive time $t_{0}$ such that
the support of the rough path $Y_{t}$ with starting point $y_{0}=g\left(
x\right)  $ restricted to $\left[  0,t_{0}\right]  $ is contained in a ball of
radius one centred at $y_{0}$.

Let $\alpha$ be any Lip-$\left(  \gamma-1\right)  $ one form on $M$ that
vanishes outside $B_{\psi}(z,\delta^{\prime}/2).$ Then there exists a
Lip-$\left(  \gamma-1\right)  $ one form $\beta$ on $R^{d}$ such that $\beta$
vanishes outside $B\left(  \psi\left(  z\right)  ,\delta^{\prime}\right)  $
and restricted to $B_{\psi}(z,\delta^{\prime}/2)$ we have $\alpha=g^{\ast
}\beta$ (consider $\left(  \left(  g%
%TCIMACRO{\U{a6}}%
%BeginExpansion
\vert
%EndExpansion
_{_{B_{\psi}(z,\delta^{\prime}/2)}}\right)  ^{-1}\right)  ^{\ast}\alpha$ and
extend suitably using Whitney extension on $R^{d}$ ). By definition of $Z $
\[
Z(g^{\ast}\beta)=\int\beta(Y_{t})dY_{t}.
\]
As by construction $g\left(  x\right)  =y_{0}\notin\overline{B\left(
\psi\left(  z\right)  ,\delta^{\prime}+1\right)  }$ and $\beta$ vanishes
outside $B\left(  \psi\left(  z\right)  ,\delta^{\prime}\right)  $ it follows
that $\beta$ vanishes on $B(y_{0},1)\supseteq$supp$\left(  Y_{t}|_{\left[
0,t_{0}\right]  }\right)  .$ We deduce that
\[
Z(g^{\ast}\beta)%
%TCIMACRO{\U{a6}}%
%BeginExpansion
\vert
%EndExpansion
_{\left[  0,t_{0}\right]  }=0.
\]
Now $\alpha$ has support contained inside the open set $B_{\psi}%
(z,\delta^{\prime}/2)$ and in addition by construction $\alpha$ and $g^{\ast
}\beta$ agree on $B_{\psi}(z,\delta^{\prime}/2).$ Thus we may apply
Proposition \ref{localisation proposition} with $U=B_{\psi}(z,\delta^{\prime
}/2)$ and deduce that $Z(\alpha)%
%TCIMACRO{\U{a6}}%
%BeginExpansion
\vert
%EndExpansion
_{\left[  0,t_{0}\right]  }=0$ for any Lip-$\left(  \gamma-1\right)  $ one
form $\alpha$ on $M$ that vanishes outside $B_{\psi}(z,\delta^{\prime}/2).$ We
have shown that for any $z\notin B_{\phi}(x,\delta)\subseteq U$ an open
neighbourhood misses the path $Z%
%TCIMACRO{\U{a6}}%
%BeginExpansion
\vert
%EndExpansion
_{\left[  0,t_{0}\right]  }$ and hence conclude $z\notin$supp$\left(  Z%
%TCIMACRO{\U{a6}}%
%BeginExpansion
\vert
%EndExpansion
_{\left[  0,t_{0}\right]  }\right)  $ as required.
\end{proof}

\subsection{Proof of Theorem \ref{push forward localisation} and Corollary
\ref{last-in-line}\label{sec-proof2}}

Before we begin the proof of Theorem we obtain an elementary lemma. The claim
of the lemma, namely that if two sets miss the rough path their union does so
as well, is intuitively obvious; but the non-linearity of the rough path
functional requires us to give a rigorous proof.

\begin{lemma}
\label{union lemma}Let $1\leq p<\gamma_{0}$ and suppose that $Z$ is a
geometric $p$-rough path on a Lipschitz-$\gamma_{0}$ manifold $M$ with
starting point $z.$ Suppose $U$ and $V$ are two open sets in $M$ such that $Z
$ misses $U$ and $V$ then $Z$ also misses the union $U\cup V$.
\end{lemma}

\begin{proof}
Suppose that $Z$ does not miss $U\cup V.$ Then, for some (Banach space valued)
Lip$-\left(  \gamma-1\right)  $ one form $\alpha$ which is compactly supported
in $U\cup V$ and some $s<t$ in $\left[  0,T\right]  $ we have $Z\left(
\alpha\right)  _{s,t}\neq0$. \ As in the proof of Proposition
\ref{localisation proposition} define
\[
t_{0}:=\sup\left\{  u\geq s:Z\left(  \alpha\right)  _{s,v}=0\text{ for all
}v\in\left[  s,u\right]  \right\}
\]
so that $Z\left(  \alpha\right)  _{t_{0},u}\neq0$ for some $u$ in $\left(
t_{0},t\right)  $ which may be taken arbitrarily close to $t_{0}.$
Furthermore, we use the construction of that theorem to represent
$\alpha=q^{\ast}\xi_{\alpha}$ so that
\[
Z\left(  \alpha\right)  =\int\xi_{\alpha}\left(  Y\right)  dY.
\]
The definition of $t_{0}$ implies that $\xi_{\alpha}$ does not vanish
identically on any neighbourhood of $y_{t_{0}}=q\left(  z\right)  +Z\left(
dq\right)  _{0,t_{0}}$ in $%
%TCIMACRO{\U{211d} }%
%BeginExpansion
\mathbb{R}
%EndExpansion
^{n}.$ We know that $Y$ is supported in $\overline{q\left(  M\right)  },$ so
that $y_{t_{0}}$ is in $\ \overline{q\left(  M\right)  }$ but in fact we must
have
\[
y_{t_{0}}\in q\left(  \text{supp }\alpha\right)  \ \subseteq q\left(
U\right)  .
\]
To see this suppose to the contrary, then $y_{t_{0}}$ must belong to the set
\[
\ \overline{q\left(  M\right)  }\setminus q\left(  \text{supp }\alpha\right)
\]
which is open in the relative topology on $\overline{q\left(  M\right)  }$
(supp $\alpha$ is compact) and the continuity of $y$ gives that for some
$t_{1}>t_{0}$
\[
\left\{  y_{t}:t\in\left[  t_{0},t_{1}\right]  \right\}  \subseteq
\ \overline{q\left(  M\right)  }\setminus q\left(  \text{supp }\alpha\right)
.
\]
On the other hand, any $y$ in $\overline{q\left(  M\right)  }\setminus
q\left(  \text{supp }\alpha\right)  $ can be written as the limit of some
sequence $\left(  y_{l}\right)  _{l=1}^{\infty}$ where $y_{l}=q\left(
m_{l}\right)  $ with $m_{l}$ in $M\setminus$supp $\alpha.$ It follows from the
proof of Proposition \ref{localisation proposition} that $\xi_{\alpha}$
vanishes on a neighbourhood of each element of the sequence $y_{l}$, which is
enough to conclude that%
\[
\int_{t_{0}}^{t}\xi_{\alpha}\left(  Y\right)  dY=0\text{ for every }t\text{ in
}\left[  t_{0},t_{1}\right]
\]
violating the definition of $t_{0.}$ \ 

Now suppose $y_{t_{0}}=q\left(  m\right)  $ for some $m$ in the support of
$\alpha$ ($\subseteq U\cup V$) and without loss of generality assume that $m$
is in $U.$ As in Proposition \ref{localisation proposition} let $J_{\alpha}$
denote the set of indices $j$ such that the support of $\alpha$ has non-empty
intersection with the chart $U_{j}$ and let $I_{\alpha}\subseteq J_{\alpha}$
denote the set of indices $j$ such that $m$ is in $U_{j}.$ Let $N_{1}\left(
m\right)  $ and $N_{2}\left(  m\right)  $ be two neighbourhoods of $m$ in $M$
such that
\[
N_{1}\left(  m\right)  \subset N_{2}\left(  m\right)  \subset\bigcap
\nolimits_{j\in I_{\alpha}}U_{j}\cap U
\]
and let $\beta$ be any Lip$-\left(  \gamma-1\right)  $ one form which
coincides with $\alpha$ on $N_{1}\left(  m\right)  $ and vanishes on
$M\setminus N_{2}\left(  m\right)  $. Define $\xi_{\beta}$ as in
(\ref{one form definition}), then using the proof of Proposition
\ref{localisation proposition} we may construct a neighbourhood of $y_{t_{0}}$
in $%
%TCIMACRO{\U{211d} }%
%BeginExpansion
\mathbb{R}
%EndExpansion
^{n}$ on which $\xi_{\alpha}$ and $\xi_{\beta}$ coincide. Hence it follows
that for $t>t_{0}$ sufficiently close to $t_{0}$
\[
Z\left(  \alpha\right)  _{t_{0},t}=\int_{t_{0}}^{t}\xi_{\alpha}\left(
Y\right)  dY=\int_{t_{0}}^{t}\xi_{\beta}\left(  Y\right)  dY=Z\left(
\beta\right)  _{t_{0},t}.
\]
But we have established that the left hand side can be made non-zero for $t$
arbitrarily close to $t_{0};$ since $\beta$ is supported in $U$ and by
hypothesis $Z$ misses $U,$ the contradiction has thus been rendered manifest.
\end{proof}

We restate Theorem \ref{push forward localisation} before we begin its proof.

\textbf{Theorem }\ref{push forward localisation} \textit{Suppose }$\gamma
_{0}>p\geq1$\textit{\ and let }$Z$\textit{\ be a geometric }$p$\textit{-rough
path on a Lip}$-\gamma_{0}$\textit{\ manifold }$M$\textit{\ with starting
point }$x\in M.$\textit{\ Then for some strictly positive time, which does not
depend on }$x,$\textit{\ and some chart }$\ \left(  \phi,U\right)
$\textit{\ we have that the support of the pushforward of }$Z$\textit{\ under
}$\phi$\textit{\ is contained in the coordinate neighbourhood }$\phi\left(
U\right)  ,$\textit{\ i.e. }%
\[
\text{supp }\left(  \phi_{\ast}Z\right)  |_{\left[  0,t_{0}\right]  }%
\subseteq\phi\left(  U\right)  .
\]
\textit{Moreover, we can recover the support of }$Z$\textit{\ restricted to
}$\left[  0,t_{0}\right]  $\textit{\ as the preimage under }$\phi$\textit{\ of
the support of the pushed forward rough path; that is, we have }%
\[
\text{supp }Z|_{\left[  0,t_{0}\right]  }=\left(  \phi|_{U}\right)
^{-1}\left(  \text{supp }\left(  \phi_{\ast}Z\right)  |_{\left[
0,t_{0}\right]  }\right)  .
\]

\begin{proof}
[Proof of Theorem \ref{push forward localisation}]Proposition
\ref{one chart localisation} gives a strictly positive time $t_{0}$ so that
supp$(Z%
%TCIMACRO{\U{a6}}%
%BeginExpansion
\vert
%EndExpansion
_{[0,t_{0}]})\subseteq U$ for some chart $\left(  \phi,U\right)  .$ From the
properties of a Lipschitz-$\gamma$ atlas we know that the coordinate
neighbourhood $\phi\left(  U\right)  $ is precompact and because $\phi|_{U}$
is a homeomorphism it follows that the closed subset
\[
\phi\left(  \text{supp }(Z%
%TCIMACRO{\U{a6}}%
%BeginExpansion
\vert
%EndExpansion
_{[0,t_{0}]})\right)  \subseteq\phi\left(  U\right)
\]
is compact. Hence we may find an open set $V\subseteq%
%TCIMACRO{\U{211d} }%
%BeginExpansion
\mathbb{R}
%EndExpansion
^{d}$ with the property that
\[
\phi\left(  \text{supp }(Z%
%TCIMACRO{\U{a6}}%
%BeginExpansion
\vert
%EndExpansion
_{[0,t_{0}]})\right)  \subseteq V\subseteq\bar{V}\subseteq\phi\left(
U\right)  .
\]
Let $\xi$ be a Lip-$\left(  \gamma-1\right)  $ one form which is compactly
supported in the open set $%
%TCIMACRO{\U{211d} }%
%BeginExpansion
\mathbb{R}
%EndExpansion
^{d}\setminus\bar{V}.$ The chart functions are compactly supported so the
pull-back $\phi^{\ast}\xi$ will be too; moreover $\phi^{\ast}\xi\left(
m\right)  =0$ for any $m$ in the closed set $\phi|_{U}^{-1}\left(  \bar
{V}\right)  $ from which we can deduce that%
\[
\text{supp }\phi^{\ast}\xi\subseteq\overline{M\setminus\phi|_{U}^{-1}\left(
\bar{V}\right)  }=\overline{M\setminus\overline{\phi|_{U}^{-1}\left(
V\right)  }\text{ }}\subseteq M\setminus\phi|_{U}^{-1}\left(  V\right)
\subseteq M\setminus\text{supp }(Z%
%TCIMACRO{\U{a6}}%
%BeginExpansion
\vert
%EndExpansion
_{[0,t_{0}]})
\]
and hence $0=Z%
%TCIMACRO{\U{a6}}%
%BeginExpansion
\vert
%EndExpansion
_{[0,t_{0}]}\left(  \phi^{\ast}\xi\right)  =\left(  \phi_{\ast}Z\right)
%TCIMACRO{\U{a6}}%
%BeginExpansion
\vert
%EndExpansion
_{[0,t_{0}]}\left(  \xi\right)  .$ It follows at once that
\[
\text{supp }\left(  \phi_{\ast}Z\right)  |_{\left[  0,t_{0}\right]  }%
\subseteq\bar{V}\subseteq\phi\left(  U\right)  .
\]

Let $K$ be the compact subset of $M$ defined by
\[
K=\left(  \phi|_{U}\right)  ^{-1}\left(  \text{supp }\left(  \phi_{\ast
}Z\right)  |_{\left[  0,t_{0}\right]  }\right)  \subseteq U,
\]
we will show that $K$ \emph{is} the support of the restricted rough path $Z%
%TCIMACRO{\U{a6}}%
%BeginExpansion
\vert
%EndExpansion
_{[0,t_{0}]}.$ To achieve this we first show that
\begin{equation}
\text{supp }Z|_{\left[  0,t_{0}\right]  }\subseteq K \label{containment}%
\end{equation}
and then we will finish by exhibiting the reverse inclusion. To start we
observe that supp$(Z%
%TCIMACRO{\U{a6}}%
%BeginExpansion
\vert
%EndExpansion
_{[0,t_{0}]})\subseteq U$ by definition implies that the open set
$M\setminus\bar{U}$ misses the support of $Z%
%TCIMACRO{\U{a6}}%
%BeginExpansion
\vert
%EndExpansion
_{[0,t_{0}]}.$ Consequently, if we can show that $U\setminus K$ also misses
the support then, by Lemma \ref{union lemma}, the union
\[
\left(  M\setminus\bar{U}\right)  \cup U\setminus K=M\setminus K
\]
does too and (\ref{containment}) will have been proven. Therefore, we consider
a $Lip\left(  \gamma-1\right)  $ one form $\alpha$ (compactly) supported in
$U\setminus K$ and aim to show that $Z|_{\left[  0,t_{0}\right]  }\left(
\alpha\right)  =0.$ This is achieved by defining a $Lip\left(  \gamma
-1\right)  $ one form on $%
%TCIMACRO{\U{211d} }%
%BeginExpansion
\mathbb{R}
%EndExpansion
^{d}$ by
\[
\zeta_{\alpha}=\left\{
\begin{array}
[c]{c}%
\left(  \phi|_{U}^{-1}\right)  ^{\ast}\alpha\text{ \ on }\phi\left(  U\right)
\text{ \ \ \ \ \ \ \ \ \ \ \ \ \ \ \ }\\
0\text{ }\ \ \ \ \ \ \ \ \text{on }%
%TCIMACRO{\U{211d} }%
%BeginExpansion
\mathbb{R}
%EndExpansion
^{d}\setminus\phi\left(  U\right)  \text{ \ \ }%
\end{array}
\right.
\]
whilst noting that the support of this one form satisfies
\begin{equation}
\text{supp }\zeta_{\alpha}\subseteq\phi\left(  U\right)  \setminus\phi\left(
K\right)  \subseteq%
%TCIMACRO{\U{211d} }%
%BeginExpansion
\mathbb{R}
%EndExpansion
^{d}\setminus\text{supp }\left(  \phi_{\ast}Z\right)  |_{\left[
0,t_{0}\right]  } \label{zeta support}%
\end{equation}
and its pull-back under $\phi$ behaves so that
\[
\phi^{\ast}\zeta_{\alpha}\equiv\alpha\text{ on }U.
\]
Then since supp$(Z%
%TCIMACRO{\U{a6}}%
%BeginExpansion
\vert
%EndExpansion
_{[0,t_{0}]})\subseteq U$ an application of Lemma \ref{localisation 2} shows
that $Z|_{\left[  0,t_{0}\right]  }\left(  \phi^{\ast}\zeta_{\alpha}\right)
=Z|_{\left[  0,t_{0}\right]  }\left(  \alpha\right)  $ and hence by
(\ref{zeta support})%
\[
Z|_{\left[  0,t_{0}\right]  }\left(  \alpha\right)  =Z|_{\left[
0,t_{0}\right]  }\left(  \phi^{\ast}\zeta_{\alpha}\right)  =\left(  \phi
_{\ast}Z\right)  |_{\left[  0,t_{0}\right]  }\left(  \zeta_{\alpha}\right)
=0.
\]
It follows that $U\setminus K$ misses the support of $Z|_{\left[
0,t_{0}\right]  }$.

To conclude we need to show that supp$(Z%
%TCIMACRO{\U{a6}}%
%BeginExpansion
\vert
%EndExpansion
_{[0,t_{0}]})\supseteq K$, which we will prove by contradiction. Therefore,
suppose that $m$ is in supp$(Z%
%TCIMACRO{\U{a6}}%
%BeginExpansion
\vert
%EndExpansion
_{[0,t_{0}]})\cap\left(  U\setminus K\right)  ,$ let $N$ be a neighbourhood of
$m$ in $M$ which is contained in $U\setminus K$ and suppose $\beta$ is a
$Lip\left(  \gamma-1\right)  $ one form supported in $N$ such that $Z\left(
\beta\right)  $ is non-zero. Define $\chi_{\beta}$ a $Lip\left(
\gamma-1\right)  $ one form on $%
%TCIMACRO{\U{211d} }%
%BeginExpansion
\mathbb{R}
%EndExpansion
^{d}$ by
\[
\chi_{\beta}=\left\{
\begin{array}
[c]{c}%
\left(  \phi|_{U}^{-1}\right)  ^{\ast}\beta\text{ \ on }\phi\left(  U\right)
\text{ \ \ \ \ \ \ \ \ \ \ \ \ \ \ \ }\\
0\text{ \ \ \ }\ \ \ \text{on }%
%TCIMACRO{\U{211d} }%
%BeginExpansion
\mathbb{R}
%EndExpansion
^{d}\setminus\phi\left(  U\right)  \text{ \ }%
\end{array}
\right.
\]
and observe that $\chi_{\beta}$ is supported in $%
%TCIMACRO{\U{211d} }%
%BeginExpansion
\mathbb{R}
%EndExpansion
^{d}\setminus$supp $\left(  \phi_{\ast}Z\right)  |_{\left[  0,t_{0}\right]  }.
$ Thus we must have that $\left(  \phi_{\ast}Z\right)  |_{\left[
0,t_{0}\right]  }\left(  \chi_{\beta}\right)  =0$, but on the other hand
$\phi^{\ast}\chi_{\beta}\equiv\beta$ on $U$ hence Lemma \ref{localisation 2}
\ implies the contradiction
\[
0\neq Z|_{\left[  0,t_{0}\right]  }\left(  \beta\right)  =Z|_{\left[
0,t_{0}\right]  }\left(  \phi^{\ast}\chi_{\beta}\right)  =\left(  \phi_{\ast
}Z\right)  |_{\left[  0,t_{0}\right]  }\left(  \chi_{\beta}\right)  =0.
\]

\end{proof}

Finally we restate and prove Corollary \ref{last-in-line}.

\textbf{Corollary \ref{last-in-line}} \textit{Suppose }$\gamma_{0}>p\geq
1$\textit{\ and let }$Z$\textit{\ be a geometric }$p$\textit{-rough path on a
Lip}$-\gamma_{0}$\textit{\ manifold }$M$\textit{\ with starting point }$x\in
M.$\textit{\ Then for some strictly positive time }$t_{0},$\textit{\ which
does not depend on }$x,$\textit{\ and some chart }$\ \left(  \phi,U\right)
$\textit{\ we have}

\begin{enumerate}
\item $\left(  \phi_{\ast}Z\right)  |_{\left[  0,t_{0}\right]  }\left(
\xi\right)  =\left(  \phi_{\ast}Z\right)  |_{\left[  0,t_{0}\right]  }\left(
\eta\right)  $\textit{\ for any compactly supported Lip-}$\left(
\gamma-1\right)  $\textit{\ one forms which agree on }$\phi\left(  U\right)
.$

\item \textit{The restriction }$Z|_{\left[  0,t_{0}\right]  }$\textit{\ is
characterised by the push forward }$\left(  \phi_{\ast}Z\right)  |_{\left[
0,t_{0}\right]  }$\textit{\ in the sense that, for any compactly supported
Banach space valued Lip}$-\left(  \gamma-1\right)  $\textit{\ one form
}$\alpha$\textit{\ on }$M,$\textit{\ \ }%
\[
Z|_{\left[  0,t_{0}\right]  }\left(  \alpha\right)  =\left(  \phi_{\ast
}Z\right)  |_{\left[  0,t_{0}\right]  }\left(  \xi_{\alpha}\right)  ,
\]
\textit{where }$\xi_{\alpha}$\textit{\ is any compactly supported
Lip}$-\left(  \gamma-1\right)  $\textit{\ one form on }$%
%TCIMACRO{\U{211d} }%
%BeginExpansion
\mathbb{R}
%EndExpansion
^{d}$\textit{\ which agrees with }$\left(  \phi|_{U}^{-1}\right)  ^{\ast
}\alpha$\textit{\ on }$\phi\left(  U\right)  .$
\end{enumerate}

\begin{proof}
[Proof of Corollary \ref{last-in-line}]Let $t_{0}$ be as provided by Theorem
\ref{push forward localisation}. Then we have that supp$\left(  \phi_{\ast
}Z\right)  |_{\left[  0,t_{0}\right]  }$ is contained in $\phi\left(
U\right)  $ and the first conclusion of the present result follows immediately
from Lemma \ref{localisation 2}. For the second part we note that Theorem
\ref{push forward localisation} also gives that the support of $Z$ is
contained in the chart $U$; that is
\[
\text{supp}Z|_{\left[  0,t_{0}\right]  }\subseteq U.
\]
If $\xi_{\alpha}$ is any (Banach space valued, compactly supported)
Lip$-\left(  \gamma-1\right)  $ one form on $%
%TCIMACRO{\U{211d} }%
%BeginExpansion
\mathbb{R}
%EndExpansion
^{d}$ which agrees with $\left(  \phi|_{U}^{-1}\right)  ^{\ast}\alpha$ on
$\phi\left(  U\right)  $ then it is trivial to check that
\[
\alpha\equiv\phi^{\ast}\xi_{\alpha}%
\]
on $U$. And Lemma \ref{localisation 2} then implies that
\[
Z|_{\left[  0,t_{0}\right]  }\left(  \alpha\right)  =Z|_{\left[
0,t_{0}\right]  }\left(  \phi^{\ast}\alpha\right)  =\left(  \phi_{\ast
}Z\right)  |_{\left[  0,t_{0}\right]  }\left(  \xi_{\alpha}\right)  .
\]

\end{proof}

\section{Rough differential equations on a manifold}

Let $N$ and $M$ be Lip-$\gamma_{0}$ manifolds of dimension $d_{1}$ and $d_{2}$
respectively and $E=M\times N.$ We may interpret $E$ as the trivial fibre
bundle with bundle map $\pi:E\rightarrow N$, where $\pi$ is the projection
onto $N$. In the following we will think of the signal $X_{t}$ as living in
the manifold $N$ and the response $Y_{t}$ in the manifold $M$. For all $x\in
N$ let $g\left(  x,\cdot\right)  $ be a linear map from $T_{x}N$ to
$\tau\left(  M\right)  $, the space of vector fields on $M.$ The map $g$
together with a rough signal $X$ on $M$ and an initial condition $y_{0}\in N$
gives rise to a differential equation%

\begin{equation}
dY_{t}=g\left(  X_{t},Y_{t}\right)  dX_{t}. \label{rde-on-manifold}%
\end{equation}

Analogous to the classical theory of rough differential equation on Banach
spaces the solution of such a rough differential equation is a rough path on
the product of the two manifolds on which signal and response live. Before we
move on to give a formal definition of the solution of an RDE\ on a manifold
in terms of a fixed point, we briefly demonstrate how the map $g$ induces a
general form of a connection on $E$ that allows us to lift the signal from $N$
to $E$ and encodes the RDE.

\subsection{Geometric interpretation of the lift as an Ehresmann connection}

Defining a controlled differential equation on the product of two manifolds
requires us to find the right lift of the signal from $N$ to the product. For
all $x\in N,y\in M$ define a linear map $\Gamma_{(x,y)}:T_{x}N\rightarrow
T_{(x,y)}E=T_{x}N\oplus T_{y}M$ by%
\begin{equation}
\Gamma_{(x,y)}(v)=(v,g(x,y)(v)) \label{connection-map}%
\end{equation}
for all $v\in T_{x}N.$ For any $y\in M$ the map $\Gamma$ lifts tangent vectors
$v\in T_{x}N$ to tangent vectors in $T_{(x,y)}E$ and will play a crucial role
in our definition of a rough differential equation on a manifold. We can show
that the map $\Gamma$ may be interpreted as a connection in the following sense.

\begin{definition}
Given a manifold $E$ a Lip-$\gamma$ connection is a choice of subspaces
$H_{x}$ of $T_{x}E$ for all $x\in E.$ We require that all subspaces $H_{x}$
are of equal dimension and the mapping $x\rightarrow H_{x}$ is Lip-$\gamma$.
By this we mean that for any smooth vector field $V$ on $E$ the vector field
defined by projecting $V(x)$ onto $H_{x}$ for all \thinspace$x\in E$ are
Lip-$\gamma$.
\end{definition}

The family of subspaces%

\[
H_{z}=\operatorname{Im}(\Gamma_{z}),\text{ }z\in E
\]

is a Lip-$\gamma$ connection and we will in the following refer to $\Gamma$ as
the connection map. Note that by construction using the linearity of $g\left(
x,\cdot\right)  $ the $H_{(x,y)}$ are linear subspaces of dimension equal to
the dimension of $N$ and%
\[
d\pi\circ\Gamma_{(x,y)}=\text{Id}_{T_{x}N}.
\]

\subsection{Formal definition of the RDE solution via the connection}

We now come to the formal definition of a solution to a rough differential
equation
\begin{equation}
dY_{t}=g\left(  X_{t},Y_{t}\right)  dX_{t} \label{RDE2}%
\end{equation}
on a manifold. We use the connection map $\Gamma$ to define the solution as a
fixed point of an integral.

\begin{definition}
Let $X$ be a geometric $p-$ rough path on a Lip-$\gamma_{0}$ manifold $N$ with
starting point $x_{0}$ and for all $x\in N$ let $g\left(  x,\cdot\right)  $ be
a linear map from $T_{x}N$ into the space of vector fields on a Lip-$\gamma
_{0}$ manifold $M.$ We say a geometric $p-$rough path $Z$ on the manifold
$E=N\times M$ with starting point $z_{0}$ is a solution to (\ref{RDE2}) with
initial condition $y_{0}$ if
\end{definition}

\begin{enumerate}
\item $z_{0}=(x_{0},y_{0});$

\item $\pi_{\ast}Z\sim X;$

\item for every compactly supported Lip-$\gamma-1$ one-form $\alpha$ on $E$
taking values in a Banach space, $\alpha^{\Gamma}$ is Lip-$\left(
\gamma-1\right)  $ and we have
\begin{equation}
Z\left(  \alpha\right)  =Z\left(  \alpha^{\Gamma}\right)  ,
\label{solution_def}%
\end{equation}
where $\alpha^{\Gamma}$ is the one form on $E$ defined by
\[
\alpha_{z}^{\Gamma}(v_{z})=\alpha_{z}\circ\Gamma_{z}\circ(\pi_{\ast})(v_{z})
\]
for all $z$ $\in E,$ $v_{z}\in T_{z}E.$
\end{enumerate}

As usual we first verify the consistency of our definition with the classical
one if the manifolds are finite dimensional normed vector spaces.

\begin{lemma}
\label{product one forms}Suppose that $V$ is a Banach space and $f:\rightarrow
L\left(  V,V\right)  $, $g:V\rightarrow L\left(  V,W\right)  $ are two
Lip-$\left(  \gamma-1\right)  $ one forms and $Z\in G\Omega_{p}\left(
V\right)  $ is a solution to the RDE%
\[
Z=\int f\left(  Z\right)  dZ.
\]
Then we have%
\[
\int g\left(  Z_{u}\right)  dZ_{u}=\int\left(  \widehat{gf}\right)  \left(
Z_{u}\right)  dZ_{u},
\]
where $\widehat{gf}:V\rightarrow L\left(  V,W\right)  $ is the Lip-$\left(
\gamma-1\right)  $ one form defined by $\widehat{gf}\left(  v_{1}\right)
\left(  v_{2}\right)  =g\left(  v_{1}\right)  \left(  f\left(  v_{1}\right)
v_{2}\right)  .$
\end{lemma}

\begin{proof}
Since $Z\in G\Omega_{p}\left(  V\right)  $ we can find a sequence $z\left(
n\right)  $ $:\left[  0,T\right]  \rightarrow V$ of paths with bounded
variation such that if $Z\left(  n\right)  \equiv S_{\lfloor p\rfloor}\left(
z\left(  n\right)  \right)  $ we have, in the rough path $p$-variation
metric,
\[
\lim_{n\rightarrow\infty}Z\left(  n\right)  =Z.
\]
Then, by the continuity of the map $Z\longmapsto\int f\left(  Z\right)  dZ$
under this metric we have that
\[
\lim_{n\rightarrow\infty}\int f\left(  Z\left(  n\right)  \right)  d\left(
Z\left(  n\right)  \right)  =\int f\left(  Z\right)  dZ.
\]
Let $y\left(  n\right)  \equiv\int f\left(  z\left(  n\right)  \right)
dz\left(  n\right)  ,$ then we have that%
\[
Y\left(  n\right)  :=S_{\lfloor p\rfloor}\left(  y\left(  n\right)  \right)
\equiv\int f\left(  Z\left(  n\right)  \right)  d\left(  Z\left(  n\right)
\right)  .
\]
\ 

By defining $\bar{f}:V\rightarrow L\left(  V,V\oplus V\right)  $
\[
\bar{f}\left(  x\right)  \left(  v\right)  =\left(  f\left(  x\right)  \left(
v\right)  ,v\right)
\]
and
\begin{equation}
\bar{Y}\left(  n\right)  :=\int\bar{f}\left(  Z\left(  n\right)  \right)
dZ\left(  n\right)  \in G\Omega_{p}\left(  V\oplus V\right)
\label{joint path}%
\end{equation}
we can use the same arguments as before to see that $\bar{Y}\left(  n\right)
=S_{\lfloor p\rfloor}\left(  \bar{y}\left(  n\right)  \right)  $ where
$\bar{y}\left(  n\right)  =\left(  y\left(  n\right)  ,z\left(  n\right)
\right)  \in V\oplus V.$ Moreover it is easy to see from (\ref{joint path})
and continuity of the rough integral that $\bar{Y}\left(  n\right)
\rightarrow\bar{Y}\in G\Omega_{p}\left(  V\oplus V\right)  $ where $\pi
_{1}\bar{Y}=Z$ and $\pi_{2}\bar{Y}=Z$ (taking the obvious projections of
$G\Omega_{p}\left(  V_{1}\oplus V_{2}\right)  $ onto $G\Omega_{p}\left(
V_{1}\right)  $ and $G\Omega_{p}\left(  V_{2}\right)  $).

Hence, since $Y\left(  n\right)  \rightarrow Z$ and $Z\left(  n\right)
\rightarrow Z$ we have that
\begin{align}
\int g\left(  Z_{u}\right)  dZ_{u}  &  =\lim_{n\rightarrow\infty}\int g\left(
Y_{u}\left(  n\right)  \right)  dY_{u}\left(  n\right) \nonumber\\
&  =\lim_{n\rightarrow\infty}S_{\lfloor p\rfloor}\left(  \int g\left(
y_{u}\left(  n\right)  \right)  dy_{u}\left(  n\right)  \right) \nonumber\\
&  =\lim_{n\rightarrow\infty}S_{\lfloor p\rfloor}\left(  \int g\left(
y_{u}\left(  n\right)  \right)  f\left(  z_{u}\left(  n\right)  \right)
dz_{u}\left(  n\right)  \right) \nonumber\\
&  =\lim_{n\rightarrow\infty}S_{\lfloor p\rfloor}\left(  \int h\left(  \bar
{y}_{u}\left(  n\right)  \right)  d\bar{y}_{u}\left(  n\right)  \right)
\nonumber\\
&  =\lim_{n\rightarrow\infty}\int h\left(  \bar{Y}\left(  n\right)  \right)
d\bar{Y}\left(  n\right) \nonumber\\
&  =\int h\left(  \bar{Y}\right)  d\bar{Y} \label{relation}%
\end{align}
where $h:V\oplus V\rightarrow L\left(  V\oplus V,W\right)  $ is the
Lip-$\left(  \gamma-1\right)  $ one-form defined by $h\left(  x,y\right)
\left(  v_{1},v_{2}\right)  =g\left(  x\right)  \left(  f\left(  y\right)
\left(  v_{2}\right)  \right)  $ with the auxiliary terms%
\[
h^{j}:V\oplus V\rightarrow L\left(  \left(  V\oplus V\right)  ^{\otimes
j},L\left(  V\oplus V,W\right)  \right)  ,j=1,....,\lfloor\gamma\rfloor-1
\]
given by
\begin{align}
&  h^{j}\left(  x,y\right)  \left(  \left(  v_{1}^{1},v_{2}^{1}\right)
\otimes...\otimes\left(  v_{1}^{j},v_{2}^{j}\right)  \right)
\label{auxilliary terms}\\
&  =\sum_{k=0}^{j}g^{k}\left(  x\right)  \left(  v_{1}^{1}\otimes...\otimes
v_{1}^{j}\right)  f^{j-k}\left(  y\right)  \left(  v_{2}^{k+1}\otimes
...\otimes v_{2}^{j}\right)  .\nonumber
\end{align}
By definition $\int h\left(  \bar{Y}\right)  d\bar{Y}$ is the unique $p$-rough
path associated with the almost $p-$rough path
\[
Q_{s,t}^{n}=\sum_{k_{1},...,k_{n}=1}^{\lfloor p\rfloor}h^{k_{1}-1}\left(
\bar{Y}_{s}\right)  \otimes...\otimes h^{k_{n}-1}\left(  \bar{Y}_{s}\right)
\sum_{\pi\in OS\left(  k_{1},....,k_{n}\right)  }\pi^{-1}\bar{Y}_{s,t}%
^{k_{1}+...+k_{n}};
\]
using (\ref{auxilliary terms}) and the fact that $\pi_{1}\bar{Y}=Z$ and
$\pi_{2}\bar{Y}=Z$ it is easy to see that this coincides with an almost
$p-$rough path associated with $\int\left(  gf\right)  \left(  Z\right)  dZ$.
The result follows from the relation (\ref{relation}).
\end{proof}

\begin{proposition}
\label{rde equivalence}Let $V$ and $W$ be finite dimensional normed vector
spaces and suppose that $f:V\oplus W\rightarrow L\left(  V,W\right)  $ is a
Lip$-\gamma$ one form for $\gamma>p.$ Let $X$ be a geometric $p$-rough path on
$V$ and $Y$ the geometric $p-$rough path obtained by solving the RDE%
\begin{equation}
dY_{t}=f\left(  X_{t},Y_{t}\right)  dX_{t},\text{ \ }Y\left(  0\right)
=Y_{0},\text{ }X\left(  0\right)  =X_{0}. \label{RDE}%
\end{equation}
That is to say $Z=\left(  X,Y\right)  $ $\in G\Omega_{p}\left(  V\oplus
W\right)  $ is such that $Z\left(  0\right)  =$ $\left(  X_{0},Y_{0}\right)
,$ $\pi_{V}Z=X$ and
\begin{equation}
Z=\int\widetilde{f}\left(  Z\right)  dZ, \label{fixed point}%
\end{equation}
where $\widetilde{f}\left(  v_{1},w_{1}\right)  \left(  v_{2},w_{2}\right)
=\left(  v_{2},f\left(  v_{1,}w_{1}\right)  v_{2}\right)  $ is a Lip-$\gamma$
one form on $V\oplus W$. \ For every $v\in V\oplus W$ let $g\left(
v,\cdot\right)  $ be the linear map from $T_{\pi_{V}(v)}V\cong V$ to the space
of Lip-$\gamma$ vector fields on $W$ obtained from $f$ via $g\left(
v,\cdot\right)  :\tilde{v}\longmapsto f\left(  \cdot\right)  \tilde{v}$ \ and
let $\Gamma$ denote the connection obtained from $g$ as in $\left(
\ref{connection-map}\right)  .$ Then, for any Banach space $E$ and any
Lip$-\gamma$ one form $\alpha:V\oplus W\rightarrow L\left(  V\oplus
W,E\right)  $ we have that
\begin{equation}
\int\alpha\left(  Z\right)  dZ=\int\alpha^{\Gamma}\left(  Z\right)  dZ.
\label{one form condition}%
\end{equation}
Conversely, if (\ref{one form condition}) holds for any such one form then $Z
$ is a solution to the classical RDE (\ref{RDE}) in the sense described. \ 
\end{proposition}

\begin{proof}
We first observe that because of the identification $T_{v}V\cong V$ the
connection $\Gamma$ has a simplified form, in fact for all $v_{1},v_{2}\in V $
and $w_{1},w_{2}\in W$ we have
\[
\Gamma_{\left(  v_{1},w_{1}\right)  }\left(  v_{2}\right)  =\left(
v_{2},f\left(  w_{1}\right)  v_{2}\right)  =\widetilde{f}\left(  v_{1}%
,w_{1}\right)  \left(  v_{2},w_{2}\right)  .
\]
Suppose (\ref{one form condition}) is verified then we can show that we have a
classical RDE solution by taking $E=V\oplus W$ and defining $\alpha\left(
v_{1},w_{1}\right)  \left(  v_{2},w_{2}\right)  =\left(  v_{2},w_{2}\right)
$, the identity one form. A simple computation then gives that $\int%
\alpha\left(  Z\right)  dZ=Z$ , hence since $\alpha^{\Gamma}\equiv
\widetilde{f}$ (\ref{fixed point}) is verified and $\pi_{V}Z=X$. \ 

To prove the statement in the other direction we suppose that $Z$ is a
solution to the RDE (\ref{RDE}), then for any Lip$-\gamma$ one form
$\alpha:V\oplus W\rightarrow L\left(  V\oplus W,E\right)  $ we have that
\begin{align*}
\alpha^{\Gamma}\left(  v_{1},w_{1}\right)  \left(  v_{2},w_{2}\right)   &
=\alpha\left(  v_{1},w_{1}\right)  \circ\Gamma_{\left(  v_{1},w_{1}\right)
}\circ\left(  \pi_{V}\right)  _{\ast}\left(  v_{2},w_{2}\right) \\
&  =\alpha\left(  v_{1},w_{1}\right)  \circ\Gamma_{\left(  v_{1},w_{1}\right)
}\left(  v_{2}\right) \\
&  =\alpha\left(  v_{1},w_{1}\right)  \widetilde{f}\left(  v_{1},w_{1}\right)
\left(  v_{2},w_{2}\right) \\
&  =\left(  \widehat{\alpha\widetilde{f}}\right)  \left(  v_{1},w_{1}\right)
\left(  v_{2},w_{2}\right)  .
\end{align*}
Using Lemma \ref{product one forms} with $V\oplus W$ for $V$ and $E$ for $W$
we then have that
\[
\int\alpha^{\Gamma}\left(  Z\right)  dZ=\int\widehat{\alpha\widetilde{f}%
}\left(  Z\right)  dZ=\int\alpha\left(  Z\right)  dZ.
\]

\end{proof}

Note that in the preceding proposition the one form $f$ had a dependence on
the position of the signal $X.$ The usual existence theorems for the solution
of rough differential equations (see e.g. Theorem
\ref{universal limit theorem}) assume that $f:V\rightarrow L\left(
V,W\right)  ,$ i.e. there is no dependence on the position of the signal $X$
in the definition of the one form $f$. In the following lemma we consider an
augmented RDE that incorporates the signal dependence. This allows us to use
the classical RDE theory to derive the existence of solutions in the sense of
$\left(  \ref{fixed point}\right)  .$ Once this is achieved Proposition
\ref{rde equivalence} immediately gives the existence of solutions of rough
differential equations in the manifold sense, i.e. satisfying the fixed point
$\left(  \ref{one form condition}\right)  $. In short, if the manifold is a
finite dimensional normed vector space sufficiently regular rough differential
equations have a rough path solution.

\begin{lemma}
\label{existence augmented rde}Let $1\leq p<\gamma,$ $V$ and $W$ be Banach
spaces and $f:V\oplus W\rightarrow L\left(  V,W\right)  $ be a Lip-$\gamma$
one form. If $X$ is a classical $p$ -rough path on $V$ the rough differential
equation
\begin{equation}
dY_{t}=f\left(  X_{t},Y_{t}\right)  dX_{t},\text{ \ }Y\left(  0\right)
=Y_{0},\text{ }X\left(  0\right)  =X_{0}.
\end{equation}
has a solution in the sense of $\left(  \ref{fixed point}\right)  .$
\end{lemma}

\begin{proof}
Consider the classical rough differential equation%
\begin{equation}
d\hat{Y}_{t}=\hat{f}\left(  \hat{Y}_{t}\right)  dX_{t},\text{ }\hat{Y}%
_{0}=\left(  X_{0},Y_{0}\right)  , \label{augmented rde}%
\end{equation}
where $\hat{f}:V\oplus W\rightarrow L\left(  V,V\oplus W\right)  $ is defined
by $\hat{f}\left(  v_{1},w_{1}\right)  \left(  v_{2}\right)  =\left(
v_{2},f\left(  v_{1},w_{1}\right)  \left(  v_{2}\right)  \right)  .$ From
Theorem \ref{universal limit theorem} equation $\left(  \ref{augmented rde}%
\right)  $ has a unique solution $\hat{Z}$ in $G\Omega_{p}\left(  V\oplus
W\oplus V\right)  .$ We need to demonstrate that
\[
Z:=\pi_{V\oplus W}\hat{Z}\in G\Omega_{p}\left(  V\oplus W\right)
\]
satisfies $\left(  \ref{fixed point}\right)  $, $\pi_{V}Z=X$ and $Z\left(
0\right)  =\left(  X_{0},Y_{0}\right)  $. To achieve this we first note in
passing that the canonical projection map
\[
\pi_{V\oplus W}:\left(  G\Omega_{p}\left(  V\oplus W\oplus V\right)
,d_{p}^{V\oplus W\oplus V}\right)  \rightarrow\left(  G\Omega_{p}\left(
V\oplus W\right)  ,d_{p}^{V\oplus W}\right)
\]
is continuous (in fact $d_{p}^{V\oplus W}\left(  \pi_{V\oplus W}X_{1}%
,\pi_{V\oplus W}X_{2}\right)  \leq d_{p}^{V\oplus W\oplus V}\left(
X_{1},X_{2}\right)  ).$

We now suppose that $\left(  x_{n},y_{n},\tilde{x}_{n}\right)  $ is a sequence
of bounded variation paths in $V\oplus W\oplus V$ whose lift $S_{\left\lfloor
p\right\rfloor }\left(  \left(  x_{n},y_{n},\tilde{x}_{n}\right)  \right)  $
converges to the solution $\widetilde{Z}$ of (\ref{augmented rde}) in
$d_{p}^{V\oplus W\oplus V}$. \ Using the continuity of $\pi_{V\oplus W}$ we
have
\[
S_{\left\lfloor p\right\rfloor }\left(  \left(  x_{n},y_{n}\right)  \right)
=\pi_{V\oplus W}S_{\left\lfloor p\right\rfloor }\left(  \left(  x_{n}%
,y_{n},\tilde{x}_{n}\right)  \right)  \rightarrow\pi_{V\oplus W}Z=Z
\]
in $d_{p}^{V\oplus W}.$ Thus, using the continuity of both $\pi_{V\oplus W}$
and the rough integral together we have that
\begin{align*}
\int\tilde{f}\left(  Z\right)  dZ  &  =\lim_{n\rightarrow\infty}\int\tilde
{f}\left(  S_{\left\lfloor p\right\rfloor }\left(  \left(  x_{n},y_{n}\right)
\right)  \right)  dS_{\left\lfloor p\right\rfloor }\left(  \left(  x_{n}%
,y_{n}\right)  \right) \\
&  =\lim_{n\rightarrow\infty}S_{\left\lfloor p\right\rfloor }\left(
\int\tilde{f}\left(  x_{n},y_{n}\right)  d\left(  x_{n},y_{n}\right)  \right)
\\
&  =\lim_{n\rightarrow\infty}S_{\left\lfloor p\right\rfloor }\left(
x_{n},\int f\left(  x_{n},y_{n},\right)  dx_{n}\right) \\
&  =\lim_{n\rightarrow\infty}\pi_{V\oplus W}\left[  S_{\left\lfloor
p\right\rfloor }\left(  x_{n},\int f\left(  x_{n},y_{n},\right)  dx_{n}%
,\tilde{x}_{n}\right)  \right] \\
&  =\lim_{n\rightarrow\infty}\pi_{V\oplus W}\left[  \int\hat{f}\left(
S_{\left\lfloor p\right\rfloor }\left(  \left(  x_{n},y_{n},\tilde{x}%
_{n}\right)  \right)  \right)  dS_{\left\lfloor p\right\rfloor }\left(
\left(  x_{n},y_{n},\tilde{x}_{n}\right)  \right)  \right] \\
&  =\pi_{V\oplus W}\int\hat{f}\left(  \hat{Z}\right)  d\hat{Z}\\
&  =\pi_{V\oplus W}\hat{Z}\\
&  =Z.
\end{align*}
We finish by noting the trivial fact that $Z\left(  0\right)  =\left(
X_{0},Y_{0}\right)  $ and remarking that $\pi_{V}Z=\pi_{V}\pi_{V\oplus W}%
\hat{Z}=X$ can again be verified by taking an approximating sequence and using
continuity of the projections.
\end{proof}

Let $U\subseteq N,$ $V\subseteq M$ be open, $\xi_{1}:U\rightarrow%
%TCIMACRO{\U{211d} }%
%BeginExpansion
\mathbb{R}
%EndExpansion
^{d_{1}}$ and $\xi_{2}:V\rightarrow%
%TCIMACRO{\U{211d} }%
%BeginExpansion
\mathbb{R}
%EndExpansion
^{d_{2}}$ be Lip-$\gamma$ diffeomorphisms onto open subsets of $%
%TCIMACRO{\U{211d} }%
%BeginExpansion
\mathbb{R}
%EndExpansion
^{d_{1}}$, $%
%TCIMACRO{\U{211d} }%
%BeginExpansion
\mathbb{R}
%EndExpansion
^{d_{2}}$ respectively and
\[
\xi=(\xi_{1},\xi_{2}):U\times V\rightarrow\xi_{1}(N)\times\xi_{2}(M).
\]
Let $\pi_{\xi}$ denote the natural projection from $\xi_{1}(U)\times\xi
_{2}(V)$ onto $\xi_{1}(U).$ Note that for the we have a commutation relation
between the maps $\xi$ and the projections given by%
\[
\pi_{\xi}\circ\xi=\xi_{1}\circ\pi
\]

\begin{definition}
Let $\xi$ be as in the previous paragraph. We define the pushforward of a
connection map $\Gamma$ on $\xi\left(  U\times V\right)  $ by letting%
\begin{equation}
(\xi_{\ast}\Gamma)(\xi(z),\left(  \xi_{1\ast}\right)  v_{z})=\xi_{\ast}%
(v_{z},g(z)v_{z}) \label{pushforward connection}%
\end{equation}
for $z\in U\times V$ and $v_{z}\in T_{\pi(z)}N.$
\end{definition}

Since $\xi$ and $\xi_{1}$ are both invertible Lip-$\gamma$ maps and $\left(
\xi_{1}\right)  _{\ast}(x):$ $T_{x}N$ $\rightarrow T_{\xi_{1}(x)}%
%TCIMACRO{\U{211d} }%
%BeginExpansion
\mathbb{R}
%EndExpansion
^{d_{1}}$ is a linear isomorphism for each $x\in U$ this gives a well defined
map
\[
(\xi_{\ast}\Gamma)(y,\cdot):T_{\pi_{\xi}(y)}%
%TCIMACRO{\U{211d} }%
%BeginExpansion
\mathbb{R}
%EndExpansion
^{d_{1}}\rightarrow T_{y}\left(
%TCIMACRO{\U{211d} }%
%BeginExpansion
\mathbb{R}
%EndExpansion
^{d_{1}}\times%
%TCIMACRO{\U{211d} }%
%BeginExpansion
\mathbb{R}
%EndExpansion
^{d_{2}}\right)
\]
for all $y\in\xi(U\times V).$

In the following we will take the map $\xi$ to be a chart map $\left(
\phi,\psi\right)  $ on the product manifold $N\times M,$ where $\left(
\phi,U\right)  ,$ $\left(  \psi,V\right)  $ are two charts on $N$ and $M$
respectively. The following proposition demonstrates how we (provided the path
is supported inside the chart set) can use the pushforward of the connection
to construct a fixed point on $\xi(U\times V)\subseteq%
%TCIMACRO{\U{211d} }%
%BeginExpansion
\mathbb{R}
%EndExpansion
^{d_{1}}\times%
%TCIMACRO{\U{211d} }%
%BeginExpansion
\mathbb{R}
%EndExpansion
^{d_{2}}$ that is equivalent to $\left(  \ref{solution_def}\right)  ,$ the
fixed point which characterises the solution of a RDE on a manifold.

\begin{proposition}
\label{fixed points proposition}Let $\left(  \phi,U\right)  ,$ $\left(
\psi,V\right)  $ be two charts on $N$ and $M$ respectively. Let $Z$ be a rough
path on $E$ with supp$\left(  Z\right)  \subseteq U$ $\times$ $V$ and
$\xi=\left(  \phi,\psi\right)  $. For any $\gamma\geq1$ and all one forms
$\alpha$ defined on $U\times V:$

\begin{enumerate}
\item $\alpha$ Lip-$\gamma$ implies $\alpha^{\Gamma}$ Lip-$\gamma$ if and only
if for all one forms $\beta$ on $\xi\left(  U\times V\right)  $ the one form
$\beta$ Lip-$\gamma$ implies $\beta\circ(\xi_{\ast}\Gamma)\circ\left(
\pi_{\xi}\right)  _{\ast}$ Lip-$\gamma.$

\item The rough path $Z$ has a fixed point
\begin{equation}
Z(\alpha)=Z(\alpha^{\Gamma}) \label{fp-id1}%
\end{equation}
for all Lip-$\gamma$ one forms $\alpha\,,\alpha^{\Gamma}$ on $U\times V$ if
and only if
\begin{equation}
\xi_{\ast}Z(\beta)=\xi_{\ast}Z\left(  \beta\circ(\xi_{\ast}\Gamma)\circ\left(
\pi_{\xi}\right)  _{\ast}\right)  \label{fp-id2}%
\end{equation}
for all Lip-$\gamma$ one forms $\beta,\beta\circ(\xi_{\ast}\Gamma)\circ\left(
\pi_{\xi}\right)  _{\ast}$ defined on $\xi\left(  U\times V\right)  .$
\end{enumerate}
\end{proposition}

\begin{proof}
First note that if a one form $\zeta$ is Lip-$\gamma$ on $U\times$ $V$ then by
definition $\left(  \xi^{-1}\right)  ^{\ast}\zeta$ is Lip-$\gamma$ on $U\times
V.$ Conversely if $\zeta$ is a Lip-$\gamma$ one form $\xi\left(  U\times
V\right)  $ the pull back $\xi^{\ast}\zeta$ is by Lemma \ref{pullback-oneform}
Lip-$\gamma$ on $U\times$ $V.$ By definition
\[
\xi_{\ast}Z(\beta)=Z\left(  \xi^{\ast}\beta\right)  .
\]
We first show that%
\begin{equation}
\xi^{\ast}(\beta\circ(\xi_{\ast}\Gamma)\circ\left(  \pi_{\xi}\right)  _{\ast
})(z,v_{z})=\xi^{\ast}\beta\circ\Gamma\circ\pi_{\ast}(z,v_{z})=(\xi^{\ast
}\beta)^{\Gamma} \label{fp-id3}%
\end{equation}
for all $z\in E,$ $v_{z}\in T_{z}E.$ We have%
\[
\xi^{\ast}(\beta\circ(\xi_{\ast}\Gamma)\circ\left(  \pi_{\xi}\right)  _{\ast
})(z,v_{z})=(\beta\circ(\xi_{\ast}\Gamma)\circ\left(  \pi_{\xi}\right)
_{\ast})(\xi(z),\xi_{\ast}v_{z})
\]%
\[
=\left[  \beta\circ(\xi_{\ast}\Gamma)\right]  (\xi\left(  z\right)  ,(\pi
_{\xi})_{\ast}\xi_{\ast}v_{z})=\left[  \beta\circ(\xi_{\ast}\Gamma)\right]
(\xi\left(  z\right)  ,\phi_{\ast}\pi_{\ast}v_{z})
\]%
\begin{align}
&  =\beta(\xi\left(  z\right)  ,\xi_{\ast}\left[  \pi_{\ast}v_{z}%
,g(z)\pi_{\ast}v_{z}\right]  )=\xi^{\ast}\beta(z,\left(  (\pi_{\ast}%
v_{z},g(z)\pi_{\ast}v_{z})\right) \nonumber\\
&  =\xi^{\ast}\beta\circ\Gamma\circ\pi_{\ast}(z,v_{z}).
\end{align}

Let $\alpha$ be a Lip-$\gamma$ one form on $U\times V$ and setting
$\beta=\left(  \xi^{-1}\right)  ^{\ast}\alpha$ in $\left(  \ref{fp-id3}%
\right)  $ we see that
\[
\xi^{\ast}(\beta\circ(\xi_{\ast}\Gamma)\circ\left(  \pi_{\xi}\right)  _{\ast
})(z,v_{z})=\alpha^{\Gamma}.
\]
Therefore we therefore deduce the reverse direction of the proposition, namely
that if $\beta$ Lip-$\gamma$ on $\xi\left(  U\times V\right)  $ implies
$\beta\circ(\xi_{\ast}\Gamma)\circ\left(  \pi_{\xi}\right)  _{\ast}$
Lip-$\gamma$ it follows that $\alpha$ Lip-$\gamma$ $U\times V$ implies
$\alpha^{\Gamma}$ Lip-$\gamma$ and the fixed point $\left(  \ref{fp-id2}%
\right)  $ implies $\left(  \ref{fp-id1}\right)  .$ To deduce the forward
direction note that $\xi$ is invertible on $U\times V$ and therefore $\left(
\xi^{-1}\right)  _{\ast}\xi_{\ast}\Gamma=\Gamma$ on $U\times V$ and that
$\left(  \xi^{-1}\right)  _{\ast}\xi_{\ast}Z=Z$ (use that support of $Z$ is
contained inside the chart and argue as in Remark \ref{one-forms-open} to see this).
\end{proof}

We are finally ready to prove the existence of solutions to rough differential
equations in the manifold sense. As usual we assume $g$ to have Lip-$\gamma$ regularity.

\begin{definition}
For all $x\in N$ let $g\left(  x,\cdot\right)  $ be a linear map from $T_{x}N$
to $\tau\left(  M\right)  $, the space of vector fields on $M.$ We say the map
$g$ is Lip-$\gamma$ with Lipschitz constant $C$ if there exists $C>0$ such
that for any charts $(\phi,U)$ on $N$ and $\left(  \psi,V\right)  $ on $M$
respectively the function $g_{\phi,\psi}:\left(  \phi\,,\psi\right)  \left(
U\times V\right)  \rightarrow L\left(
%TCIMACRO{\U{211d} }%
%BeginExpansion
\mathbb{R}
%EndExpansion
^{d_{1}},%
%TCIMACRO{\U{211d} }%
%BeginExpansion
\mathbb{R}
%EndExpansion
^{d_{2}}\right)  $ defined by%
\[
g_{\phi,\psi}(z,\nu)=\psi_{\ast}\left[  g\left(  \left(  \phi,\psi\right)
^{-1}\left(  z\right)  \right)  \left(  \phi_{\ast}^{-1}\nu\right)  \right]
,\text{ }z\in\left(  \phi\,,\psi\right)  \left(  U\times V\right)  ,~\nu\in%
%TCIMACRO{\U{211d} }%
%BeginExpansion
\mathbb{R}
%EndExpansion
^{d_{1}}%
\]
is Lip-$\gamma$ with Lipschitz constant at most $C$.
\end{definition}

To prove existence of a solution $Z$ to a rough differential equation on a
manifold we will - all technicalities aside - proceed in four simple steps. We
first push the equation forward onto a chart and identify the pushed forward
equation and noise with the analogous classical rough differential equation.
We deduce using the theory of classical rough differential equations the
existence of a solution $\widehat{Z}$ to this equation for some (uniform) time
interval $\left[  0,t_{0}\right]  .$ The time $t_{0}$ is chosen such that the
solution remains inside the image of the chart on the product manifold $E$ and
is independent of the starting point of the equation. We identify the
classical rough path $\widehat{Z}$ with a rough path in the manifold sense,
demonstrate it satisfies a fixed point identity and eventually pull it back
onto the manifold to obtain a solution to $\left(  \ref{RDE2}\right)  $ over
the time interval $\left[  0,t_{0}\right]  $ . Finally, repeating the process
finitely many times and concatenating the solutions over the subintervals we
obtain a solution over $\left[  0,T\right]  .$ In the proof of the following
theorem we will arguing as in Remark \ref{one-forms-open} frequently (and
slightly abusing notation) consider the push forwards of rough paths under
maps that are only defined on open subsets containing the support of the path.

\begin{theorem}
Let $1\leq p<\gamma\leq\gamma_{0},$ $M$ and $N$ are Lip-$\gamma_{0}$ manifolds
and $X$ a geometric $p$ -rough path on $N$ with starting point $x_{0}.$ For
each $x\in N$ let $g\left(  x,\cdot\right)  $ be a linear map from $T_{x}N$ to
$\tau\left(  M\right)  $ and suppose that $g$ is Lip-$\gamma.$ Then the
equation
\begin{equation}
dY_{t}=g\left(  X_{t},Y_{t}\right)  dX_{t},\text{ }Y_{0}=y_{0}
\label{rde-thm-id}%
\end{equation}
has a solution in the sense of Definition $\ref{solution_def}$.
\end{theorem}

\begin{proof}
By definition of a Lip-$\gamma$ atlas given $\left(  x_{0},y_{0}\right)  \in
N\times M$ there is a chart $\left(  (\phi,\psi),\tilde{U}\times\tilde
{V}\right)  $ such that $B\left(  (\phi,\psi)\left(  x_{0},y_{0}\right)
,\delta\right)  \subseteq(\phi,\psi)(\tilde{U}\times\tilde{V}).$ Let
$\xi=(\phi,\psi)$, $U=\tilde{U}^{\delta/2}$ and $V=\tilde{V}^{\delta/2}$
(recall $\left(  \ref{delta-sets}\right)  $ for the definition of these
sets)$.$ We will prove the existence of a solution to the RDE up to some
$t_{0}>0$ that only depends on $\left\Vert g\right\Vert _{Lip-\gamma}$, $p,$
$\delta,$ $L$ and the control of the signal $\omega.$ In particular $t_{0}$ is
going to be independent of the starting point $\left(  x_{0},y_{0}\right)  .$
For $z\in(\phi,\psi)(U\times V),$ $v\in%
%TCIMACRO{\U{211d} }%
%BeginExpansion
\mathbb{R}
%EndExpansion
^{d_{1}}$ let $h\left(  z,\nu\right)  =\psi_{\ast}g(\xi^{-1}(z),\phi_{\ast
}^{-1}v).$ Then $h$ is by assumption Lip-$\gamma$ with constant $C$ on
$(\phi,\psi)(U\times V)$ and taking values in the vector fields. Arguing as in
Remark \ref{one-forms-open} we may extend $h$ to a Lip-$\gamma$ one form on
$R^{d_{1}}\times R^{d_{2}}$ with Lipschitz constant depending only on $C$ and
the Lip-$\gamma_{0}$ atlas of $N\times M.$ Let $\widehat{X}$ denote the
classical rough path corresponding to $\phi_{\ast}X$ obtained by Lemma
\ref{bijection-bspace}. By Lemma \ref{existence augmented rde} there exists a
classical (based) rough path $\widehat{Z}$ on $%
%TCIMACRO{\U{211d} }%
%BeginExpansion
\mathbb{R}
%EndExpansion
^{d_{1}}\times%
%TCIMACRO{\U{211d} }%
%BeginExpansion
\mathbb{R}
%EndExpansion
^{d_{2}}$ solving
\[
d\widehat{Z}_{t}=h\left(  \widehat{Z},\widehat{X}\right)  d\widehat{X}_{t}%
\]
with initial conditions $(\phi,\psi)(x_{0},y_{0})$, i.e. if $\widetilde{h}$ is
the one form defined by%
\[
\widetilde{h}\left(  v_{1},w_{2}\right)  \left(  v_{2},w_{2}\right)  =\left(
v_{2},h\left(  v_{1},w_{1}\right)  v_{2}\right)  ;\text{ }v_{1,}v_{2}\in%
%TCIMACRO{\U{211d} }%
%BeginExpansion
\mathbb{R}
%EndExpansion
^{d_{1}};w_{1},w_{2}\in%
%TCIMACRO{\U{211d} }%
%BeginExpansion
\mathbb{R}
%EndExpansion
^{d_{2}}%
\]
the path $\widehat{Z}$ satisfies $\widehat{Z}\left(  0\right)  =$ $\left(
x_{0},y_{0}\right)  ,$ $\pi_{\xi}\widehat{Z}=\widehat{X}$ and
\begin{equation}
\widehat{Z}=\int\widetilde{h}\left(  \widehat{Z}\right)  d\widehat{Z}.
\end{equation}
Note that by the universal limit theorem (Theorem
\ref{universal limit theorem}) for classical rough paths, there exists a time
$t_{0}$ depending on $\delta,$ the control of the signal and $\left\Vert
\widetilde{h}\right\Vert _{Lip},$ but independent of the starting point
$(\phi,\psi)\left(  x_{0},y_{0}\right)  $ such that supp$\left(
\widehat{Z}\right)  \subseteq B\left(  (\phi,\psi)\left(  x_{0},y_{0}\right)
,\delta\right)  \subseteq(\phi,\psi)(U\times V).$ We may readily identify
$\widetilde{h}$ as $\widetilde{h}=\xi_{\ast}\Gamma$ (where $\Gamma$ is the
connection map corresponding to $g$ and $\xi_{\ast}\Gamma$ was defined in
$\left(  \ref{pushforward connection}\right)  $). We deduce by Proposition
\ref{rde equivalence} that if $\beta$ is a Lip-$\gamma$ one form on
$(\phi,\psi)(U\times V),$ $\beta\circ(\xi_{\ast}\Gamma)\circ\left(  \pi_{\xi
}\right)  _{\ast}$ is Lip-$\gamma$ and%
\[
\int\beta\left(  \widehat{Z}\right)  d\widehat{Z}=\int\beta\circ(\xi_{\ast
}\Gamma)\circ\left(  \pi_{\xi}\right)  _{\ast}\left(  \widehat{Z}\right)
d\widehat{Z}.
\]
Setting $\widetilde{Z}_{0}=(\phi,\psi)(x_{0},y_{0})$ and letting
$\widetilde{Z}(\alpha)=\int ad\widehat{Z}$ we obtain a rough path in the
manifold sense over the interval $\left[  0,t_{0}\right]  $ with
\begin{equation}
supp\left(  \widetilde{Z}\right)  \subseteq(\phi,\psi)(U\times V).
\label{lt-id1}%
\end{equation}
From the properties of $\widehat{Z}$ it is clear that $\widetilde{Z}$
satisfies
\begin{equation}
\left(  \pi_{\xi}\right)  _{\ast}\widetilde{Z}\sim\phi_{\ast}X \label{lt-id2}%
\end{equation}
and
\begin{equation}
\widetilde{Z}(\beta)=\widetilde{Z}\left(  \beta\circ(\xi_{\ast}\Gamma
)\circ\left(  \pi_{\xi}\right)  _{\ast}\right)  \label{lt-id3}%
\end{equation}
for any Lip-$\gamma$ one form $\beta$ on $(\phi,\psi)(U\times V)$ . As
$\widetilde{Z}$ satisfies $\left(  \ref{lt-id1}\right)  $ and $\left(
\ref{lt-id3}\right)  $ we may apply Proposition \ref{fixed points proposition}
and we get that $Z=\left(  \xi^{-1}\right)  _{\ast}\widetilde{Z}$ satisfies%
\[
Z(\alpha)=Z(\alpha^{\Gamma})
\]
for any Lip-$\gamma$ one form $\alpha$. It follows from $\left(
\ref{lt-id2}\right)  $ that $\pi_{\ast}\left(  \xi^{-1}\right)  _{\ast
}\widetilde{Z}\sim X$ (recall that $\pi$ denotes the projection from $E$ onto
$N$) and we deduce that $Z$ solves the rough differential equation
(\ref{rde-thm-id}) over the time interval $\left[  0,t_{0}\right]  .$ Finally,
by Lemma \ref{consistency of localisations} the path $Z$ has an endpoint.
Hence, we may repeat the entire process from time $t_{0}$ finitely many times
on different charts. We then concatenate the solution paths over the
subintervals and thus obtain a solution over the full time interval $\left[
0,T\right]  .$
\end{proof}

\section{Perspectives}

The development of RDEs presented here allows us to study stochastic processes
on manifolds interpreted as rough paths in the sense described. The authors
would like to finish this article by drawing attention to certain deficiencies
and loose ends in our development. Firstly, our definition of rough paths is
too focused on the finite dimensional case. Whilst this covers a satisfactory
number of examples, it would nevertheless be interesting to re-work things so
that one-forms and Lip-$\gamma$ functions are defined on closed sets (\`{a} la
Stein) \ rather than open sets; such an analysis will probably be most
accurately formulated using the language of jets. Finally, our proof that a
rough path on $M$ is parametrised by a compact interval and had compact trace
is annoyingly complex; much effort could be conserved by building in the
assumption of the compactness of the support to definition of a rough path on
$M.$

\end{document}